\documentclass[a4paper,11pt]{article}
\usepackage[utf8]{inputenc}
\usepackage[hmargin={25mm,25mm},vmargin={25mm,30mm}]{geometry}
\usepackage{amsmath,amsfonts,amsthm,amssymb,mathtools}
\usepackage{graphicx}
\usepackage{subcaption}
\usepackage{float}
\usepackage[ocgcolorlinks,allcolors={blue}]{hyperref}
\usepackage{authblk}
\usepackage{multirow}
\usepackage{booktabs}
\usepackage{placeins}

\theoremstyle{definition}
\newtheorem{theorem}{Theorem}
\newtheorem{proposition}[theorem]{Proposition}
\newtheorem{lemma}[theorem]{Lemma}
\newtheorem{definition}[theorem]{Definition}

\theoremstyle{remark}
\newtheorem{remark}[theorem]{Remark}

\usepackage{newtxtext,newtxmath}

\newcommand{\bu}{\boldsymbol{u}}
\newcommand{\bv}{\boldsymbol{v}}
\newcommand{\bw}{\boldsymbol{w}}
\newcommand{\ff}{\boldsymbol{f}}
\newcommand{\bg}{\boldsymbol{g}}
\newcommand{\bn}{\boldsymbol{n}}
\newcommand{\bsigma}{\boldsymbol{\sigma}}
\newcommand{\btau}{\boldsymbol{\tau}}
\newcommand{\beps}{\boldsymbol{\varepsilon}}
\def\b{\boldsymbol}
\newcommand{\R}{\mathbb{R}}
\newcommand{\GRAD}{\b\nabla}
\newcommand{\divs}{\nabla{\cdot}}
\newcommand{\DIV}{\b\nabla{\cdot}}
\DeclareMathOperator{\dev}{dev}
\DeclareMathOperator{\tr}{tr}
\newcommand{\dd}{\,\mathrm{d}}
\newcommand{\bz}{\boldsymbol{z}}
\newcommand{\by}{\boldsymbol{y}}
\newcommand{\eps}{\beps}
\newcommand{\Sym}{\mathbb{S}}
\newcommand{\dual}[2]{\left\langle #1,#2\right\rangle}
\newcommand{\email}[1]{\href{mailto:#1}{#1}}

\title{Efficiency and robustness of finite element discretizations and linearization strategies for the primal and dual $p$-Stokes problem}
\author[2]{Andrea Munno\thanks{\protect\hypertarget{note:equal}{}These authors contributed equally and are listed in alphabetical order.}\textsuperscript{,}\thanks{\email{andrea.munno@epfl.ch}}}
\author[3]{Jacopo Senoner\textsuperscript{\protect\hypersetup{hidelinks}\protect\hyperlink{note:equal}{*}}\textsuperscript{,}\thanks{\email{jacopo.senoner@epfl.ch}}}
\author[1]{Michele Botti\thanks{\email{michele.botti@polimi.it}}}

\affil[1]{MOX, Dipartimento di Matematica, Politecnico di Milano,
Piazza Leonardo da Vinci 32 - 20133 Milano, Italy}
\affil[2]{École polytechnique fédérale de Lausanne (EPFL),
Mathematics Section, 1015 Lausanne, Switzerland}
\affil[3]{École polytechnique fédérale de Lausanne (EPFL),
Section of Computer Science, 1015 Lausanne, Switzerland}
\date{}

\begin{document}

\maketitle

\begin{abstract}
We present a systematic numerical investigation of finite element
discretizations and nonlinear solution strategies for incompressible power-law Stokes flows, covering both shear-thinning and shear-thickening regimes. We consider primal and Fenchel--dual formulations and compare the MINI, $\mathbb{P}_2$--$\mathbb{P}_1$, and Kouhia--Stenberg discretizations with Newton and Ka\v{c}anov iterations. The selected discretizations ensure convergence across all regimes and provide a good balance between implementation simplicity and accuracy. We derive pressure-explicit representations of the primal and dual energy functionals and a total-stress formulation that naturally incorporates the pressure through the stress. A set of manufactured benchmarks is used to assess the effects of the nonlinear constitutive law, localized loss of regularity, and extended zero-strain degeneracy. 
%Numerical results show that primal Ka\v{c}anov is more robust for shear-thinning problems, while Newton and dual methods perform better in the shear-thickening regime. 
%The experiments also show that, when the velocity is small compared to the pressure, the velocity error of the considered elements is dominated by a pressure-induced consistency error, which a first test with a pressure-robust load reconstruction largely removes. 
We further clarify the distinct roles of regularization and $r$-continuation: the former involves an accuracy--robustness trade-off, while the latter improves the initial guess without altering the target discrete problem. Obstacle-flow and glaciological tests illustrate the applicability of the framework beyond manufactured solutions.
\end{abstract}

\section{Introduction}

The accurate numerical simulation of non-Newtonian incompressible flows is a critical requirement for several industrial and scientific applications, as these models capture the behavior of complex fluids such as blood, cosmetics, and glacial ice, among many others. This work addresses the numerical approximation of the steady non-linear Stokes equations for non-Newtonian fluids described by a power-law constitutive model, covering both shear-thinning and shear-thickening regimes. A central challenge in this context is the handling of singularities inherent in the viscosity model for shear-thinning fluids.

Numerical methods for non-Newtonian flows in the shear-thinning (pseudo-plastic) regime have a long history, beginning with the seminal work of \cite{BN:1990}, which proposed a Finite Element approximation of a non-Newtonian flow model governed by either the Carreau or the power-law model. Still in the Finite Element framework, sharp error estimates were subsequently established in \cite{Sandri:93, BL:1993, Barrett.Liu:94}. In particular, the pioneering studies \cite{BL:1993} and \cite{Barrett.Liu:94} derived (in some cases, optimal) velocity and pressure error bounds in appropriate quasi-norms for models incorporating Carreau or power-law models. We also refer to, e.g., \cite{Belenki.Berselli.ea:12, Hirn:13, Kreuzer.Suli:16, Kalte.Ruzi:23} for more recent contributions.

In this work, we present a comparative study of primal and dual finite element discretizations, assessing the performance of the MINI, lowest-order Taylor--Hood,
and Kouhia--Stenberg elements. The nonlinear problem is solved using different linearization strategies depending on the exponent $r\in(1,\infty)$ of the power-law model: fixed-point iterations (Picard/Ka\v{c}anov) and the Newton method. To improve robustness in the shear-thinning regime, where the effective viscosity may become strongly degenerate, we investigate two regularization/relaxation techniques. The first adds a small positive shift to the denominator to avoid singular behavior. The second employs a continuation strategy in the power-law parameter: starting from a milder nonlinearity and progressively approaching the target value, using each converged solution as the initial guess for the next step.
Numerical investigations are conducted to assess the convergence rates of considered finite element discretizations, as well as computational efficiency and robustness of the linearization solvers with respect to extreme nonlinearities and singular solutions. To further investigate this last aspect, we have implemented a degenerate test case in which the velocity gradient vanishes on two intersecting lines in the domain, using it as a demanding benchmark for assessing the robustness of the solver in shear-thinning regimes.

Related numerical comparisons are available in the literature, although their scope differs from that of the present work. Schmidt et al.~\cite{SchmidtHumbertSlawig2024} compare Picard and Newton-type nonlinear solvers, together with different step-size strategies, for full-Stokes ice-flow simulations in the shear-thinning regime. Hirn~\cite{Hirn2013} studies the
finite element approximation of singular shear-thinning power-law systems, with particular emphasis on regularization and the numerical stability of Newton's method. More recently, Helanow and
Ahlkrona~\cite{HelanowAhlkrona2025} consider both Taylor--Hood and MINI elements, but their comparison is primarily concerned with the spectral properties of a block preconditioner for regularized power-law ice flow. These works therefore do not simultaneously incorporate all three features considered here: both primal and dual formulations, a focus on inf-sup stable finite element discretizations as an effective compromise between simplicity and performance, and coverage of both the shear-thinning ($r\le2$) and shear-thickening ($r\ge 2$) regimes.

\subsection{Applications}

We present three main applications of the power-law model: hemodynamics, glaciology and impact-protection technologies. We include some remarks to motivate our tests and the development of robust numerical strategies for power-law models, particularly in the presence of singularities.

\paragraph{Hemodynamics}

While the Carreau--Yasuda model is widely adopted to describe blood rheology over broad shear-rate ranges~\cite{Seehanam2023},
the Walburn--Schneck model \cite{WalburnSchneck1976} provides a physiologically motivated alternative that may be especially attractive for patient-specific computational fluid dynamics.
Unlike general rheological models, whose parameters must typically be identified through dedicated rheometric experiments, the Walburn--Schneck formulation directly incorporates standard clinical hematological measurements into the constitutive law.
From a mathematical perspective, this model can be interpreted as a parameter-dependent power-law fluid.
Moreover, it has been observed that rheological parameters are not always reliably identifiable from experimental data \cite{Gallagher2019}.
It should be noted, however, that the power-law model does not reproduce the viscosity plateaus exhibited by blood at very low and very high shear rates, which are more faithfully captured by the Carreau--Yasuda model. Nevertheless, the comparison reported in Appendix~\ref{app:obstacle-blood-comparison} shows that, for the hemodynamically motivated benchmark considered here, this limitation has only a minor effect on the global velocity and pressure fields, despite appreciable local differences in apparent viscosity.

\paragraph{Glaciology}

Glaciology provides an important example of a shear-thinning power-law rheology. In ice-sheet and glacier modelling, the internal deformation of ice is commonly described by Glen's flow law~\cite{Glen1955}. The conventional stress exponent $n_G=3$~\cite[Sec.~3.4.6]{CuffeyPaterson2010} corresponds to $r=1+1/n_G=\frac{4}{3}$ in our notation.
The unregularized effective viscosity is singular at zero strain. Isaac et al.~\cite{IsaacStadlerGhattas2015} used a Newton method with line search for the nonlinear Stokes equations with regularized Glen rheology.
Ice rheology can depend on temperature and anisotropy~\cite{Fan2025FlowLawsIce}. Millstein et al.~\cite{Millstein2022IceViscosity} inferred $n_G=4.1\pm0.4$ ($r\approx\frac{5}{4}$) in fast-flowing extensional regions of Antarctic ice shelves. Recent studies relying on Glen-family rheologies for ice-sheet modeling setups are reported in \cite{Cuzzone2025PolythermalISSM,Leger2025GreenlandPISM}.

\paragraph{Impact-protection technologies}

The shear-thickening regime of the power-law model is naturally related to impact-protection technologies, where materials should remain flexible under ordinary motion but become significantly more resistant under sudden loading. A prominent example is given by shear-thickening fluids, typically dense suspensions of hard particles in a carrier liquid, which have been investigated for protective fabrics 
% and helmets
\cite{Lee2003BallisticSTF,Decker2007StabSTF}.
Reported effective exponents for silica/PEG shear-thickening fluids span approximately \(r \simeq 2.5\text{–}15.9\), covering moderate to very steep local thickening regimes \cite{Manukonda2020,BoscoCaladoMaia2019}. This range supports the \(r\) values selected in our tests.

\subsection{Notation}

Throughout this paper, let $\Omega \subset \mathbb{R}^2$ denote a bounded, connected, polygonal open set with Lipschitz boundary $\partial\Omega$, and let $\mathbf{n}$ denote the outward unit normal vector to $\partial\Omega$.
For ease of exposition, we restrict our presentation to the two-dimensional setting, noting that analogous results can be obtained for the three-dimensional case with suitable modifications. In particular, suitable modifications of the finite element spaces are required, as illustrated by the extension of the two-dimensional Kouhia--Stenberg element~\cite{KouhiaStenberg1995} to three dimensions proposed by~\cite{HuSchedensack2019}.
We assume that the boundary is partitioned into two measurable disjoint subsets $\Gamma_D$ and $\Gamma_N$ such that $\partial\Omega = \bar{\Gamma}_D \cup \bar{\Gamma}_N$ and $\Gamma_D \cap \Gamma_N = \emptyset$. We further assume that $|\Gamma_D| > 0$.
For any two vectors $\mathbf{v}, \mathbf{w} \in \mathbb{R}^2$, we denote their scalar product by $\mathbf{v} \cdot \mathbf{w}$ and the Euclidean norm by $|\mathbf{v}|$.
Similarly, for tensor fields $\btau, \boldsymbol{\eta} \in \mathbb{R}^{2\times 2}$, the Frobenius inner product is defined as $\btau : \boldsymbol{\eta} = \sum_{i,j=1}^{2} \btau_{ij} \boldsymbol{\eta}_{ij}$, and the induced norm is given by $\|\btau\| = (\btau : \btau)^{\frac{1}{2}}$.
Denoting by $\mathbf I$ the identity matrix of $\mathbb{R}^{2\times 2}$, we additionally define the trace and deviatoric operators by
${\rm tr}(\boldsymbol{\tau}) \coloneq \boldsymbol\tau : \mathbf I$ and $\boldsymbol{\rm dev}(\boldsymbol\tau)\coloneq \boldsymbol\tau - 2^{-1}{\rm tr}(\boldsymbol{\tau})\mathbf I$, respectively.
For a scalar $r>1$, we denote its conjugate by
$
r' := \frac{r}{r-1}
$. $A \Subset B$ means that $A$ is compactly contained in $B$, i.e.,
$\overline{A}$ is compact and $\overline{A} \subset B$.

For $s\ge0$ and $q\in[1,\infty]$, we denote by $W^{s,q}(X)$ the Sobolev spaces of order $s$, including fractional orders, with the conventions $W^{0,q}(X)=L^q(X)$ and $W^{s,2}(X)=H^s(X)$.
Given $\Gamma\subset\partial\Omega$, we denote by $W^{1,q}_{0,\Gamma}(\Omega)$ the subspace of $W^{1,q}(\Omega)$ consisting of functions having zero trace on $\Gamma$.
The subscript $\rm{sym}$
denotes a space of symmetric tensor fields.
The symbol $\nabla$ denotes the gradient for scalar functions, while
$\GRAD$, $\b{\epsilon}(\cdot) = 2^{-1}(\GRAD(\cdot) + \GRAD^{\rm T}(\cdot))$, and $\divs$ denote the gradient, the symmetric gradient operator, and the divergence operator for vector fields, respectively, whereas $\DIV$ denotes the vector-valued divergence operator for tensor fields.
Finally, for further use, we adopt the shorthand notation $\mathsf{a} \lesssim \mathsf{b}$ to denote the inequality $\mathsf{a} \leq C \mathsf{b}$, for a generic constant $C>0$ that may depend on the computational domain $\Omega$ and the exponent of the nonlinearity $r$, but is independent of the discretization parameters.

\subsection{Structure of the paper}

The remainder of this paper is organized as follows. Section~\ref{sec:continuouspb} introduces the power-law Stokes model and its primal and Fenchel-dual variational formulations. Section~\ref{sec:discretepb} presents the finite element discretizations and the nonlinear solution strategies considered in this work. The benchmark problems and error measures are described in Section~\ref{sec:experimentalsetting}. Section~\ref{sec:numerical-results} reports and discusses the numerical results. Section~\ref{sec:glaciologytests} illustrates the application of the proposed framework to a glaciological test based on the Haut Glacier d'Arolla geometry. Finally, Section~\ref{sec:conclusions} summarizes the main findings and draws conclusions. The appendices collect complementary material. Appendix~\ref{app:duality-gap} derives the primal and dual energy functionals and proves the primal--dual equality of Theorem~\ref{thm:primal-dual-equality}, together with a numerical illustration of the primal--dual gap. Appendix~\ref{app:proof-equivalence-dual-schemes} contains the proof of Proposition~\ref{prop:dual-kacanov-equivalence} on the equivalence of the two dual Ka\v{c}anov iterations on the Kouhia--Stenberg discretization. Appendix~\ref{app:obstacle-domain-test} reports a flow past a circular obstacle, for which no exact solution is available, to assess the solvers on a more complex geometry, and Appendix~\ref{app:obstacle-blood-comparison} compares, on the same geometry, the power-law model with a Carreau--Yasuda model fitted to the same blood rheology data.

\section{Model problem}\label{sec:continuouspb}

We consider the steady flow of an incompressible non-Newtonian fluid occupying the domain \(\Omega\), subject to a body force \(\mathbf{f} \colon \Omega \to \mathbb{R}^2\) and a boundary traction \(\mathbf{g} \colon \Gamma_N \to \mathbb{R}^2\). A prescribed velocity \(\bu_D \colon \Gamma_D \to \mathbb{R}^2\) is enforced on the Dirichlet boundary.
The flow is governed by the nonlinear Stokes system
\begin{equation}\label{eq:problem}
\begin{cases}
-\DIV \bsigma(\boldsymbol{\varepsilon}(\bu)) + \nabla p = \mathbf{f} & \text{in } \Omega, \\
\nabla \cdot \bu = 0 & \text{in } \Omega, \\
\bsigma(\boldsymbol{\varepsilon}(\bu))\,\mathbf{n} - p\,\mathbf{n} = \mathbf{g} & \text{on } \Gamma_N, \\
\bu = \bu_D & \text{on } \Gamma_D,
\end{cases}
\end{equation}
where \(\bu \colon \Omega \to \mathbb{R}^2\) and \(p \colon \Omega \to \mathbb{R}\) denote the velocity and pressure fields, respectively, and
\(
\boldsymbol{\varepsilon}(\bu)
\)
is referred to as the strain-rate tensor.

\subsection{Constitutive law}

To describe the rheological behavior of the fluid, we first introduce the \emph{Carreau--Yasuda} model \cite{yasuda1981shear}.
The constitutive law is given by

\begin{equation}\label{eq:carreau-yasuda}
\bsigma(\boldsymbol{\varepsilon}(\bu))
:=
\mu \left( \delta^\alpha + \|\boldsymbol{\varepsilon}(\bu)\|^\alpha \right)^{\frac{r-2}{\alpha}} \boldsymbol{\varepsilon}(\bu),
\end{equation}
where \(\mu \colon \Omega \to [\mu_{-},\mu_{+}]\), with $0<\mu_{-}\le \mu_{+}$, is the viscosity coefficient, \(\delta \ge 0\) is a regularization parameter, \(r \in (1,\infty)\) is the flow index, and \(\alpha > 0\) controls the transition between the zero-shear-rate and the power-law regimes.
The constitutive law \eqref{eq:carreau-yasuda} includes several classical models as special cases:
\begin{itemize}
    \item if \(\alpha = 2\), one recovers the standard \emph{Carreau model};
    \item if \(\delta = 0\), the model reduces to the \emph{power-law model}, for which the viscosity depends on the shear rate as \(\|\boldsymbol{\varepsilon}(\bu)\|^{\,r-2}\);
    \item if \(r = 2\), problem~\eqref{eq:problem} reduces to the linear Stokes system for a Newtonian fluid.
\end{itemize}

The majority of non-Newtonian fluids of interest exhibit \emph{shear-thinning} behavior, corresponding to the regime \(1 < r < 2\), in which the effective viscosity decreases as the shear rate increases.
In what follows, we restrict our attention to the power-law case, setting \(\delta = 0\). This case is of particular interest because it leads to a more challenging numerical problem than the regularized Carreau--Yasuda model. Indeed, for shear-thinning fluids (\(r<2\)), the effective viscosity becomes singular as \(\|\boldsymbol{\varepsilon}(\bu)\| \to 0\).

\subsection{Weak formulation}

We introduce the weak formulation of the problem in both primal and dual form, the latter obtained by applying Fenchel duality to the velocity--pressure energy functional. We investigate below whether this formulation offers theoretical or numerical advantages.

\subsubsection{Velocity-pressure formulation}\label{subsubsec:primal-weak-form}

First, we introduce a stable lifting of the Dirichlet datum.
Assuming $\bu_D \in W^{\frac1{r'},r}(\Gamma_D)^2$, owing to the surjectivity of the trace operator, there exists $\bw_1\in W^{1,r}(\Omega)^2$ such that $\bw_1=\bu_D$ on $\Gamma_D$ and
$$
\|\bw_1\|_{W^{1,r}(\Omega)^2}\lesssim \|\bu_D \|_{W^{\frac1{r'},r}(\Gamma_D)^2}.
$$
An explicit construction of $\bw_1$ can be obtained by solving a vector Laplacian problem with non-homogeneous boundary data.
If $\Gamma_D=\partial\Omega$, we assume that $\bu_D$ satisfies the compatibility condition $\int_{\partial\Omega} \bu_D\cdot\bn =0$ and, as a result, we have $\int_{\Omega} \nabla\cdot\bw_1 = \int_{\partial\Omega}\bw_1\cdot\bn =0$, i.e.
$$
\nabla\cdot\bw_1 \in L^{r}_0(\Omega) \coloneqq \left\{q \in L^{r}(\Omega)\ : \ \textstyle\int_\Omega q = 0 \right\}.
$$
Using \cite[Lemma III.3.1]{Galdi:2011}, it is possible to construct~$\bw_2$ in~$W^{1,r}_{0,\partial\Omega}(\Omega)$ satisfying $\nabla\!\cdot \bw_2=\nabla\!\cdot \bw_1$ and
\begin{equation}\label{eq:stab-est}
\|\bw_2\|_{W^{1,r}(\Omega)^2}\lesssim
\|\bw_1 \|_{W^{1,r}(\Omega)^2} \lesssim
\|\bu_D \|_{W^{\frac1{r'},r}(\Gamma_D)^2}.
\end{equation}
In the case $\Gamma_D\neq \partial\Omega$, instead we can use \cite[Theorem 1.4]{Botti-Mascotto:2026} to find
$\bw_2$ in~$W^{1,r}_{0,\Gamma_D}(\Omega)$ such that $\nabla\!\cdot \bw_2=\nabla\!\cdot \bw_1$ and inequality \eqref{eq:stab-est} holds true.
Thus, taking $\bw_D=\bw_1-\bw_2\in W^{1,r}(\Omega)$, we have $\bw_D=\bu_D$ on $\Gamma_D$, $\nabla\!\cdot \bw_D=0$ in $\Omega$, and, applying the triangle inequality,
$$
\|\bw_D\|_{W^{1,r}(\Omega)^2}\lesssim \|\bu_D \|_{W^{\frac1{r'},r}(\Gamma_D)^2}.
$$

We define the velocity and pressure spaces incorporating the boundary condition and the zero-average constraint in the case $\Gamma_D=\partial\Omega$:
\[
\boldsymbol{V} = \b{W}^{1,r}_{0,\Gamma_D}(\Omega),
\qquad
Q = \begin{cases} L^{r'}(\Omega)&\text{if } |\Gamma_D|<|\partial\Omega| \\ L^{r'}_0(\Omega) \coloneqq \left\{q \in L^{r'}(\Omega)\ : \ \textstyle\int_\Omega q = 0 \right\} &\text{if } \Gamma_D =\partial\Omega. \end{cases}
\]
Assuming $\mathbf{f} \in \b{L}^{r'}(\Omega)$ and $\mathbf{g}\in \b{L}^{r'}(\Gamma_N)$, the weak formulation of problem~\eqref{eq:problem} reads as follows: find $(\bu,p)\in \boldsymbol{V}\times Q$ such that
\begin{equation}
\int_\Omega \bsigma(\boldsymbol{\varepsilon}(\bu)+\boldsymbol{\varepsilon}(\bw_D)) : \boldsymbol{\varepsilon}(\bv)\,\mathrm{d}x
- \int_\Omega p\,\nabla\!\cdot \bv\,\mathrm{d}x
- \int_\Omega q\,\nabla\!\cdot \bu\,\mathrm{d}x
=\int_\Omega \mathbf{f} \cdot \bv\,\mathrm{d}x
+ \int_{\Gamma_N} \mathbf{g} \cdot \bv\,\mathrm{d}s
\label{eq:weak_formulation}
\end{equation}
for all \((\bv,q)\in \boldsymbol{V}\times Q\). The velocity field satisfying~\eqref{eq:problem} with the inhomogeneous Dirichlet condition is
$$
\bu+\bw_D\in \boldsymbol{V}_D := \bigl\{ \bv+\bw_D \in W^{1,r}(\Omega)^2 \;:\; \bv \in\boldsymbol{V} \bigr\}.
$$

The problem can be interpreted as the Euler--Lagrange system associated with the minimization of a convex energy functional. Given a convex potential \(\varphi \colon [0,+\infty) \to \mathbb{R}\), we write
\begin{equation}
\label{eq:primal-saddle-functional}
\mathcal{L}(\bu,p) :=
\int_\Omega \varphi\bigl(\|\boldsymbol{\varepsilon}(\bu)\|\bigr)\,\mathrm{d}x
- \ell(\bu)
- \int_\Omega p\,\nabla\cdot \bu\,\mathrm{d}x
\qquad
\ell(\bv) :=
\int_\Omega \mathbf{f}\cdot \bv\,\mathrm{d}x
+ \int_{\Gamma_N} \mathbf{g}\cdot \bv\,\mathrm{d}s,
\end{equation}
with feasible set \((\bu,p)\in \boldsymbol{V}_D\times Q\). The problem thus becomes to find \((\bu, p)\) such that:
\[
\mathcal{L}(\bu, p)=\inf_{v\in \boldsymbol{V}_D}\,\sup_{q \in Q} \,\mathcal{L}(\bv,q).
\]

We note that the linearization schemes for the velocity--pressure formulation presented in Section~\ref{subsec:linearization} below consistently satisfy the (weak) incompressibility constraint.
Consequently, the primal energy functional \eqref{eq:primal-saddle-functional} may be minimized over a divergence-free feasible set
\[
\mathcal{A}
:=
\Bigl\{
\bv \in \boldsymbol{V}_D \;:\;
\nabla \cdot \bv = 0\, \text{ weakly in } \Omega
\Bigr\}.
\]
This allows us to change it from a saddle point functional to a convex energy functional $\mathcal{J}(\bu)$:
\begin{equation}
\label{eq:primalenergy}
\bu = \operatorname*{arg\,min}_{\bv \in \mathcal{A}}\, \mathcal{J}(\bv),
\qquad
\mathcal{J}(\bv)
:=
\int_\Omega \varphi\bigl(\|\boldsymbol{\varepsilon}(\bv)\|\bigr)\,\mathrm{d}x
-\ell(\bv).
\end{equation}
The stress tensor is related to the potential through
\begin{equation}
\bsigma(\boldsymbol{\varepsilon}(\bu))
=
\varphi'\bigl(\|\boldsymbol{\varepsilon}(\bu)\|\bigr)\,
\frac{\boldsymbol{\varepsilon}(\bu)}{\|\boldsymbol{\varepsilon}(\bu)\|}
\qquad
\text{for } \boldsymbol{\varepsilon}(u) \neq 0.
\end{equation}
In the power-law case, we have $\varphi(s) = \frac{\mu}r s^r$ and $\varphi'(s) = \mu\, s^{r-1}$.

\subsubsection{Three-field Fenchel--dual formulation}

We apply Fenchel duality to the primal velocity--pressure energy functional and derive the corresponding weak formulation. To this end, we introduce the space of symmetric tensor fields
\[
\mathbf{\Sigma} :=
\left\{
\bsigma \in L^{r'}(\Omega;\mathbb{R}^{2\times 2}) \;:\;
\bsigma^T = \bsigma
\quad \text{a.e. in } \Omega
\right\},
\]
along with the dual admissible set
\[
\mathcal{E}
:=
\left\{
(\bsigma,p)\in\mathbf{\Sigma}\times Q:
\int_\Omega
(\bsigma-p\mathbf{I}):\boldsymbol{\varepsilon}(\bv)\,\mathrm{d}x
=
\ell(\bv)
\quad
\forall\,\bv\in \boldsymbol{V}
\right\}.
\]
Let $\varphi^*\colon [0,+\infty) \to \mathbb{R}$ denote the convex conjugate of $\varphi$; in the power-law case, this is given by $\varphi^*(s) := {(r')}^{-1} \mu^{-\frac{1}{r-1}} s^{r'}$. We then define the dual energy functional
\[
\mathcal{J}^*(\btau,q) :=
\int_\Omega \varphi^*(\|\btau\|)\,\mathrm{d}x
- \int_\Omega (\btau-q \mathbf{I}):\beps(\bw_D)\,\mathrm{d}x
+ \ell(\bw_D),
\]
where $\bw_D$ is the Dirichlet lifting introduced in the previous section. As proved in Appendix~\ref{app:duality-gap}, strong duality holds between the primal and dual problems.

\begin{theorem}[Primal--dual equivalence]
\label{thm:primal-dual-equality}
For the continuous power-law Stokes problem, under the notation and assumptions introduced in Section~\ref{sec:continuouspb}, we have
\[
\inf_{\bu\in\mathcal{A}}
\mathcal{J}(\bu)
=
-\inf_{(\bsigma,p)\in\mathcal{E}}
\mathcal{J}^*(\bsigma,p).
\]
\end{theorem}

This identity allows one to define a duality gap and adapt the dual Ka\v{c}anov iteration to the primal scheme following \cite{diening2025guaranteed} (see Section~\ref{subsubsec:primal-dual-kacanov} below). Note that the term $\ell(\bw_D)$ is constant and can therefore be omitted when only the minimizer is sought.

Equivalently, the power-law Stokes problem admits the following three-field dual weak formulation: find $(\bsigma,\bu,p)\in \mathbf{\Sigma}\times \boldsymbol{V}\times Q$ such that
\begin{equation}\label{eq:dual-three-field-continuous}
\begin{aligned}
&\int_\Omega
\mu^{-\frac{1}{r-1}}
\|\bsigma\|^{r'-2}
\bsigma:\btau\,\mathrm{d}x
-
\int_\Omega
\boldsymbol{\varepsilon}(\bu):\btau\,\mathrm{d}x
+
\int_\Omega
\bsigma:\boldsymbol{\varepsilon}(\bv)\,\mathrm{d}x
-
\int_\Omega
p\,\nabla\!\cdot \bv\,\mathrm{d}x
+
\int_\Omega
q\,\nabla\!\cdot \bu\,\mathrm{d}x
\\
&\quad =
\int_\Omega
\mathbf{f}\cdot\bv\,\mathrm{d}x
+
\int_{\Gamma_N}
\mathbf{g}\cdot\bv\,\mathrm{d}s
+ \int_\Omega \boldsymbol{\varepsilon}(\bw_D) : \btau\,\mathrm{d}x
\qquad\text{for all }\,
(\btau,\bv,q)\in \mathbf{\Sigma}\times \boldsymbol{V}\times Q.
\end{aligned}
\end{equation}

\subsubsection{Two-field Fenchel-dual formulation}

We introduce the symmetric total Cauchy stress
$\mathbf{S}
:=
\bsigma\bigl(\boldsymbol{\varepsilon}(\bu)\bigr)-p\mathbf{I}$ and observe that, unlike the effective stress \(\bsigma\), the total stress \(\mathbf{S}\) is
not trace-free. Its deviatoric and volumetric parts satisfy
\[
\operatorname{dev}\mathbf{S}
=
\bsigma\bigl(\boldsymbol{\varepsilon}(\bu)\bigr),
\qquad
p=-\frac{1}{2}\operatorname{tr}(\mathbf{S}).
\]
The dual admissible set, written in terms of the total stress, becomes
\[
\mathcal E_S
:=
\left\{
\mathbf T\in\mathbf\Sigma:
\int_\Omega
\mathbf T:\boldsymbol{\varepsilon}(\bv)\,\mathrm{d}x
=
\ell(\bv)
\quad
\forall\,\bv\in \boldsymbol{V}
\right\}.
\]
Then, we can write the dual minimization problem by expressing $\mathcal{J}^*$ with respect to the total stress, i.e.
\begin{equation}
\mathbf{S}
=
\operatorname*{arg\,min}_{\mathbf{T}\in\mathcal{E}_S}
\mathcal{J}^{*}_S(\mathbf{T}), \qquad
\mathcal{J}^{*}_S(\mathbf{T})
:=
\int_\Omega
\varphi^*\bigl(\|\operatorname{dev}\mathbf{T}\|\bigr)\,\mathrm{d}x
-
\int_\Omega
\mathbf{T}:\boldsymbol{\varepsilon}(\bw_D)\,\mathrm{d}x
+
\ell(\bw_D).
\end{equation}
The corresponding two-field Fenchel-dual formulation reads as follows: find
\(
(\mathbf{S},\bu)\in\mathbf{\Sigma}\times \boldsymbol{V}
\)
such that
\begin{equation}
\label{eq:dual-total-stress-continuous}
\begin{aligned}
&
\int_\Omega
\mu^{-\frac{1}{r-1}}
\left\|\operatorname{dev}\mathbf{S}\right\|^{r'-2}
\operatorname{dev}\mathbf{S}:\mathbf{T}\,\mathrm{d}x
-
\int_\Omega
\boldsymbol{\varepsilon}(\bu):\mathbf{T}\,\mathrm{d}x
+
\int_\Omega
\mathbf{S}:\boldsymbol{\varepsilon}(\bv)\,\mathrm{d}x
\\
&\quad=
\int_\Omega
\mathbf{f}\cdot\bv\,\mathrm{d}x
+
\int_{\Gamma_N}
\mathbf{g}\cdot\bv\,\mathrm{d}s
+ \int_\Omega \boldsymbol{\varepsilon}(\bw_D) : \mathbf{T}\,\mathrm{d}x
\qquad\text{for all }\,
(\mathbf{T},\bv)\in \mathbf{\Sigma}\times \boldsymbol{V}.
\end{aligned}
\end{equation}
Taking
\((\mathbf{T},\bv)=(-q\mathbf{I},\mathbf{0})\) in
\eqref{eq:dual-total-stress-continuous} gives
\[
\int_\Omega q\,\nabla\!\cdot\bu\,\mathrm{d}x=
\int_\Omega q\,\nabla\!\cdot\bw_D\,\mathrm{d}x=0
\qquad
\forall q\in Q,
\]
because \(\operatorname{dev}\mathbf{S}:\mathbf{I}=0\). Thus the
incompressibility equation is contained in the volumetric part of the tensor equation. Likewise,
$
\mathbf{S}:\boldsymbol{\varepsilon}(\bv)
=
\operatorname{dev}\mathbf{S}:
\boldsymbol{\varepsilon}(\bv)
-
p\,\nabla\!\cdot\bv
$,
so the pressure contribution to the momentum equation is contained in the trace of the total stress.

The motivation behind introducing the dual formulation in this context is twofold: the possibility of deriving a simple Ka\v canov scheme that converges for \(r>2\)
and the availability of an explicit piecewise-constant representation of the dual unknown, namely the symmetric deviatoric stress tensor. From a numerical perspective, we investigate whether this formulation can offer a useful
alternative to the primal formulation. On the KS discretization, the Ka\v{c}anov scheme on the
Fenchel-dual formulation (Section~\ref{subsubsec:fenchel-dual-kacanov}) coincides with the dual Ka\v{c}anov
scheme on the primal formulation~\cite{diening2025guaranteed} (Proposition~\ref{prop:dual-kacanov-equivalence}).

\section{Numerical schemes and discretizations}\label{sec:discretepb}

From now on, we restrict our attention to the power-law model, namely \eqref{eq:carreau-yasuda} with $\delta=0$ and
$\mu= \frac 1 2$.
The constitutive law then reads
\begin{equation}\label{eq:constitutive}
\bsigma(\beps(\bu))
=
\nu(\beps(\bu))\,\beps(\bu),
\qquad
\nu(\beps(\bu))
=
\mu\,\|\beps(\bu)\|^{r-2}.
\end{equation}

\subsection{Finite element discretization}

Let $\mathcal T_h$ be a shape-regular triangulation of $\Omega$, and let 
$\boldsymbol{V}_{h}$ and $Q_h$ denote the discrete velocity and pressure spaces. For conforming discretizations, we have
\(
\boldsymbol{V}_{h}\subset \boldsymbol{V}
\text{ and }
Q_h\subset Q,
\)
whereas for the Kouhia--Stenberg method the velocity space is nonconforming, thus it is a subspace of the broken counterpart of $\boldsymbol{V}$.
In this case, differential operators are defined elementwise, i.e., for a piecewise $W^{1,r}$-regular vector field $\bv_h$, we define the broken gradient, symmetric gradient, and divergence by
\[
(\nabla_h\bv_h)|_K
:=
\nabla(\bv_h|_K),
\quad
\boldsymbol{\varepsilon}_h(\bv_h)
:=
\frac12
\left(
\nabla_h\bv_h+\nabla_h\bv_h^T
\right),
\quad
(\nabla_h\!\cdot\bv_h)|_K
:=
\nabla\!\cdot(\bv_h|_K),
\quad
\forall K\in\mathcal T_h.
\]
For conforming finite element functions these operators coincide with the usual derivatives. To ease the notation, in what follows we continue to write
$\nabla\bv_h$,
$\boldsymbol{\varepsilon}(\bv_h)$, and
$\nabla\!\cdot\bv_h$;
whenever a nonconforming space is used, these symbols are understood in the broken sense.
Additionally, we introduce the vector field $\bw_{D,h}$ defined as the interpolation of the lifting $\bw_D$ on the discrete space $\boldsymbol{V}_h$.
A direct construction of 
$$\bw_{D,h}\in \boldsymbol{V}_{D,h}
:= \bigl\{ \bv_h+\bw_{D,h} \;:\; \bv_h \in\boldsymbol{V}_h \bigr\},$$ 
such that \((\nabla\cdot \bw_{D,h}, q_h)_\Omega =0
\text{ for all } q_h\in Q_h\) and 
\(
\bw_{D,h}(\boldsymbol{x}_{D}) 
= \bw_{D}(\boldsymbol{x}_{D}) 
\)
for each boundary degree of freedom $\boldsymbol{x}_D \in \Gamma_D$, is also possible, e.g. by solving a linear, discrete Stokes problem.

We introduce the nonlinear form $a$ together with the pressure--velocity coupling form $b$ defined as
\begin{equation}
a(\bw;\bu,\bv)
:=
\int_\Omega \nu(\beps(\bw))\,\beps(\bu):\beps(\bv)\,\mathrm{d}x,
%=\int_\Omega \mu\,\|\beps(\bw)\|^{r-2}\,\beps(\bu):\beps(\bv)\,\mathrm{d}x,
\qquad
b(\bv,q)
:=
-\int_\Omega q\,\nabla\cdot \bv\,\mathrm{d}x.
\end{equation}
The finite element approximation of the primal problem reads: find $(\bu_h,p_h)\in \boldsymbol{V}_{h}\times Q_h$ such that
\begin{equation}\label{eq:FEM-discrete-abstract}
a(\bu_h + \bw_{D,h};\bu_h + \bw_{D,h},\bv_h)
+b(\bv_h,p_h)
+b(\bu_h,q_h)
=
\ell(\bv_h),
\qquad\forall
(\bv_h,q_h)\in \boldsymbol{V}_{h}\times Q_h.
\end{equation}

\subsection{Finite element spaces}\label{subsec:finite-element-spaces}

All nonlinear solvers presented in Section~\ref{subsec:linearization} below are implemented considering three alternative discretization methods for the primal problem: Taylor--Hood ($\mathbb{P}_2-\mathbb{P}_1$), MINI~\cite{ArnoldBrezziFortin1984}, and Kouhia--Stenberg (KS) finite element pairs. In addition, for the mixed dual formulation we introduce a suitable discrete stress space. The choice of these spaces is motivated as follows:
\begin{itemize}
    \item we do not focus on high-order finite element methods, as they may perform poorly in the presence of singular solutions and thus may offer an unfavorable balance between computational cost and accuracy, especially in highly nonlinear problems;
    \item we include both conforming ($\mathbb{P}_2-\mathbb{P}_1$, MINI) and nonconforming (KS) discretizations;
    \item among the conforming methods, $\mathbb{P}_2-\mathbb{P}_1$ exhibits empirical convergence rates higher than MINI, when the solution is regular ``enough'';
   \item all the considered discretizations are inf-sup stable and do not require additional stabilization terms. None of them is pressure robust; pressure-robust variants, such as Crouzeix--Raviart elements with a reconstruction operator applied to the load
\cite{diening2025pressure}, are not considered here (see Section~\ref{subsec:degenerate-benchmark}).
\end{itemize}

\begin{definition}[Finite element spaces]
Let $\mathcal{T}_h$ be a conforming, shape-regular triangulation of $\Omega\subset\mathbb{R}^2$, and let $\mathbb{P}_k(K)$ denote the space of polynomials of total degree at most $k$ on an element $K\in\mathcal{T}_h$. We define:

\begin{itemize}
    \item \textbf{Lowest-order Taylor--Hood ($\mathbb{P}_2-\mathbb{P}_1$) element pair:}
    \(
    \boldsymbol{W}_h^{\mathrm{P2\text{-}P1}} := \boldsymbol{V}_h^{\mathrm{P2}} \times Q_h^{\mathrm{P1}},
    \)
    where
    \[
    \boldsymbol{V}_h^{\mathrm{P2}}
    :=
    \Bigl\{
    \bv_h \in [H^1(\Omega)]^2
    :
    \bv_h|_K \in [\mathbb{P}_2(K)]^2
    \quad \forall K\in\mathcal{T}_h
    \Bigr\},
    \]
    \[
    Q_h^{\mathrm{P1}}
    :=
    \Bigl\{
    q_h \in H^1(\Omega)
    :
    q_h|_K \in \mathbb{P}_1(K)
    \quad \forall K\in\mathcal{T}_h
    \Bigr\}.
    \]

    \item \textbf{MINI element pair:}
    \(
    \boldsymbol{W}_h^{\mathrm{MINI}} := \boldsymbol{V}_h^{\mathrm{MINI}} \times Q_h^{\mathrm{P1}},
    \)
    where
    \[
    \boldsymbol{V}_h^{\mathrm{MINI}}
    :=
    \Bigl\{
    \bv_h \in [H^1(\Omega)]^2
    :
    \bv_h|_K \in [\mathbb{P}_1(K)\oplus \mathrm{span}\{b_K\}]^2
    \quad \forall K\in\mathcal{T}_h
    \Bigr\},
    \]
    and $b_K$ denotes the cubic scalar bubble function on the triangle $K$.

    \item \textbf{Kouhia--Stenberg (KS) element pair:}
    \(
    \boldsymbol{V}_h^{\mathrm{KS}} := V_{1,h}^{\mathrm{CR1}} \times V_{2,h}^{\mathrm{P1}},
    \;
    \boldsymbol{W}_h^{\mathrm{KS}} := \boldsymbol{V}_h^{\mathrm{KS}} \times Q_h^{\mathrm{P0}},
    \)
    where
    \[
    V_{1,h}^{\mathrm{CR1}}
    :=
    \Bigl\{
    v_h \in L^2(\Omega)
    :
    v_h|_K \in \mathbb{P}_1(K)
    \quad \forall K\in\mathcal{T}_h,
    \ \text{$v_h$ is continuous at edge midpoints}
    \Bigr\},
    \]
    \[
    V_{2,h}^{\mathrm{P1}}
    :=
    \Bigl\{
    v_h \in H^1(\Omega)
    :
    v_h|_K \in \mathbb{P}_1(K)
    \quad \forall K\in\mathcal{T}_h
    \Bigr\},
    \]
    \[
    Q_h^{\mathrm{P0}}
    :=
    \Bigl\{
    q_h \in L^2(\Omega)
    :
    q_h|_K \in \mathbb{P}_0(K)
    \quad \forall K\in\mathcal{T}_h
    \Bigr\}.
    \]
\end{itemize}
\end{definition}

\subsubsection{Dual discrete space for the Fenchel dual formulation}

For the Fenchel dual formulation, which in our implementation is considered for the KS scheme, we retain the same discrete velocity space and introduce the following symmetric tensor-valued stress space.

\begin{definition}[Finite element spaces for the dual problem]
    We let
   \[
    \mathbf{\Sigma}_h^{\mathrm{KS}}
    :=
    \Bigl\{
    \btau_h \in L^2(\Omega;\mathbb{R}^{2\times 2}_{\mathrm{sym}})
    :
    \btau_h|_K \in [\mathbb{P}_0(K)]^{2\times 2}_{\mathrm{sym}}
    \quad \forall K\in\mathcal{T}_h
    \Bigr\}.
    \]
    Accordingly, the discrete spaces associated with the two dual formulations are
\[
\boldsymbol{W}_{h,\mathrm{dual}}^{(2)}
:=
\mathbf{\Sigma}_h^{\mathrm{KS}}\times \boldsymbol{V}_h^{\mathrm{KS}},
\qquad
\boldsymbol{W}_{h,\mathrm{dual}}^{(3)}
:=
\mathbf{\Sigma}_h^{\mathrm{KS}}\times \boldsymbol{V}_h^{\mathrm{KS}}
\times Q_h^{\mathrm{P0}},
\]
corresponding to the unknowns
$(\mathbf{S}_h,\bu_h)$ and
$(\bsigma_h,\bu_h,p_h)$, respectively.
\end{definition}

\subsubsection{Convergence rates}

We now summarize the reference convergence rates used in the numerical experiments. Since the available theoretical results do not cover all the discretization methods considered in this work in a uniform way, we distinguish between literature-based estimates and empirical reference slopes.

\begin{proposition}[Convergence rates for conforming methods]\label{prop:literature-rates}
For inf-sup stable conforming approximations satisfying the assumptions of the abstract framework in \cite{Belenki.Berselli.ea:12}, the velocity error is estimated in the natural distance associated with the nonlinear stress tensor, which is equivalent to the quasi-norm introduced in Definition~\ref{def:quasinorm} below. Accordingly, the theoretical reference convergence rates for the MINI element are
\[
    \|\bu_{\mathrm{ex}}-\bu_h\|_{\mathrm{QN}}
    = \mathcal{O}\!\left(h^{\min\{1,r'/2\}}\right),
    \qquad
    \|p_{\mathrm{ex}}-p_h\|_{L^{r'}(\Omega)}
    = \mathcal{O}\!\left(h^{\min\{2/r',r'/2\}}\right).
\]
\end{proposition}

\begin{remark}[Reference slopes for the Kouhia--Stenberg element]\label{rem:ks-empirical}
The Kouhia--Stenberg element is nonconforming and is therefore not directly covered by the conforming abstract theory of \cite{Belenki.Berselli.ea:12}. Consequently, we do not claim a separate theoretical convergence estimate for this discretization; in the numerical comparisons, the reference slopes are taken to match the theoretical orders of the MINI element. For the nonconforming velocity component, the quasi-norm is evaluated using the broken symmetric gradient.
\end{remark}

\begin{remark}[Convergence behavior and reference slopes for $\mathbb{P}_2$--$\mathbb{P}_1$]\label{rem:p2p1-empirical}
The abstract framework of \cite{Belenki.Berselli.ea:12} also applies to conforming inf-sup stable $\mathbb{P}_2$--$\mathbb{P}_1$ approximations and guarantees the baseline rates stated in Proposition~\ref{prop:literature-rates}. Because the Taylor--Hood element has a higher approximation order than the MINI element, it systematically exhibits faster convergence when applied to smooth solutions. Nevertheless, for benchmarks with low regularity, the observed convergence rate can be limited by the regularity of the exact solution rather than by the polynomial degree of the finite element space; for this reason, we retain the MINI reference slopes in our plots.

Related higher-order \emph{a priori} error estimates are available for dual-mixed formulations of nonlinear generalized Stokes and quasi-Newtonian flow problems under stronger Sobolev regularity assumptions \cite{ErvinHowellStanculescu2008,Howell2009DualMixed,FarhloulZine2017DualMixed}. In particular, these analyses rely on $W^{2,q}$-type regularity assumptions on the solution fields, with some results extending to general higher-order $W^{m,q}$ settings. Although these works address different formulations, unknowns, and discrete spaces, they support the expectation of higher-order convergence rates in our primal setting as well whenever the exact solution is sufficiently regular. 
%Deriving rigorous higher-order error estimates for $\mathbb{P}_2$--$\mathbb{P}_1$ in the present framework is left for future investigation.
\end{remark}

\subsection{Linearization strategies}\label{subsec:linearization}

We employ two nonlinear solvers: the \emph{Newton} method and the \emph{Ka\v canov} iteration. The choice is motivated by their simplicity and by the fact that, taken together, they cover the whole range of values of $r$ considered in this work. In addition to the primal Newton solver, we also implemented a Newton method for the two-field
Fenchel-dual formulation. For the Ka\v canov scheme, which is a fixed-point iteration, we consider
two different dual viewpoints. The first one is a dual-weighted iteration acting on the primal
velocity--pressure formulation \cite{diening2025guaranteed}, which we refer to as \emph{Dual
Ka\v canov iteration on the primal scheme}. The second one is a Ka\v canov iteration applied
directly to the Fenchel-dual formulation. On the KS discretization the two iterations coincide
(cf. Proposition~\ref{prop:dual-kacanov-equivalence} below).

In the form used here, the Newton schemes should be regarded as locally convergent.
The Ka\v canov iterations have a different convergence mechanism, since they require an additional condition \cite{diening2025guaranteed} on $\varphi$ for $r < 2$ and on $\varphi^*$ for $r\ge 2$.

\begin{remark}[Degeneracy as $r\to 1^+$]
\label{rem:deg}
The $N$-function associated with the power-law model is
\[
\varphi_r(t):=\frac{\mu}{r}t^r,
\qquad t\geq 0.
\]
For every $r>1$, the function $\varphi_r$ is strictly convex. However, as $r\to 1^+$, its convexity degenerates: $\varphi_r$ converges pointwise to the linear function $\varphi_1(t)=\mu t$, and the family $\{\varphi_r\}_{r>1}$ is not uniformly convex with constants independent of $r$. This loss of coercivity reflects the increasing degeneracy of the strongly shear-thinning regime and helps explain the deterioration in the performance of some nonlinear solvers.
\end{remark}

\subsubsection{Primal Ka\v canov iteration}

Given $\bu_h^k\in \boldsymbol{V}_{h}$, find $(\bu_h^{k+1},p_h^{k+1})\in \boldsymbol{V}_{h}\times Q_h$ such that
\begin{equation}
a(\bu_h^k+\bw_{D,h};\bu_h^{k+1}+\bw_{D,h},\bv_h)
+b(\bv_h,p_h^{k+1})
+b(\bu_h^{k+1},q_h)
=
\ell(\bv_h), \quad
\forall (\bv_h,q_h)\in \boldsymbol{V}_{h}\times Q_h.
\end{equation}
Equivalently,
$$
\begin{aligned}
&\int_\Omega \nu(\beps(\bu_h^k+\bw_{D,h}))\,\beps(\bu_h^{k+1}+\bw_{D,h}) : \beps(\bv_h)\,\mathrm{d}x
- \int_\Omega p_h^{k+1}\,\nabla\cdot \bv_h\,\mathrm{d}x
\\
&\qquad- \int_\Omega q_h\,\nabla\cdot \bu_h^{k+1}\,\mathrm{d}x
= \ell(\bv_h), \quad
\forall (\bv_h,q_h)\in \boldsymbol{V}_{h}\times Q_h.
\end{aligned}
$$

\subsubsection{Dual Ka\v canov iteration on the primal scheme}\label{subsubsec:primal-dual-kacanov}

Following \cite{diening2025guaranteed}, we may define a dual-weighted Ka\v canov iteration while still solving a primal velocity--pressure problem. Given a stress-like field $\bsigma_h^k$, find $(\bu_h^{k+1},p_h^{k+1})\in \boldsymbol{V}_{h}\times Q_h$ such that
\begin{equation}\label{eq:dual-kacanov-primal}
\int_\Omega
\frac{\|\bsigma_h^k\|}{(\varphi^*)'(\|\bsigma_h^k\|)}
\,\beps(\bu_h^{k+1}+\bw_{D,h}):\beps(\bv_h)\,\mathrm{d}x
-
\int_\Omega p_h^{k+1}\,\nabla\cdot \bv_h\,\mathrm{d}x
-
\int_\Omega q_h\,\nabla\cdot \bu_h^{k+1}\,\mathrm{d}x
=\ell(\bv_h),
\end{equation}
for all $(\bv_h,q_h)\in \boldsymbol{V}_{h}\times Q_h$.
The stress-like field is then updated by
\begin{equation}\label{eq:duality-u-sigma}
\bsigma_h^{k+1}
:=
\frac{\|\bsigma_h^k\|}{(\varphi^*)'(\|\bsigma_h^k\|)}
\,\beps(\bu_h^{k+1}+\bw_{D,h}).
\end{equation}
This scheme keeps the primal velocity--pressure structure, while using a dual quantity to update the viscosity coefficient.

\subsubsection{Ka\v canov iteration on the two-field Fenchel-dual formulation}
\label{subsubsec:fenchel-dual-kacanov}

Given $\mathbf{S}_h^k\in\mathbf{\Sigma}_h^{\mathrm{KS}}$, find
\(
(\mathbf{S}_h^{k+1},\bu_h^{k+1})
\in
\mathbf{\Sigma}_h^{\mathrm{KS}}\times \boldsymbol{V}_{h}^{\mathrm{KS}}
\)
such that
\begin{equation}\label{eq:dual-two-field-kacanov}
\begin{aligned}
&\int_\Omega
\mu^{-\frac{1}{r-1}}
\|\dev(\mathbf{S}_h^k)\|^{r'-2}
\dev(\mathbf{S}_h^{k+1}):\mathbf{T}_h\,\mathrm{d}x
\\
&\qquad
-\int_\Omega
\beps(\bu_h^{k+1}):\mathbf{T}_h\,\mathrm{d}x
+
\int_\Omega
\mathbf{S}_h^{k+1}:\beps(\bv_h)\,\mathrm{d}x
=
\ell(\bv_h)
+
\int_\Omega
\beps(\bw_{D,h}):\mathbf{T}_h\,\mathrm{d}x
\end{aligned}
\end{equation}
for all
\(
(\mathbf{T}_h,\bv_h)
\in
\mathbf{\Sigma}_h^{\mathrm{KS}}\times \boldsymbol{V}_{h}^{\mathrm{KS}}.
\)
The pressure is reconstructed as
$
p_h^{k+1}
=
-\frac{1}{2}\tr(\mathbf{S}_h^{k+1})
$.
Duality is essential in order to satisfy one of the sufficient conditions for the convergence of the Ka\v canov scheme \cite{diening2025guaranteed}. In particular, for $r\geq 2$, this condition may fail for $\varphi$, whereas it holds for the convex conjugate $\varphi^*$. For a general $N$-function $a$, the relevant condition reads
$
a''(t)\le t^{-1}a'(t)
$.

\begin{proposition}[Equivalence of the two dual Ka\v{c}anov iterations on KS]
\label{prop:dual-kacanov-equivalence}
Assume that $\mu$ is constant on each $K\in\mathcal{T}_h$. Let $\bsigma_h^0\in\mathbf{\Sigma}_h^{\mathrm{KS}}$ and
$\mathbf{S}_h^0\in\mathbf{\Sigma}_h^{\mathrm{KS}}$ satisfy $\dev(\mathbf{S}_h^0)=\bsigma_h^0$, and suppose that
$\bsigma_h^k|_K\neq\mathbf 0$ for every $K\in\mathcal{T}_h$ and every $k$, so that
both iterations are well defined. Then, for every $k\ge 0$, the iterate $(\bu_{h,~3F}^{k+1},p_h^{k+1},\bsigma_h^{k+1})$ of the dual-weighted Ka\v{c}anov iteration \eqref{eq:dual-kacanov-primal}-\eqref{eq:duality-u-sigma} and the iterate $(\bu_{h,~2F}^{k+1}, \mathbf{S}_h^{k+1})$ of the Fenchel-dual Ka\v{c}anov iteration \eqref{eq:dual-two-field-kacanov} are related by
\[
  \bu_{h, ~3F}^{k+1}=\bu_{h, ~2F}^{k+1},\qquad
  \bsigma_h^{k+1}=\dev(\mathbf{S}_h^{k+1}),\qquad
  p_h^{k+1}=-\tfrac12\tr(\mathbf{S}_h^{k+1}).
\]
\end{proposition}
The proof is given in Appendix~\ref{app:proof-equivalence-dual-schemes}.
\subsubsection{Newton method on the primal scheme}

We write the discrete nonlinear problem in residual form: find $(\bu_h,p_h)\in \boldsymbol{V}_{h}\times Q_h$ such that
\begin{equation}
\mathcal{R}(\bu_h,p_h;\bv_h,q_h)=0
\qquad
\forall (\bv_h,q_h)\in \boldsymbol{V}_{h}\times Q_h,
\end{equation}
where
\begin{equation}
\mathcal{R}(\bu,p;\bv,q)
:=
a(\bu+\bw_{D,h};\bu+\bw_{D,h},\bv)
+b(\bv,p)
+b(\bu,q)
-\ell(\bv).
\end{equation}
Given an iterate $(\bu_h^k,p_h^k)$, the Newton step $(\delta \bu_h^k,\delta p_h^k)\in \boldsymbol{V}_{h}\times Q_h$ is defined by
\begin{equation}
D\mathcal{R}(\bu_h^k,p_h^k)[\delta \bu_h^k,\delta p_h^k](\bv_h,q_h)
=
-\mathcal{R}(\bu_h^k,p_h^k;\bv_h,q_h)
\end{equation}
for all $(\bv_h,q_h)\in \boldsymbol{V}_{h}\times Q_h$,
and the iterates are updated as
$\bu_h^{k+1}=\bu_h^k+\delta \bu_h^k$
and $p_h^{k+1}=p_h^k+\delta p_h^k$.
Therefore, the Newton increment satisfies the following problem, obtained by linearizing the nonlinear term via its Fr\'echet derivative:
\begin{align}
&\int_\Omega \mu\,\|\beps(\bu_h^k+\bw_{D,h})\|^{r-2}\,\beps(\delta \bu_h^k):\beps(\bv_h)\,\mathrm{d}x
\nonumber\\
&\qquad+
(r-2)\int_\Omega \mu\,\|\beps(\bu_h^k+\bw_{D,h})\|^{r-4}
\bigl(\beps(\bu_h^k+\bw_{D,h}):\beps(\delta \bu_h^k)\bigr)
\bigl(\beps(\bu_h^k+\bw_{D,h}):\beps(\bv_h)\bigr)\,\mathrm{d}x
\nonumber\\
&\qquad
-\int_\Omega \delta p_h^k\,\nabla\cdot \bv_h\,\mathrm{d}x
-\int_\Omega q_h\,\nabla\cdot \delta \bu_h^k\,\mathrm{d}x
=
-\mathcal{R}(\bu_h^k,p_h^k;\bv_h,q_h),
\qquad \forall(\bv_h,q_h)\in \boldsymbol V_{h}\times Q_h.
\end{align}

\subsubsection{Newton method on the two-field Fenchel-dual formulation}

We introduce
\begin{equation}
\begin{aligned}
\mathcal{R}_2^*(\mathbf{S},\bu;\mathbf{T},\bv)
:={}&
\int_\Omega
\mu^{-\frac{1}{r-1}}
\|\dev(\mathbf{S})\|^{r'-2}
\dev(\mathbf{S}):\mathbf{T}\,\mathrm{d}x
-
\int_\Omega \beps(\bu):\mathbf{T}\,\mathrm{d}x
\\
&+
\int_\Omega \mathbf{S}:\beps(\bv)\,\mathrm{d}x
-
\ell(\bv)
-
\int_\Omega \beps(\bw_{D,h}):\mathbf{T}\,\mathrm{d}x.
\end{aligned}
\end{equation}
Given the solution at the previous iteration $(\mathbf{S}_h^k,\bu_h^k)$,
%(or initial guess for $k=0$),
find
\(
(\delta\mathbf{S}_h^k,\delta\bu_h^k)
\in
\mathbf{\Sigma}_h^{\mathrm{KS}}\times \boldsymbol{V}_{h}^{\mathrm{KS}}
\)
such that
\begin{align}
&\int_\Omega
\mu^{-\frac{1}{r-1}}
\|\dev(\mathbf{S}_h^k)\|^{r'-2}
\dev(\delta\mathbf{S}_h^k):\mathbf{T}_h\,\mathrm{d}x
\nonumber\\
&\quad+
(r'-2)\int_\Omega
\mu^{-\frac{1}{r-1}}
\|\dev(\mathbf{S}_h^k)\|^{r'-4}
\bigl(
\dev(\mathbf{S}_h^k):
\dev(\delta\mathbf{S}_h^k)
\bigr)
\bigl(
\dev(\mathbf{S}_h^k):\mathbf{T}_h
\bigr)\,\mathrm{d}x
\nonumber\\
&\quad-
\int_\Omega
\beps(\delta\bu_h^k):\mathbf{T}_h\,\mathrm{d}x
+
\int_\Omega
\delta\mathbf{S}_h^k:\beps(\bv_h)\,\mathrm{d}x
\nonumber =
-\mathcal{R}_2^*
(\mathbf{S}_h^k,\bu_h^k;\mathbf{T}_h,\bv_h)
\end{align}
for all
\(
(\mathbf{T}_h,\bv_h)
\in
\mathbf{\Sigma}_h^{\mathrm{KS}}\times \boldsymbol{V}_{h}^{\mathrm{KS}}
\). The iterates are updated according to
\[
\mathbf{S}_h^{k+1}
=
\mathbf{S}_h^k+\delta\mathbf{S}_h^k,
\qquad
\bu_h^{k+1}
=
\bu_h^k+\delta\bu_h^k,
\qquad
p_h^{k+1}
=
-\frac{1}{2}\tr(\mathbf{S}_h^{k+1}).
\]

% \subsubsection{Equivalence of the two dual Ka\v{c}anov iterations on KS}

\subsection{Regularization and relaxation techniques}

We present two regularization and relaxation strategies, both acting only on the viscosity; unless otherwise stated, the rest of the model remains unchanged. These techniques are introduced for large or close to $1$ values of $r$, in order to reduce the dependence of the nonlinear solver on a carefully tuned initial guess, as discussed in Section~\ref{sec:experimentalsetting} below.

\subsubsection{\texorpdfstring{$\delta_h$}{δₕ}-regularization}

The $\delta_h$-regularization is motivated by the fact that, for $r<2$, the constitutive law becomes singular as $\|\beps(\bu)\|\to 0$. Introducing a positive parameter $\delta_h$ prevents this degeneracy and improves stability.

\begin{definition}[$\delta_h$-regularization of \eqref{eq:constitutive}]
We define the $\delta_h$-regularized viscosity by
\[
\nu_{\delta_h}(\beps(\bu))
:=
\mu\bigl(\|\beps(\bu)\|+\delta_h\bigr)^{r-2},
\qquad \text{with }
\delta_h := k\,h^2,~
k\ll \frac{1}{h}
\]
\end{definition}

The scaling $\delta_h=\mathcal{O}(h^2)$ ensures that the regularization error vanishes as $h\to0$.
On a fixed mesh, however, it may dominate the discretization error for large values of $k$, in
particular for the $\mathbb{P}_2$--$\mathbb{P}_1$ element, whose observed rates on smooth
solutions are close to $\mathcal{O}(h^2)$; see Section~\ref{sec:numerical-results}.

\subsubsection{\texorpdfstring{$r$}{r}-continuation}

The continuation strategy is based on two heuristic observations: first, problems with more extreme values of $r$ tend to be harder to solve because of the stronger nonlinearity; second, problems corresponding to nearby values of $r$ are expected to have similar solutions. Therefore, the target exponent can be approached gradually through a sequence of intermediate exponents.

\begin{definition}[$r$-continuation for \eqref{eq:constitutive}]
Let
\(
    r_0,r_1,\ldots,r_N=r
\)
be a sequence of intermediate exponents, typically starting from $r_0=2$. At the $k$-th continuation step, we consider the constitutive law
\[
    \nu_{r_k}(\beps(\bu))
    :=
    \mu \|\beps(\bu)\|^{r_k-2},
\]
and use the converged discrete solution corresponding to $r_{k-1}$ as the initial guess for the nonlinear problem corresponding to $r_k$.
\end{definition}
The following result provides a mathematical justification for this strategy at the discrete level.

\begin{proposition}[Continuity of the discrete solution with respect to $r$]
\label{prop:r-continuity}
For a fixed mesh $\mathcal{T}_h$, consider any of the finite element discretizations introduced in Section \ref{subsec:finite-element-spaces}.
Let
\[
    I=[r_-,r_+]\Subset(1,\infty),
    \qquad
    0<\mu_0\leq\mu(x)\leq\mu_1<\infty
    \quad\text{for a.e. }x\in\Omega,
\]
and assume that the load functional $\ell$ is independent of the exponent. For every $s\in I$, let
\(
    (\boldsymbol{u}_h(s),p_h(s))
    \in \boldsymbol{V}_{h}\times Q_h
\)
be the solution of~\eqref{eq:FEM-discrete-abstract}, with the viscosity $\nu$ replaced by
\(
    \nu_s(\boldsymbol{\varepsilon}(\boldsymbol{u}))
    :=
    \mu
    \lVert\boldsymbol{\varepsilon}(\boldsymbol{u})\rVert^{s-2},
\)
where
$\nu_s(\boldsymbol{\varepsilon}(\boldsymbol{u}))
 \boldsymbol{\varepsilon}(\boldsymbol{u})$
is understood to be zero wherever
$\boldsymbol{\varepsilon}(\boldsymbol{u})=\boldsymbol{0}$. If $r_k\to r$, with $r_k,r\in I$, then
\[
    \boldsymbol{u}_h(r_k)
    \rightarrow
    \boldsymbol{u}_h(r)
    \qquad\text{in }\boldsymbol{V}_h.
\]
Moreover, if the pressure is uniquely determined, then
\(
    p_h(r_k)
    \rightarrow
    p_h(r)
    \;\text{in }Q_h
\).
\end{proposition}

\begin{proof}
We introduce the discrete divergence-free space
\(
    \boldsymbol{Z}_{h}
    :=
    \left\{
        \boldsymbol{v}_h\in \boldsymbol{V}_{h}:
        b(\boldsymbol{v}_h,q_h)=0
        \quad\forall q_h\in Q_h
    \right\}
\).
For $s\in I$, we define the operator
$A_h(\cdot,s):\boldsymbol{Z}_{h}\to \boldsymbol{Z}_{h}^{\ast}$ by
\[
\begin{aligned}
    \left\langle
        A_h(\boldsymbol{w}_h,s),\boldsymbol{v}_h
    \right\rangle
    &:=
    \int_\Omega
        \mu
        \left\lVert
            \boldsymbol{\varepsilon}
            (\boldsymbol{w}_h+\bw_{D,h})
        \right\rVert^{s-2}
        \boldsymbol{\varepsilon}
        (\boldsymbol{w}_h+\bw_{D,h})
        :
        \boldsymbol{\varepsilon}(\boldsymbol{v}_h)
    \,\mathrm{d}x .
\end{aligned}
\]
The discrete velocity problem can then be written as
\(
    A_h(\boldsymbol{w}_h,s)
    =
    \ell|_{\boldsymbol{Z}_{h}}
\).
For every $s>1$, let
\[
    F_s(\boldsymbol{\xi})
    :=
    \begin{cases}
        \lVert\boldsymbol{\xi}\rVert^{s-2}\boldsymbol{\xi},
        &\boldsymbol{\xi}\neq\boldsymbol{0},\\[2mm]
        \boldsymbol{0},
        &\boldsymbol{\xi}=\boldsymbol{0}.
    \end{cases}
\]
Since
\(
    \lVert F_s(\boldsymbol{\xi})\rVert
    =
    \lVert\boldsymbol{\xi}\rVert^{s-1},
\)
the mapping
\(
    (\boldsymbol{\xi},s)
    \mapsto F_s(\boldsymbol{\xi})
\)
is continuous on
$\mathbb{R}_{\mathrm{sym}}^{2\times 2}\times I$.
Because $\boldsymbol{Z}_{h}$ is finite-dimensional, it follows that
\(
    (\boldsymbol{w}_h,s)
    \mapsto A_h(\boldsymbol{w}_h,s)
\)
is continuous from $\boldsymbol{Z}_{h}\times I$ into $\boldsymbol{Z}_{h}^{\ast}$.
For every fixed $s\in I$, the mapping $F_s$ is strictly monotone.
Consequently,
\[
\begin{aligned}
    &\left\langle
        A_h(\boldsymbol{w}_h,s)
        -
        A_h(\boldsymbol{z}_h,s),
        \boldsymbol{w}_h-\boldsymbol{z}_h
    \right\rangle
    \\
    &=
    \int_\Omega
        \mu
        \left[
            F_s\left(
                \boldsymbol{\varepsilon}
                (\boldsymbol{w}_h+\bw_{D,h})
            \right)
            -
            F_s\left(
                \boldsymbol{\varepsilon}
                (\boldsymbol{z}_h+\bw_{D,h})
            \right)
        \right]
        :
        \boldsymbol{\varepsilon}
        (\boldsymbol{w}_h-\boldsymbol{z}_h)
    \,\mathrm{d}x
    >0
\end{aligned}
\]
whenever
$\boldsymbol{w}_h,\boldsymbol{z}_h\in \boldsymbol{Z}_{h}$ and
$\boldsymbol{w}_h\neq\boldsymbol{z}_h$.
The lower bound on $\mu$, together with the discrete Korn inequality~\cite{CarstensenSchedensack2015},
also implies that the family
\(
    \bigl(A_h(\cdot,s)\bigr)_{s\in I}
\)
is uniformly coercive on $\boldsymbol{Z}_{h}$. Thus, for every $s\in I$, the
discrete velocity problem admits a unique solution.
The assumptions of the parameter-dependent monotonicity result
of \cite{Beldzinski2022} are therefore
satisfied. In particular, the parameter-dependent condition
$(S_2)$ is automatic because $\boldsymbol{Z}_{h}$ is finite-dimensional and
weak and strong convergence coincide in finite-dimensional spaces.
It follows that
\(
    \boldsymbol{w}_h(r_k)
    \rightarrow
    \boldsymbol{w}_h(r)
\)
in $\boldsymbol{Z}_{h}$.
To prove the convergence of the pressure, we subtract the discrete
momentum equations corresponding to $r_k$ and $r$. This gives
\[
\begin{aligned}
    b\bigl(
        \boldsymbol{v}_h,
        p_h(r_k)-p_h(r)
    \bigr)
    =
    -\int_\Omega
        \mu
        \Big[
            F_{r_k}\bigl(
                \boldsymbol{\varepsilon}
                (\bw_h(r_k)+\bw_{D,h})
            \bigr)
            -
            F_r\bigl(
                \boldsymbol{\varepsilon}
                (\bw_h(r)+\bw_{D,h})
            \bigr)
        \Big]
        :
        \boldsymbol{\varepsilon}(\boldsymbol{v}_h)
    \,\mathrm{d}x
\end{aligned}
\]
for every $\boldsymbol{v}_h\in \boldsymbol{V}_{h}$.
The right-hand side converges to zero by the continuity of
$(\boldsymbol{\xi},s)\mapsto F_s(\boldsymbol{\xi})$ and the
convergence of the discrete velocities.
The discrete inf--sup condition yields
\(
    p_h(r_k)
    \longrightarrow
    p_h(r)
\)
in $Q_h$.
%This concludes the proof of Proposition~\ref{prop:r-continuity}.
\end{proof}

Proposition~\ref{prop:r-continuity} shows that, for a fixed finite element discretization, the discrete velocity and pressure depend continuously on the power-law exponent. In particular, if the continuation increments satisfy $|r_{k+1}-r_k|\to0$, then
\(
    |\,\boldsymbol{u}_h(r_k)
    -
    \boldsymbol{u}_h(r_{k+1})\,|
    \rightarrow
    0.
\)
This justifies using the converged solution at $r_k$ as the initial guess for solving the nonlinear problem at $r_{k+1}$. 
% The proof of the proposition is based on the parameter-dependent monotonicity framework of \cite{Beldzinski2022}.

\section{Experimental setting}\label{sec:experimentalsetting}

\subsection{Test problems}

We consider two classes of numerical tests. First, we use three manufactured benchmark solutions in order to assess accuracy and convergence rates: a \emph{trigonometric} solution, a \emph{power-type} solution taken from \cite{diening2025pressure}, and a \emph{degenerate} solution on a cross-shaped set in the middle of the domain, inspired by \cite{dipietro2020improved}. Second, in Appendix~\ref{app:obstacle-domain-test}, we consider an obstacle-domain flow problem without an exact solution, in order to assess the robustness of the proposed methods on a geometrically nontrivial domain.

\subsubsection{Manufactured benchmark solutions}
\label{subsubsec:benchmark-problems}

We consider three manufactured solutions designed to test the numerical methods under different regularity and degeneracy conditions.

\paragraph{Trigonometric benchmark.}
Let
\(
\Omega_{\mathrm{trig}}=[0,1]^2
\).
The exact solution is
\[
\boldsymbol u_{\mathrm{ex}}(x,y)
=
\begin{pmatrix}
x\sin(xy)\\
-y\sin(xy)
\end{pmatrix},
\qquad
p_{\mathrm{ex}}(x,y)=\arctan(xy).
\]
This smooth benchmark is used as a reference case for assessing the
accuracy of the finite element discretizations and the robustness of the
nonlinear solvers.

\paragraph{Power-type benchmark.}
Let
\(
\Omega_{\mathrm{power}}=[-1,1]^2,\,
\text{ and }\,
\rho(x,y)=\sqrt{x^2+y^2}.
\)
We consider
\[
\boldsymbol u_{\mathrm{ex}}(x,y)
=
\rho(x,y)^{0.01}
\begin{pmatrix}
y\\
-x
\end{pmatrix},
\qquad
p_{\mathrm{ex}}(x,y)=\eta\,\rho(x,y)^\gamma,
\]
where $\eta$ is defined by
\[
\eta :=
\begin{cases}
0.01, & 1<r<2,\\[0.3em]
1.0, & r\geq 2.
\end{cases}
\]
and
\(
\gamma=\frac{2}{r}-1+0.01
\).
The solution has limited regularity at the interior point $(0,0)$ and is
therefore used to investigate the effect of a localized singularity on
the discretization error and nonlinear convergence.

\paragraph{Degenerate benchmark.}
Let
\(
\Omega_{\mathrm{deg}}=[0,1]^2,
\;
X=x-\frac12,
\text{ and }\,
Y=y-\frac12
\).
For $\xi:=r+\frac{3}{4}$, we define
\[
\begin{aligned}
\boldsymbol u_{\mathrm{ex}}
&=
\begin{pmatrix}
u_{\mathrm{ex},1}\\
u_{\mathrm{ex},2}
\end{pmatrix},
\\
u_{\mathrm{ex},1}
&=
-\frac{1}{10}
Y|Y|^\xi |X|^{\xi+2}
\left(\xi+2-10\xi|Y|^\xi\right)
e^{-10(|X|^\xi+|Y|^\xi)},
\\
u_{\mathrm{ex},2}
&=
\frac{1}{10}
X|X|^\xi |Y|^{\xi+2}
\left(\xi+2-10\xi|X|^\xi\right)
e^{-10(|X|^\xi+|Y|^\xi)}.
\end{aligned}
\]
The exact pressure is
\(
p_{\mathrm{ex}}(x,y)=|X|^\xi
\), normalized to zero mean.
The velocity is divergence-free and its symmetric gradient vanishes on
the cross-shaped set
\(
\left\{x=1/2\right\}
\cup
\left\{y=1/2\right\}
\).
This benchmark therefore contains an extended zero-strain region, where
the effective viscosity is singular for $1<r<2$ and vanishes for $r>2$.
It is the most challenging test considered for assessing solver robustness.
The different domains are intentional. In the power-type
benchmark, the loss of regularity is localized at a single interior
point, whereas in the degenerate benchmark the zero-strain set extends
along two intersecting lines.
\subsection{Error metrics}

The error metrics employed are the \(L^{r'}\)-error for the pressure \(p\), and the quasi-norm (QN) error for the velocity \(\bu\). Following \cite{diening2025pressure}, the quasi-norm is defined as follows.

\begin{definition}[Quasi-norm]
\label{def:quasinorm}
Let $\bu_{\mathrm{ex}}$ be the exact solution and $\bu_h$ its discrete approximation. We define the quasi-norm error by
\[
\|\bu_{\mathrm{ex}}-\bu_h\|_{\mathrm{QN}}
:=
\left(
\int_\Omega
\bigl(
\bsigma(\beps(\bu_{\mathrm{ex}}))
-
\bsigma(\beps(\bu_h))
\bigr)
:
\bigl(
\beps(\bu_{\mathrm{ex}})
-
\beps(\bu_h)
\bigr)\,dx
\right)^{1/2}.
\]
Moreover, we also define
\(
\|\bu\|_{\mathrm{QN}}
:=
\left(
\int_\Omega
\bsigma(\beps(\bu)):\beps(\bu)\,dx
\right)^{1/2}
\).
\end{definition}

For values of \(r\in(1,1.1)\), we use the \(L^2\)-norm in place of the \(L^{r'}\)-norm in the pressure error in order to avoid numerical instabilities in the post-processing of the results.

\subsection{Implementation details}

All finite element discretizations and nonlinear solvers were implemented in \textsc{FEniCS} using DOLFIN 2019.1.0~\cite{LoggWells2010}. The linear systems arising at each nonlinear iteration were solved with MUMPS \cite{MUMPS:1,MUMPS:2}, which allowed us to treat the finest meshes considered in this work while retaining the robustness of a direct algebraic solver. All experiments were run on a MacBook Pro equipped with an M4 Pro chip.

In the implementation, the tolerance $\mathrm{tol}=\max\{10^{-8},h^3\}$ is used to define the stopping criteria described below; its mesh-dependent part avoids oversolving on coarse meshes.
We use as an initial guess, for a given problem with parameter $r$, the solution of the corresponding linear problem obtained by freezing the viscosity at the constant value $\nu_{\mathrm{guess}}:=0.2$, corresponding to the Newtonian case $r=2$. This provides a sufficiently accurate and standard initial guess, while avoiding a more complicated search procedure.
For the tests regularized via \(\delta_h\), we set \(\delta_h := k h^2\), where \(k\in\mathbb{R}\) is a constant. Since its value may vary among different tests, it is specified each time.

\subsection{Stopping criteria}

Different stopping criteria are implemented, tested, and compared for alternative scenarios. In particular, the exact-solution-based criterion is used only when a manufactured solution is available. The increment-based criterion is simple to implement, but has some limitations, whereas the energy-step criterion provides a meaningful residual-type indicator and proves to be particularly effective in practice, although this improved performance might come with a trade-off in terms of higher computational cost.
The increment-based criterion stops the iteration when, for both fields,
\begin{equation}\label{eq:increment-criteria}
\|\bu_h^k-\bu_h^{k-1}\|_{\mathrm{QN}}
\le \mathrm{atol}+\mathrm{tol}\,\|\bu_h^{k-1}\|_{\mathrm{QN}},
\qquad
\|p_h^k-p_h^{k-1}\|_{L^{r'}(\Omega)}
\le \mathrm{atol}+\mathrm{tol}\,\|p_h^{k-1}\|_{L^{r'}(\Omega)},
\end{equation}
with $\mathrm{atol}=10^{-8}$ and $\mathrm{tol}=\max\{10^{-8},h^3\}$. For the KS element, the velocity increment is measured in the sum of the $H^1$ norms of its two components.
% Although this criterion is very easy to implement, it may lead to either premature stopping, in particular when the increment is not monotone, or late stopping, as shown in the following.
% \subsubsection{Energy-step criterion}
\noindent
The energy-step criterion is based on the primal energy functional defined in \eqref{eq:primalenergy}. 
% The purpose of this criterion is to stop the nonlinear iteration when the discrete energy has essentially stabilized, namely when two consecutive iterates produce only a negligible relative variation of the energy.
We introduce
\begin{equation}
   \mathcal{E}_{\mathrm{energy}}^k
=
\frac{\left| \mathcal{J}(\bu_h^k) - \mathcal{J}(\bu_h^{k-1})\right|}
{\left|\mathcal{J}(\bu_h^{k-1})\right|}
\end{equation}
\noindent
The nonlinear iteration is stopped when
\(
\mathcal{E}_{\mathrm{energy}}^k
\leq \mathrm{tol}
\).
The numerical experiments show that this stopping rule is effective provided that the tolerance is chosen consistently with the discretization error. In
particular, the energy-step criterion can stop the solver at the point where the nonlinear error has become negligible with respect to the finite element error. Beyond this point, further nonlinear iterations only reduce the energy, but do not lead to a visible improvement in the final error curve.

\begin{figure}[ht]
    \begin{subfigure}{0.48\textwidth}
    \centering
    \includegraphics[width=\linewidth, trim=0.5cm 2cm 0cm 0cm, clip]{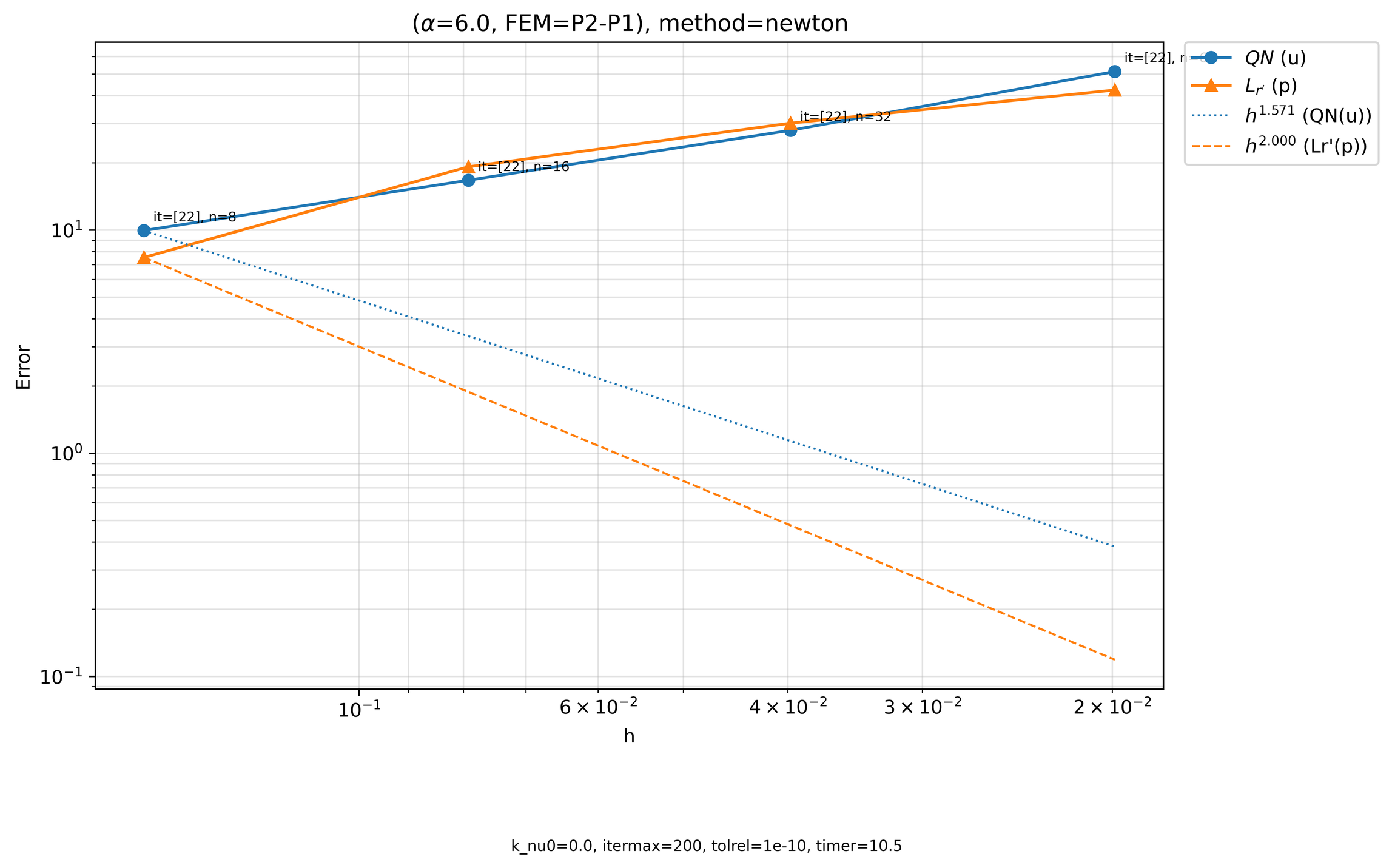}
    \caption{Default \textsc{FEniCS} residual stopping criterion.}
\end{subfigure}\hfill
\begin{subfigure}{0.48\textwidth}
    \centering
    \includegraphics[width=\linewidth, trim=0.5cm 2cm 0cm 0cm, clip]{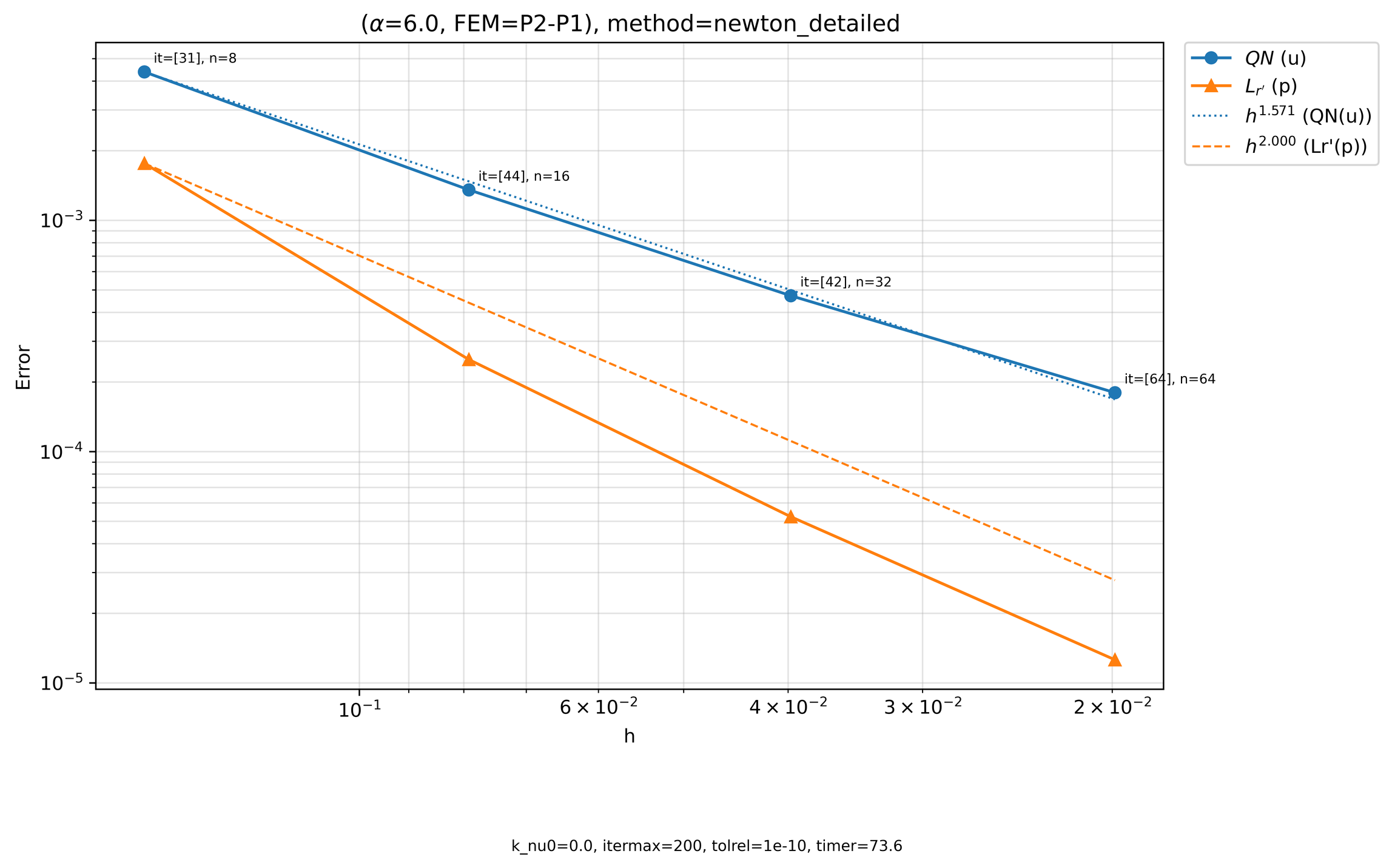}
    \caption{Our energy-step stopping criterion.}
\end{subfigure}
\caption{Comparison of stopping strategies in Newton, $r=8$, \(\mathbb{P}_2-\mathbb{P}_1\), Newton linearization. Plots show convergence (error) with respect to mesh refinement ($h$).}
\label{fig:newton-l2-stop}
\end{figure}

%\subsection{Need for a custom Newton stopping wrapper}
In Figure~\ref{fig:newton-l2-stop}, we show that we indeed require either the increment criterion or the energy step criterion, and thus need a custom wrapper of \textsc{FEniCS}' Newton method, because its standard weak residual or increment prevents our method from converging.

\section{Numerical results}
\label{sec:numerical-results}

We assess the accuracy and robustness of the finite element discretizations and nonlinear solvers on the three manufactured benchmarks introduced in Section~\ref{subsubsec:benchmark-problems}.
For each test, we compare the direct non-regularized solve with $\delta_h$-regularization, $r$-continuation, and their combination.
In some cases, such as in Figures~\ref{fig:trigonometric-summary}, \ref{fig:power-summary}, \ref{fig:degenerate-summary}, as well as in selected remarks, we report only the $\mathbb{P}_2$--$\mathbb{P}_1$ and KS methods, since these are sufficient to provide representative examples of a conforming and a nonconforming FEM, respectively.

\subsection{Trigonometric benchmark}
\label{subsec:trigonometric-results}

\begin{figure}[hbtp]
    \centering

    \includegraphics[width=0.98\linewidth]{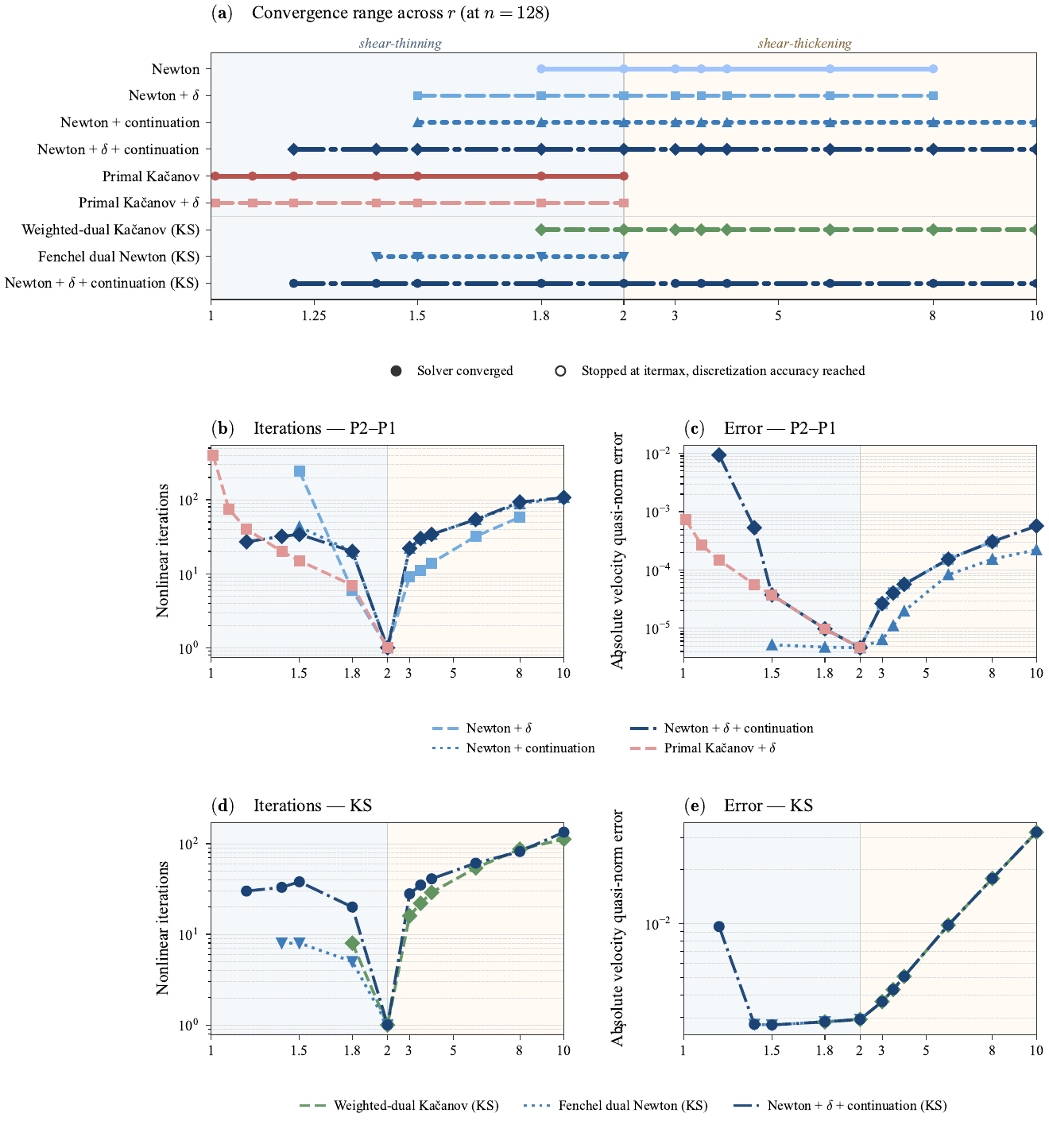}
    \caption{
    Summary of the nonlinear solver performance for the trigonometric benchmark on the finest mesh, with $n_h=128$. Panel~(a) reports the sampled exponents for which each strategy converges. Panels~(b)--(c) show, for the $\mathbb{P}_2$--$\mathbb{P}_1$ discretization, the nonlinear iteration counts and the velocity quasi-norm errors, respectively. Panels~(d)--(e) report the corresponding quantities for the KS discretization and the dual formulations. The blue and beige backgrounds distinguish the shear-thinning and shear-thickening regimes. Lines connect sampled values only and are included as visual guides.
    }
    \label{fig:trigonometric-summary}
\end{figure}

The trigonometric benchmark provides a smooth reference problem for assessing both the finite element discretizations and the nonlinear solvers. Since the exact solution has no interior singularities or extended zero-strain regions, the behaviour observed in this test is mainly determined by the power-law nonlinearity. It therefore provides the baseline against which the power and degenerate benchmarks are interpreted.
Figure~\ref{fig:trigonometric-summary} summarizes the solver robustness, nonlinear iteration counts, and velocity errors for the trigonometric benchmark on the finest mesh, with mesh refinement parameter $n_h=128$. At $r=2$, the constitutive law becomes linear and all schemes terminate after one nonlinear iteration.

\paragraph{Agreement between converged solutions.}

A first observation, also confirmed by the other manufactured benchmarks, is that the nonlinear strategy does not significantly affect the final discrete solution. Whenever Newton, primal Ka\v{c}anov, or Newton with continuation converge for the same $\mathbb{P}_2$--$\mathbb{P}_1$ problem, their velocity and pressure errors differ by less than one percent.
The differences between the methods lie in the nonlinear robustness and computational cost. The only exception is the use of $\delta_h$-regularization, which changes the constitutive law and may introduce a regularization error.
By Proposition~\ref{prop:dual-kacanov-equivalence}, the dual-weighted and Fenchel-dual Ka\v{c}anov
iterations on KS produce identical iterates in exact arithmetic; we therefore report only the
former. 
\paragraph{Shear-thinning regime.}

Without regularization, primal Newton converges for the tested values $r\geq1.8$, but fails at $r=1.5$. Primal Ka\v{c}anov is considerably more robust and remains convergent down to $r=1.01$. Its iteration count, however, increases rapidly as $r\to1^+$.
For the $\mathbb{P}_2$--$\mathbb{P}_1$ discretization, the number of iterations is
\(
7,\ 15,\ 41,\ 76,\ 420
\)
for
\(
r=1.8,\ 1.5,\ 1.2,\ 1.1,\ 1.01,
\)
respectively. A logarithmic least-squares fit gives
\(
N_{\mathrm{it}}
\approx C(r-1)^{-0.91}.
\)
More specifically, the \emph{ratios of consecutive velocity increments} $q_{\mathrm{primal}}$ are consistent with the benchmark-independent relation
\begin{equation}\label{eq:emp-rel-prim}
q_{\mathrm{primal}}\approx 2-r,
\qquad 1<r<2.
\end{equation}
The measured velocity-increment ratio thus approaches one as $r\to1^+$. A similar dependence on $r$ appears in the energy convergence analysis of Heid and S\"uli~\cite{HeidSuli2022} for the relaxed power-law model. This explains why robustness close to the limiting exponent does not imply uniform computational efficiency and it is consistent with the loss of uniform convexity discussed in Remark~\ref{rem:deg}. Note the same deterioration is observed in the power and degenerate benchmarks.

\paragraph{Shear-thickening regime.}

The relative performance of the methods is reversed for $r>2$. The primal Ka\v{c}anov iteration does not converge already at $r=3$, whereas Newton remains convergent up to $r=8$. Newton with $r$-continuation extends this range to the largest tested exponent, $r=10$.
The dual-weighted Ka\v{c}anov iteration on the KS discretization is also convergent throughout the tested shear-thickening range. Its iteration count increases from $16$ at $r=3$ to $29$, $56$, $85$, $121$ at $r=4$, $6$, $8$, $10$, respectively. As for the primal iteration, the same qualitative behaviour is observed across all three benchmarks.
The \emph{asymptotic contraction factors} $q_\mathrm{dual}$ are well described by
\begin{equation}\label{eq:emp-rel-dual}
q_{\mathrm{dual}}
\approx
\frac{r-2}{r-1},
\qquad r>2.
\end{equation}
Thus the dual iteration is linearly convergent, but becomes progressively less contractive as $r$ increases. Over the tested interval, this corresponds to an approximately linear growth of the iteration count with $r$.
These two empirical relations \eqref{eq:emp-rel-dual}, \eqref{eq:emp-rel-prim} provide a quantitative explanation for the complementary roles of the
primal and dual Ka\v{c}anov schemes in the shear-thinning and
shear-thickening regimes.

\paragraph{Finite element convergence rates.}

Table~\ref{tab:trig-rates-summary} reports the observed orders for representative values of $r$. The $\mathbb{P}_2-\mathbb{P}_1$ element exhibits  second-order convergence in both velocity and pressure.
Both MINI and KS display approximately first-order convergence for the velocity and pressure over the entire tested range.

\paragraph{Regularization and continuation.}

Regularization and continuation improve Newton's robustness through two different mechanisms. The continuation strategy leaves the target discrete problem unchanged and only provides a sequence of improved initial guesses. Accordingly, whenever both the direct and continued solvers converge, they produce essentially identical final errors. Regularization, instead, replaces the singular or degenerate viscosity by a shifted one and therefore modifies the discrete problem.
With continuation alone, Newton converges over the tested interval $1.5\leq r\leq10$. However, continuation increases the total iteration number when the direct solve is already successful.

The effect of $\delta_h$-regularization is particularly evident near the shear-thinning threshold. At $r=1.5$, direct Newton fails without regularization, whereas the choices
\(
k=1, k=10, k=100
\)
require $244$, $41$, and $14$ nonlinear iterations, respectively. The corresponding velocity errors are
\(
3.63\times10^{-5},
3.46\times10^{-4},
3.09\times10^{-3}
\).
Increasing $k$ therefore improves robustness but also increases the error.
Combining regularization and continuation extends the range further to
\(
1.2\leq r\leq10.
\)
The combination is therefore the most robust Newton-based strategy for this test.
Overall, the trigonometric benchmark identifies three trends that remain visible in the subsequent tests: the complementary roles of the primal and dual Ka\v{c}anov iterations, the independence of the attained discrete solution from the convergent nonlinear strategy, and the distinction between continuation as an initialization technique and regularization as a modification of the discrete operator. The power and degenerate benchmarks will show how these general mechanisms are altered by limited regularity and by an extended zero-strain set.

\begin{table}[H]
\centering
\caption{Observed and reference convergence rates for the trigonometric benchmark on the finest mesh ($n_h=128$).
The reference columns give the MINI velocity and pressure rates from Proposition~\ref{prop:literature-rates} for all three elements.
PN denotes primal Newton and PC denotes primal Ka\v{c}anov.}
\label{tab:trig-rates-summary}
\begin{tabular}{llccccl}
\toprule
FEM & $r$ & $\mathrm{rate}_u^{\,\mathrm{obs}}$ & $\mathrm{rate}_u^{\,\mathrm{ref}}$ & $\mathrm{rate}_p^{\,\mathrm{obs}}$ & $\mathrm{rate}_p^{\,\mathrm{ref}}$ & solver \\
\midrule
\multirow{4}{*}{$\mathbb{P}_2$--$\mathbb{P}_1$} & 1.5 & 2.03 & 1.00 & 1.64 & 0.67 & PC \\
 & 2.0 & 2.00 & 1.00 & 1.99 & 1.00 & PN \\
 & 3.0 & 1.93 & 0.75 & 2.00 & 0.75 & PN \\
 & 5.0 & 1.50 & 0.62 & 2.01 & 0.62 & PN \\
\midrule
\multirow{4}{*}{MINI} & 1.5 & 1.01 & 1.00 & 1.07 & 0.67 & PC \\
 & 2.0 & 1.01 & 1.00 & 1.05 & 1.00 & PN \\
 & 3.0 & 1.01 & 0.75 & 1.04 & 0.75 & PN \\
 & 5.0 & 1.01 & 0.62 & 1.04 & 0.62 & PN \\
\midrule
\multirow{4}{*}{KS} & 1.5 & 1.01 & 1.00 & 1.01 & 0.67 & PC \\
 & 2.0 & 1.01 & 1.00 & 1.01 & 1.00 & PN \\
 & 3.0 & 1.00 & 0.75 & 1.00 & 0.75 & PN \\
 & 5.0 & 0.97 & 0.62 & 0.98 & 0.62 & PN \\
\bottomrule
\end{tabular}
\end{table}

\subsection{Power benchmark}
\label{subsec:power-results}

\begin{figure}[H]
    \centering
    \begin{subfigure}[t]{0.45\textwidth}
        \centering
        \includegraphics[width=\linewidth]{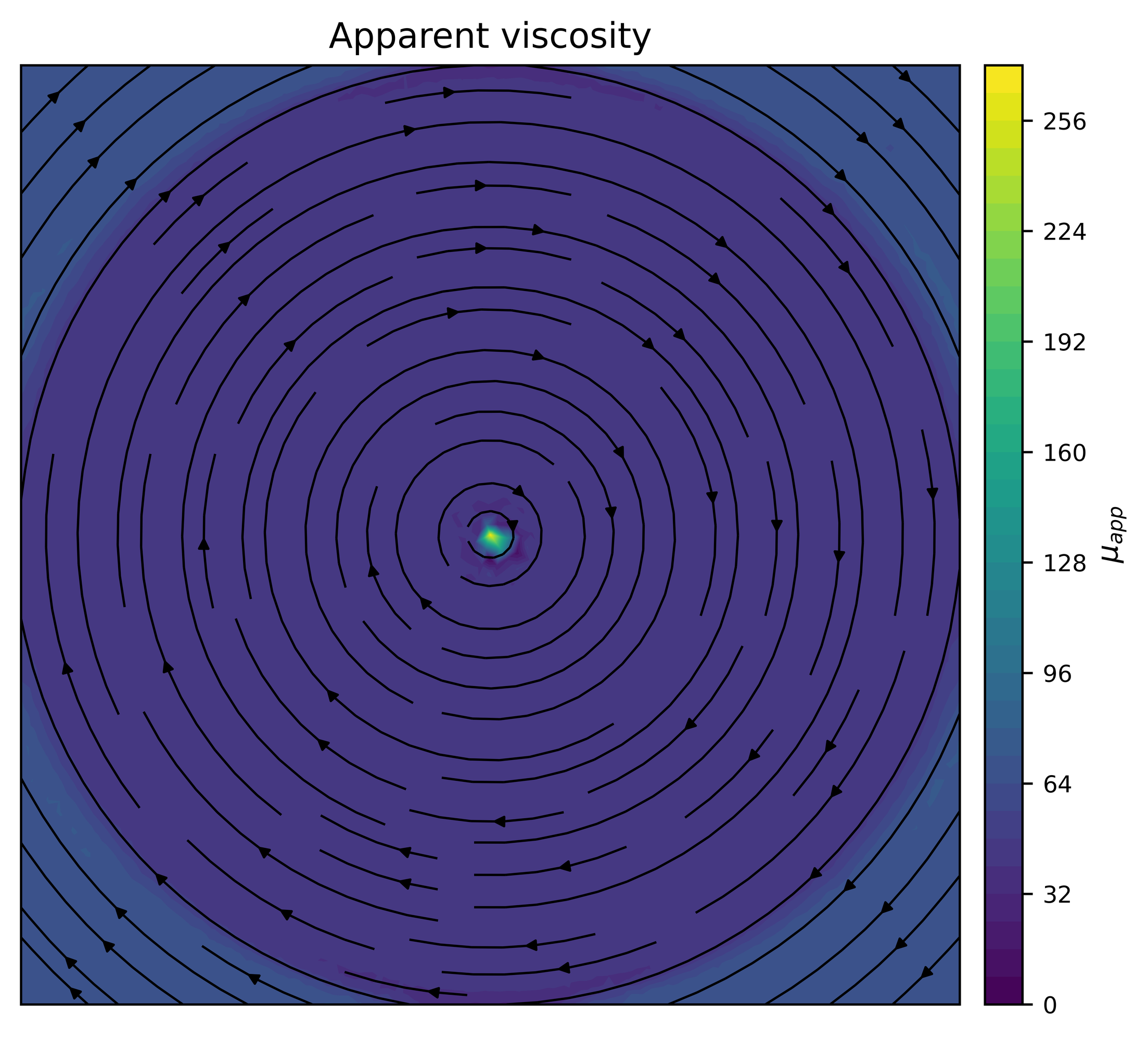}
        \caption{Apparent viscosity $\nu(\beps(\bu))$ and velocity streamlines at $r=1.01$.}
        \label{fig:power-viscosity}
    \end{subfigure}
    \hfill
    \begin{subfigure}[t]{0.48\textwidth}
        \centering
        \includegraphics[width=\linewidth, trim=0 0 0 0,
        clip]{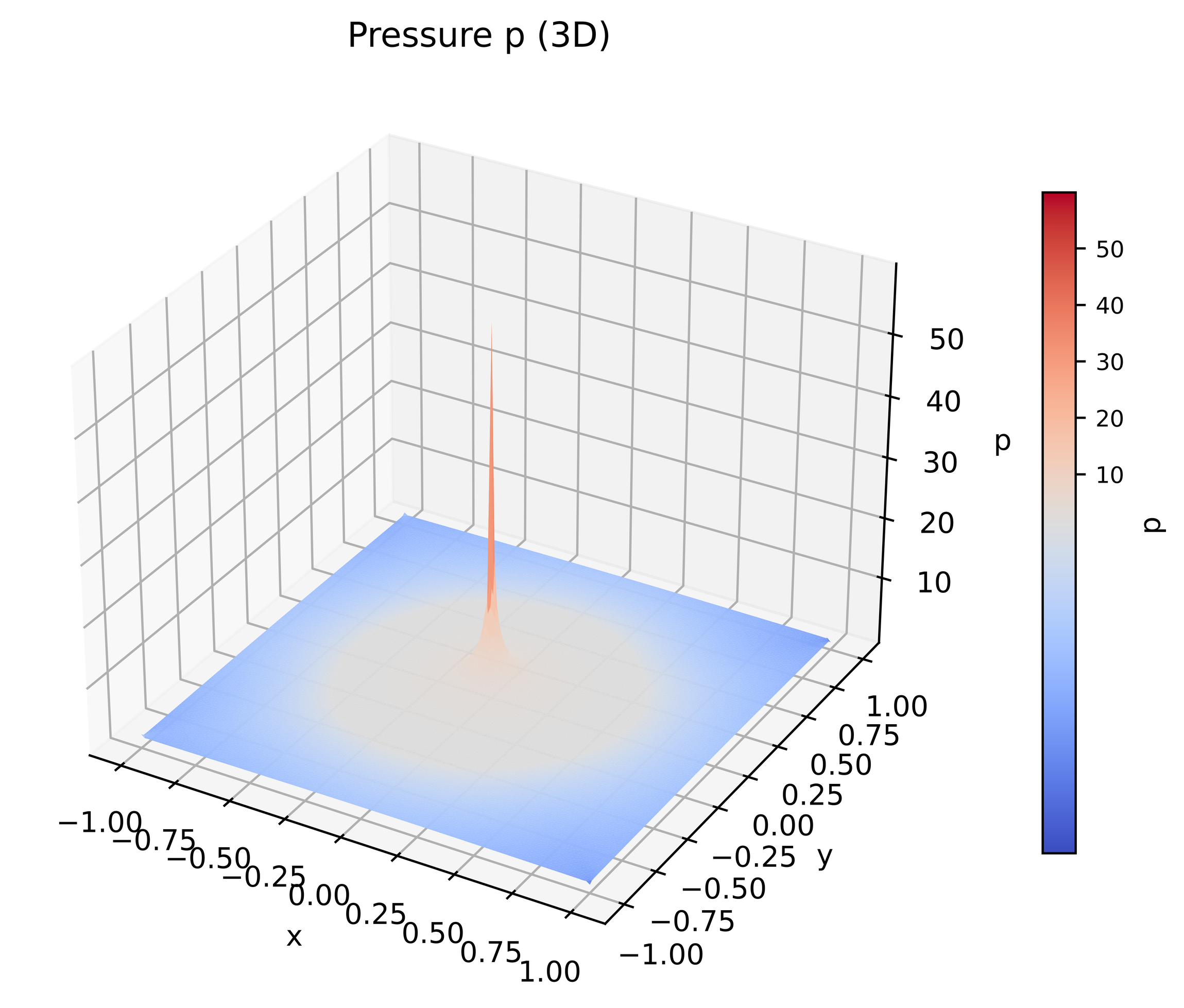}
        \caption{Pressure field $p$ at $r=7$.}
        \label{fig:power-pressure}
    \end{subfigure}
    \caption{Computed fields for the power benchmark.}
    \label{fig:power-fields}
\end{figure}

The power benchmark has limited regularity at the interior point $(0,0)$. It is therefore used to distinguish the effect of a localized singularity from that of the extended zero-strain set considered in Section~\ref{subsec:degenerate-benchmark}. As shown below, the shear-thickening regime is characterized by large velocity errors, which are mainly due to the singular pressure combined with the lack of pressure robustness of the considered elements, whereas the loss of solver robustness is substantially milder.
Figure~\ref{fig:power-viscosity} shows the apparent viscosity and velocity streamlines, and Figure~\ref{fig:power-pressure} shows the pressure field.

\begin{figure}[hbtp]
    \centering
    \includegraphics[width=0.98\linewidth]{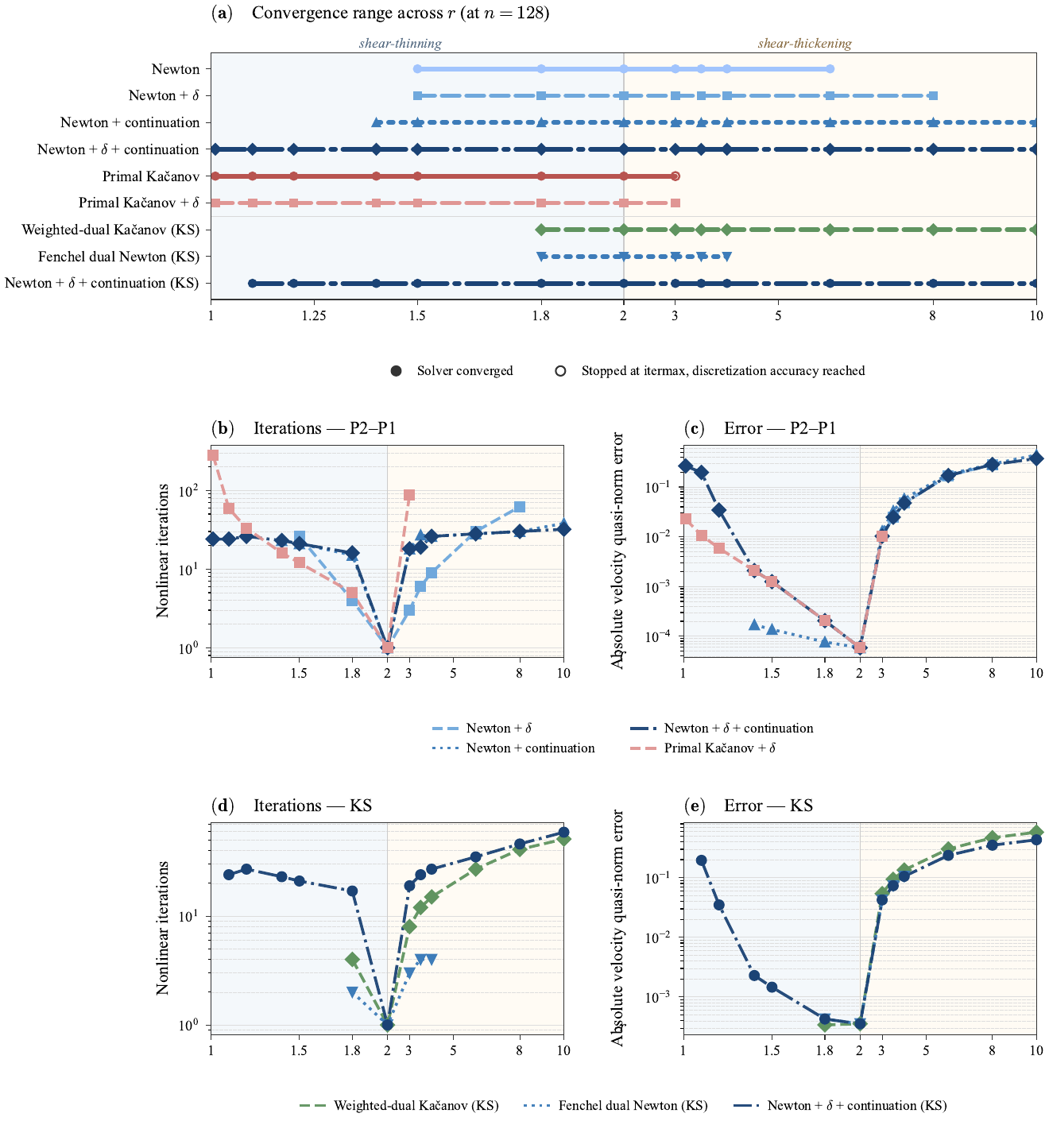}
    \caption{
    Summary of the nonlinear solver performance for the power benchmark on the finest mesh, with $n_h=128$. Panel~(a) reports the sampled exponents for which each strategy converges. Panels~(b)--(c) show, for the $\mathbb{P}_2$--$\mathbb{P}_1$ discretization, the nonlinear iteration counts and the velocity quasi-norm errors, respectively. Panels~(d)--(e) report the corresponding quantities for the KS discretization and the dual formulations. The blue and beige backgrounds distinguish the shear-thinning and shear-thickening regimes. Lines connect sampled values only and are included as visual guides.
    }
    \label{fig:power-summary}
\end{figure}

Figure~\ref{fig:power-summary} summarizes the convergence ranges, iteration counts, and velocity quasi-norm errors on the finest mesh, with $n_h=128$.

\paragraph{Primal formulations.}

Without regularization, Newton on the $\mathbb{P}_2$--$\mathbb{P}_1$ scheme and primal Ka\v{c}anov iterations preserve the same roles as in Section~\ref{subsec:trigonometric-results}. Continuation is particularly effective for $r>2$.
%in the shear-thickening regime, as in Section~\ref{subsec:trigonometric-results}.

\paragraph{Accuracy in the shear-thickening regime.}

For the non-regularized $\mathbb{P}_2$--$\mathbb{P}_1$
problem, the velocity quasi-norm error increases from
\(
6.0\times10^{-5}
~\text{at }r=2
\)
to
\(
1.32\times10^{-2},
5.81\times10^{-2},
2.00\times10^{-1},
3.38\times10^{-1},
4.51\times10^{-1}
\)
at $r=3$, $4$, $6$, $8$, and $10$, respectively.
Over the sampled shear-thickening values, the errors are well described, with respect to $r$, by the empirical relations
\[
\|\bu_{\mathrm{ex}}-\bu_h\|_{\mathrm{QN}}
\approx
1.56\times10^{-2}(r-2)^{1.71}
\qquad
\|p_{\mathrm{ex}}-p_h\|_{L^{r'}}
\approx
9.77\times10^{-4}(r-2)^{1.58}
\]

These errors are larger than the quasi-norm of the exact velocity,
$\|\bu_{\mathrm{ex}}\|_{\mathrm{QN}}\approx 8.4\cdot10^{-4}$, $7.0\cdot10^{-5}$ and $2.5\cdot10^{-11}$ at $r=3$, $4$ and $10$, respectively, i.e.\ than the error of the trivial approximation $\bv_h=\mathbf 0$. As in the degenerate benchmark (Section~\ref{subsec:degenerate-benchmark}), this indicates that the velocity error is dominated by the pressure-dependent consistency error, here amplified by the singular pressure $p_{\mathrm{ex}}\sim\rho^{\gamma}$ with $\gamma<0$ for $r>2$, rather than by the limited regularity of the velocity.

\paragraph{Regularization and continuation.}

Regularization also enlarges Newton's convergence range, but is less effective than continuation in the shear-thickening regime, as shown in Section~\ref{subsec:trigonometric-results}.
Combining regularization and continuation gives the widest Newton convergence range, covering all sampled exponents from $r=1.01$ to $r=10$, as in Section~\ref{subsec:trigonometric-results}. This robustness may, however, be accompanied by a large regularization error. At $r=1.01$, for example, the combined method with $k=100$ gives a velocity error of approximately $2.68\times10^{-1}$, whereas the non-regularized primal Ka\v{c}anov solution has an error of $7.79\times10^{-4}$.

\subsection{Degenerate benchmark}
\label{subsec:degenerate-benchmark}

Unlike the trigonometric benchmark, the degenerate test shows that robustness cannot be characterized solely in terms of the power-law exponent. The extended zero-strain set affects both rheological regimes. For $1<r<2$, the effective viscosity becomes singular as $\|\varepsilon(\bu)\|\to0$. For $r>2$, it vanishes on the cross-shaped zero-strain set, and the linearized constitutive operator loses uniform ellipticity.

\begin{figure}[hbtp]
    \centering
    \includegraphics[width=0.98\linewidth]{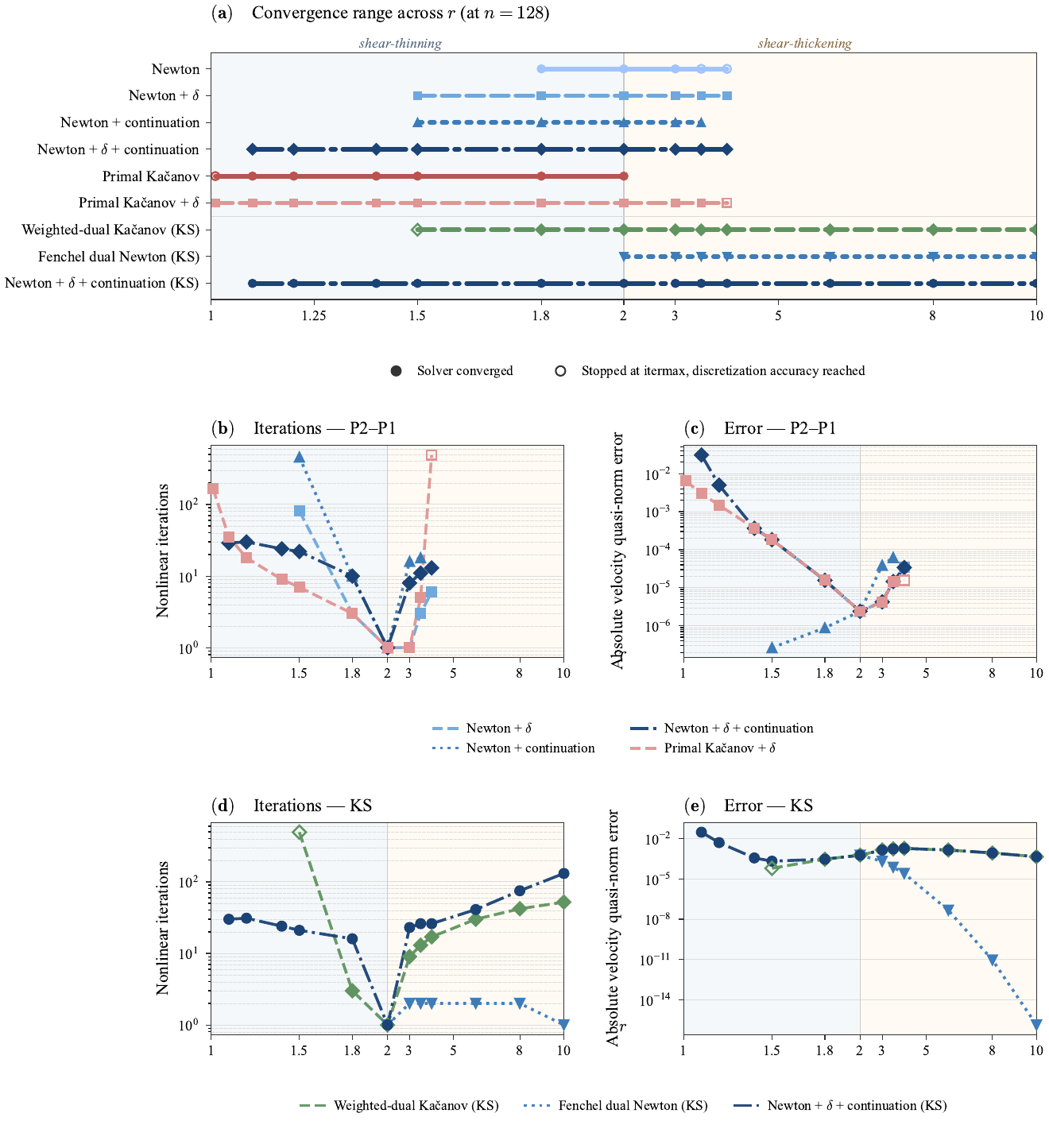}
        \caption{
    Summary of the nonlinear solver performance for the degenerate benchmark on the finest mesh, with $n_h=128$. Panel~(a) reports the sampled exponents for which each strategy converges. Panels~(b)--(c) show, for the $\mathbb{P}_2$--$\mathbb{P}_1$ discretization, the nonlinear iteration counts and the velocity quasi-norm errors, respectively. Panels~(d)--(e) report the corresponding quantities for the KS discretization and the dual formulations. The blue and beige backgrounds distinguish the shear-thinning and shear-thickening regimes. Lines connect sampled values only and are included as visual guides.
    For large $r$, the Fenchel-dual Newton iterates did not reach the discrete solution; see Section~\ref{subsec:degenerate-benchmark}.
    }
    \label{fig:degenerate-summary}
\end{figure}

The $\delta_h$-regularization acts structurally in both regimes: it bounds the viscosity at zero strain for $r<2$ and, for $r>2$, restores a positive ellipticity constant on each fixed mesh. This constant is of order $\delta_h^{r-2}$ and therefore still degenerates as $h\to0$.
Figure~\ref{fig:degenerate-summary} summarizes the sampled convergence ranges, iteration counts, and velocity quasi-norm errors on the finest mesh, with $n_h=128$.

\paragraph{Behaviour on the $\mathbb{P}_2$--$\mathbb{P}_1$ space.}

Without regularization, direct Newton converges only for the sampled exponents from $r=1.8$ to $r=3$. Continuation alone extends this range only modestly, reaching $r=1.5$ in the shear-thinning regime and $r=3.5$ in the shear-thickening regime. With regularization, for $k=100$, direct Newton converges for $r\in[1.5,4]$, while continuation extends convergence to $r=1.1$. At $r=1.5$, increasing $k$ from $1$ to $100$ reduces the iteration count from $81$ to $5$, at the cost of increasing the velocity quasi-norm error from $1.85\times10^{-4}$ to $1.94\times10^{-3}$.
Primal Ka\v{c}anov retains greater robustness on $\mathbb{P}_2$--$\mathbb{P}_1$ in the shear-thinning regime, as it converges down to $r=1.1$, while with regularization and $k=100$, convergence is restored at $r=1.01$.

\paragraph{Dependence on the finite element discretization.}

The loss of robustness observed for continued Newton is strongly discretization-dependent. On the $\mathbb{P}_2$--$\mathbb{P}_1$ space, the regularized continued method converges only up to $r=3.5$. On the KS space, the same regularized Newton continuation strategy remains convergent from $r=1.1$ to $r=10$.

\paragraph{Interaction between formulation and linearization on the KS discretization.}
The KS velocity discretization provides a common setting in which to compare the
primal, dual-weighted, and Fenchel-dual strategies. In the shear-thickening
regime, both the continued primal Newton strategy and the dual-weighted
Ka\v{c}anov iteration remain convergent up to the largest tested exponent.
The dual-weighted iteration requires fewer nonlinear iterations than the
continued primal Newton strategy shown in Figure~\ref{fig:degenerate-summary}.

\paragraph{Fenchel-dual Newton.}
As $r$ increases, the Fenchel-dual Newton iteration terminates after one or two steps, with
quasi-norm errors that are several orders of magnitude smaller than those of the dual-weighted
Ka\v{c}anov iteration, which targets the same non-regularized KS problem.
Since the discrete solution is unique, the two results cannot both correspond to it. The errors
of the dual-weighted iteration coincide with those of the independent, regularized
Newton--continuation strategy on the same space, which indicates that the Fenchel-dual Newton
iterates have not reached the discrete solution. We therefore do not draw conclusions on the
efficiency of Fenchel-dual Newton from this benchmark.

\paragraph{Accuracy at large $r$ and pressure robustness.}
For $r>2$ the velocity error of the degenerate benchmark does not follow the size of the exact solution. On the finest mesh ($n_h=128$), $\|\bu_{\mathrm{ex}}\|_{\mathrm{QN}}$ drops from $\approx 6.1\cdot10^{-5}$ at $r=2$ to $1.8\cdot10^{-8}$ at $r=4$ and $1.5\cdot10^{-31}$ at $r=10$, whereas the computed Kouhia--Stenberg errors remain of order $10^{-3}$ over $2\le r\le 10$, i.e.\ they exceed the error of the trivial approximation $\bv_h=\mathbf 0$. This cannot happen for a pressure-robust discretization: there the pressure does not enter the velocity error, which is bounded, up to a constant depending on $r$, by the best-approximation error in the quasi-norm~\cite{diening2025pressure}. Choosing $\bv_h=\mathbf 0$ as competitor then gives $\|\bu_{\mathrm{ex}}-\bu_h\|_{\mathrm{QN}}\le C(r)\,\|\bu_{\mathrm{ex}}\|_{\mathrm{QN}}$. None of the MINI, $\mathbb{P}_2$--$\mathbb{P}_1$ and Kouhia--Stenberg elements has this property, and their velocity error contains a consistency term that depends on the exact pressure but not on the size of the velocity~\cite{John2017divergence}.

In this benchmark the load is, for large $r$, essentially the gradient of $p_{\mathrm{ex}}$, so that this term is the only significant source of velocity error. Since all convergent solvers produce the same error, we attribute it to this pressure-dependent consistency error rather than to the linearization. The same mechanism dominates the shear-thickening regime of the power benchmark (Section~\ref{subsec:power-results}), whereas in the trigonometric benchmark the exact velocity is of order one and the effect is not visible.

\begin{figure}[htbp]
  \centering
  \includegraphics[width=\textwidth]{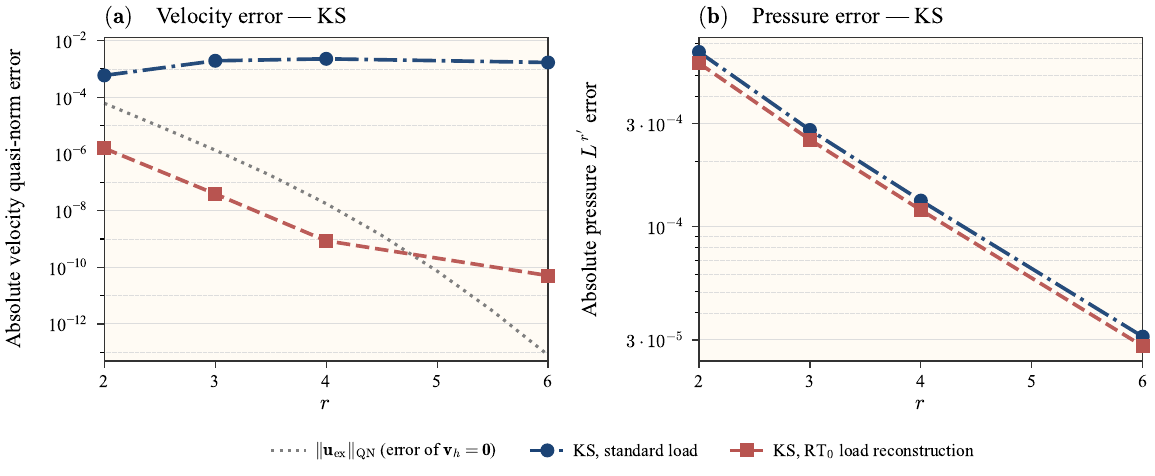}
  \caption{Degenerate benchmark on the finest mesh ($n_h=128$): Kouhia--Stenberg discretization with the standard load and with the load evaluated on the lowest-order Raviart--Thomas reconstruction of the test functions, both computed without regularization and with the same solver settings. (a) Absolute velocity quasi-norm error; the dotted grey line is $\|\bu_{\mathrm{ex}}\|_{\mathrm{QN}}$, i.e.\ the error of $\bv_h=\mathbf 0$. (b) Absolute pressure $L^{r'}$ error. Ka\v{c}anov is used at $r=2$ and primal Newton for $r>2$, continued until a further step changes the error by less than $0.1\%$; the case $r=6$ is reached by continuation in $r$ from $r=4$.}
  \label{fig:degenerate-pressure-robust}
\end{figure}

% For pressure-robust discretizations, we observe that, if the discretely divergence-free velocities are exactly divergence-free, the pressure drops out of the error equation and the velocity error is quasi-optimal in the natural distance~\cite{Belenki.Berselli.ea:12,diening2025pressure}. 
% None of the MINI, $\mathbb{P}_2$--$\mathbb{P}_1$ and Kouhia--Stenberg elements has this property, and their velocity error contains a consistency term that depends on the exact pressure but not on the size of the velocity~\cite{John2017divergence}. In this benchmark the load is, for large $r$, essentially the gradient of $p_{\mathrm{ex}}$, so that this term is the only significant source of velocity error. Since all convergent solvers produce the same error, we attribute it to this pressure-dependent consistency error rather than to the linearization. The same mechanism dominates the shear-thickening regime of the power benchmark (Section~\ref{subsec:power-results}), whereas in the trigonometric benchmark the exact velocity is of order one and the effect is not visible.
\paragraph{A first test with a pressure-robust variant.}
A first test in the same setting with a pressure-robust variant of the Kouhia--Stenberg element supports this interpretation. Following~\cite{Linke2014role}, only the load is modified: the test function in $\int_\Omega \ff\cdot\bv_h\,\mathrm dx$ is replaced by its lowest-order Raviart--Thomas interpolant, which preserves the edge means of the normal component and hence the discrete divergence, so that discretely divergence-free test functions are mapped to exactly divergence-free fields; the viscous and pressure terms are unchanged. As shown in Figure~\ref{fig:degenerate-pressure-robust}, the velocity error decreases by two to seven orders of magnitude (e.g.\ from $2.3\cdot10^{-3}$ to $8.6\cdot10^{-10}$ at $r=4$) and, up to $r=4$, remains below $\|\bu_{\mathrm{ex}}\|_{\mathrm{QN}}$, while the pressure error is essentially unchanged. At $r=6$ the computed error ($5.1\cdot10^{-11}$) lies above $\|\bu_{\mathrm{ex}}\|_{\mathrm{QN}}\approx 8\cdot10^{-14}$; this value is limited by the accuracy of the numerical integration rather than by the discretization, since raising the quadrature degree from $10$ to $20$ lowers it to $2.7\cdot10^{-11}$. A systematic comparison including pressure-robust elements is left to future work.

\begin{figure}[H]
    \centering
    \begin{subfigure}[t]{0.47\textwidth}
        \centering
        \includegraphics[width=\linewidth]
        {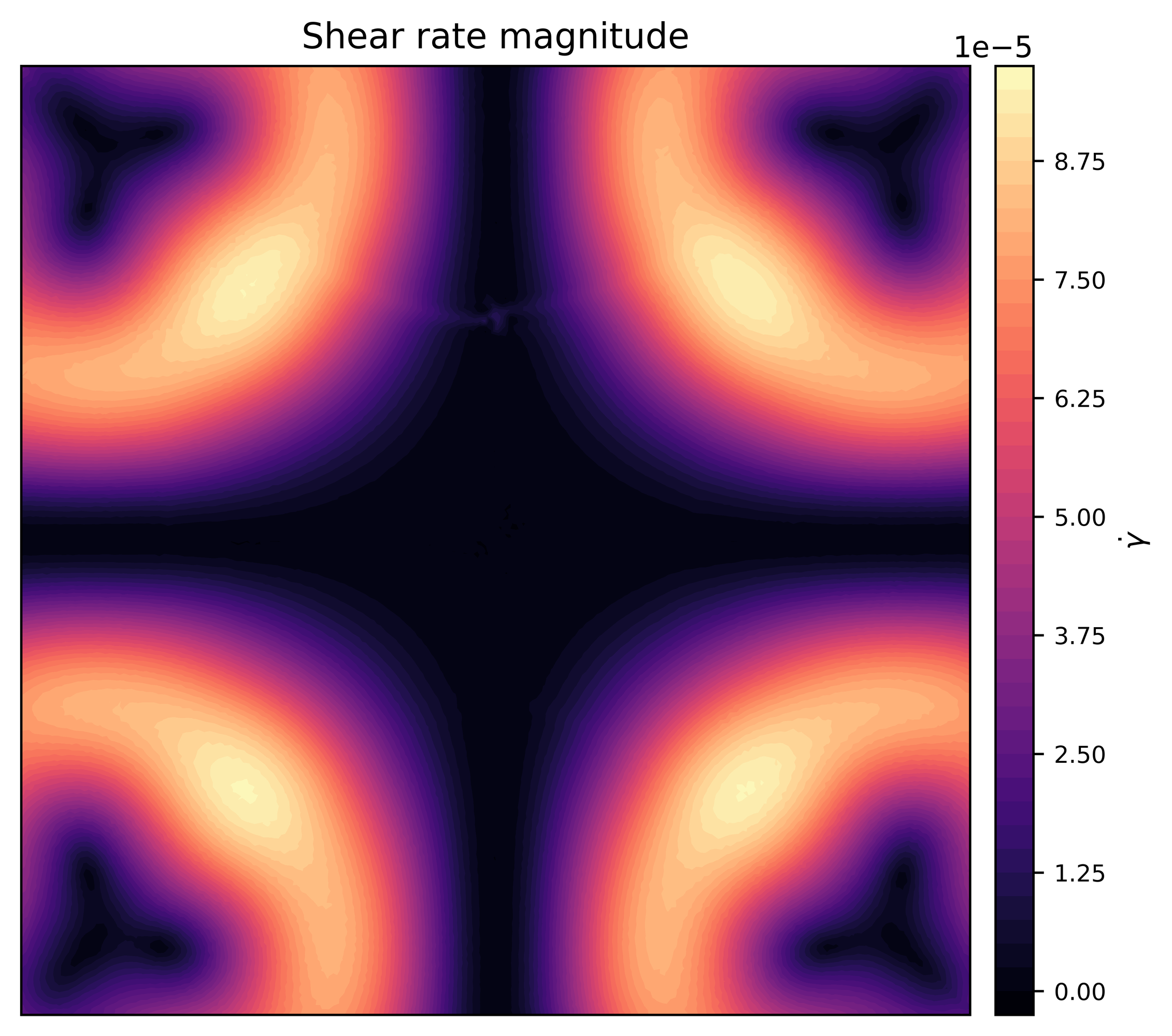}
        \caption{
        Shear-rate magnitude
        $\dot{\gamma}(\bu)
        =\sqrt{2\,\varepsilon(\bu):\varepsilon(\bu)}$.
        }
    \end{subfigure}
    \hfill
    \begin{subfigure}[t]{0.47\textwidth}
        \centering
        \includegraphics[width=\linewidth]
        {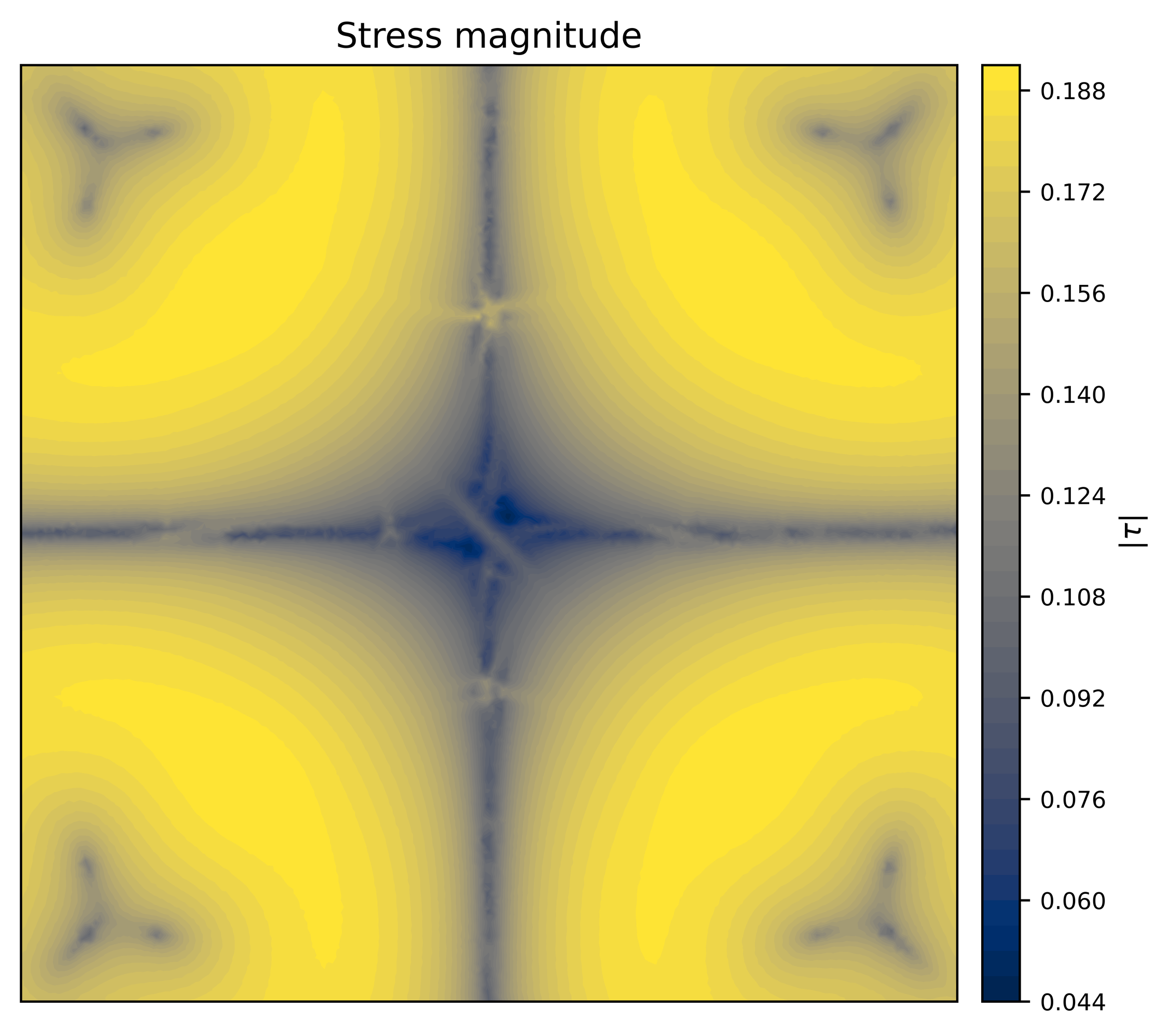}
        \caption{
        Deviatoric-stress magnitude
        $\|\bsigma\|$, with
        $\bsigma=\nu(\beps(\bu))\varepsilon(\bu)$.
        }
    \end{subfigure}

    \caption{
    Shear-rate and deviatoric-stress magnitudes for the degenerate
    benchmark at $r=1.1$, computed with the non-regularized primal
    Ka\v{c}anov method, the $\mathbb{P}_2$--$\mathbb{P}_1$
    discretization, and $n_h=128$.
    }
    \label{fig:degenerate-fields}
\end{figure}

Figure~\ref{fig:degenerate-fields} displays the shear-rate and deviatoric-stress magnitudes for the strongly shear-thinning case $r=1.1$. The cross-shaped zero-strain set is clearly visible in the shear-rate field. Although the effective viscosity is singular on this set, the stress remains finite, since
\(
\|\bsigma\|
\sim
\|\beps(\bu)\|^{r-1}.
\)
It therefore vanishes on the cross, although only slowly for $r=1.1$, where $r-1=0.1$.

\section{Glaciology test}\label{sec:glaciologytests}

\begin{figure}[ht!]
    \centering
    \includegraphics[width=0.7\linewidth]{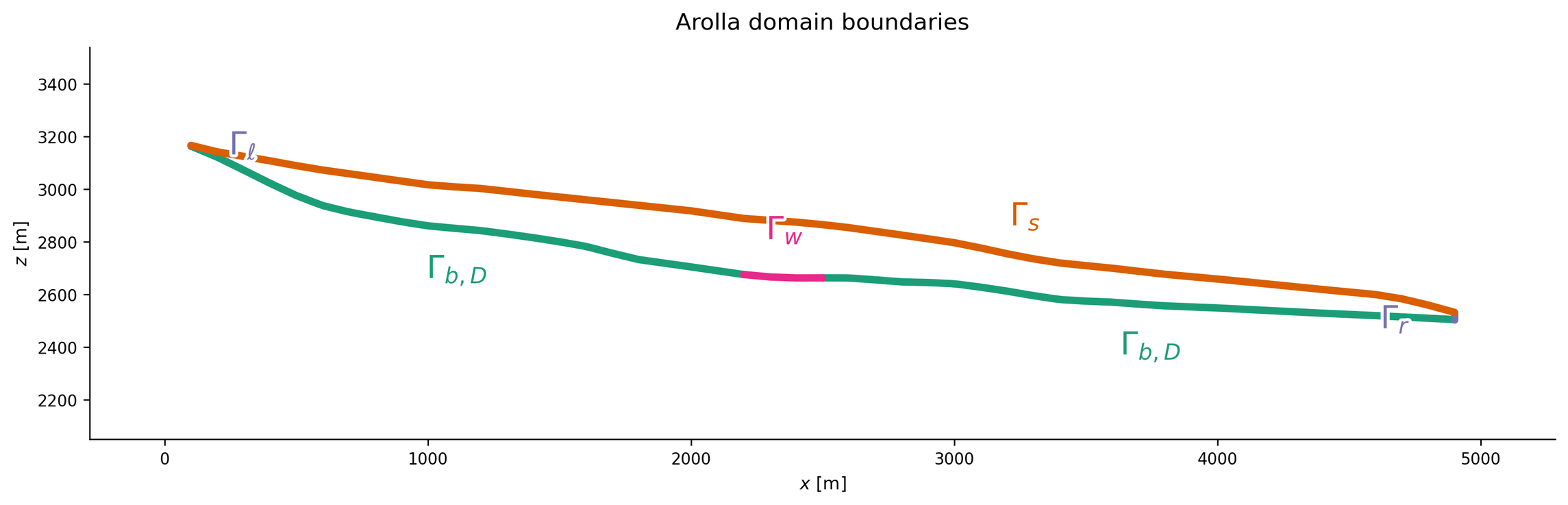}
    \caption{Boundary partition for the fitted Arolla mesh. The coordinates are shown in metres.}
    \label{fig:glaciology-arolla-bed-slip-bc}
\end{figure}

We use one glaciology experiment to check that the velocity--pressure formulation in Section~\ref{sec:continuouspb} also handles a fitted glacier flowline mesh. The mesh and boundary partition are the ones shown in Figure~\ref{fig:glaciology-arolla-bed-slip-bc}.
The test is based on the Haut Glacier d'Arolla profile \cite{Gagliardini2008,Pattyn2008} and is motivated by
\cite{Ahlkrona2021}. In contrast to the CutFEM discretization used there, our computation uses a fitted triangular mesh: the bed and free-surface elevations are read from the Arolla profile, zero-thickness endpoint columns are discarded, and each vertical column is interpolated between bed and surface before being split into triangles. The mesh is made by \(n_x=37\) points along the flowline and \(n_z=14\) points through the thickness. Also note, compared to \cite{Ahlkrona2021}, that the finite element spaces and stabilization mechanisms are different: in our computations we use the finite element discretizations described in Section~\ref{sec:discretepb}, without ghost penalties or CutFEM pressure stabilization.

\subsection{Geometry, boundary conditions, weak formulation}
\label{subsec:arolla-bed-slip-geometry}

We denote the fitted Arolla flowline domain by \(\Omega_{\mathrm{A}}\). Its
boundary is decomposed as
\[
    \partial\Omega_{\mathrm{A}}
    =
    \Gamma_s
    \cup
    \Gamma_{b,D}
    \cup
    \Gamma_w
    \cup
    \Gamma_{\ell}
    \cup
    \Gamma_r,
\]
where \(\Gamma_s\) is the upper free surface, \(\Gamma_{\ell}\) and
\(\Gamma_r\) are the left and right vertical sections,
\(\Gamma_{b,D}\) is the no-slip part of the bed, and \(\Gamma_w\) is a basal sliding patch, where no penetration is imposed weakly. The full bed is therefore
\[
    \Gamma_b = \Gamma_{b,D}\cup\Gamma_w,
    \qquad
    \text{with}~
    \Gamma_w
    :=
    \left\{
        x\in\Gamma_b
        :
        2200\,\mathrm{m}\le x_1\le 2500\,\mathrm{m}
    \right\}.
\]
We use the same notation as in the model problem. The outward unit normal is \(\bn\), and on \(\Gamma_w\) we set
\(
    \boldsymbol{t} := (n_2,-n_1).
\)
The body force is the gravity load
\(
    \ff = \rho\,\boldsymbol{g}_{\mathrm{ice}},
    \rho = 9.1380\cdot 10^{-19},
    \boldsymbol{g}_{\mathrm{ice}}
    =
    (0,-9.7692\cdot 10^{15})^T .
\)
For Glen's exponent \(n_G=3\), the corresponding power-law index is
\(
    r = 1+\frac{1}{n_G}=\frac{4}{3}
\).

\subsubsection{Strong formulation}

The equations are the nonlinear Stokes equations from
\eqref{eq:problem}, restricted to \(\Omega_{\mathrm{A}}\). The boundary conditions represented by the mesh partition are
\[
\begin{cases}
\bu = \boldsymbol{0}
& \text{on } \Gamma_{b,D}\cup\Gamma_{\ell}\cup\Gamma_r,\\
\bsigma(\beps(\bu))\bn - p\bn = \boldsymbol{0}
& \text{on } \Gamma_s,\\
\bu\cdot\bn = 0
& \text{on } \Gamma_w,\\
\bigl(\bsigma(\beps(\bu))\bn\bigr)\cdot\boldsymbol{t}
+ \beta_{\mathrm{bc}}(\bu\cdot\boldsymbol{t}) = 0
& \text{on } \Gamma_w.
\end{cases}
\]
The first line is imposed strongly by Dirichlet conditions. The free surface $\Gamma_s$ uses the homogeneous Neumann condition \(\bg=\boldsymbol{0}\) in the notation of \eqref{eq:weak_formulation}. On \(\Gamma_w\), the normal condition \(\bu\cdot\bn=0\) is enforced weakly using an \(r\)-power penalty with numerical parameter \(\alpha_{\mathrm{bc}}=9.1\).
The tangential condition is a Robin-type law whenever the coefficient \(\beta_{\mathrm{bc}}\) is different from zero. 
In what follows, we take \(\beta_{\mathrm{bc}}=0\), hence the normal penalty remains active, whereas the tangential Robin contribution vanishes.
%, giving free tangential traction. 
Note that \(\Gamma_{b,D}\) represents the portion of the glacier bed where the ice is frozen to the underlying bedrock, whereas \(\Gamma_w\) denotes the region where the presence of basal water allows the ice to slide over the bedrock.

\subsubsection{Weak formulation}

Let
\(
    \Gamma_D
    :=
    \Gamma_{b,D}\cup\Gamma_{\ell}\cup\Gamma_r.
\)
The continuous trial and test spaces are the same as in Subsection~\ref{subsubsec:primal-weak-form}, with $\Gamma_D$ as the Dirichlet set. The corresponding discrete spaces are denoted by  \(\boldsymbol{V}_{h}^{\mathrm{A}}\) and \(Q_h^{\mathrm{A}}\).
On the left-hand side of the discrete problem \eqref{eq:FEM-discrete-abstract} we add the operator
\begin{equation}
\begin{aligned}
\mathcal{B}_w(\bu_h,p_h;\bv_h)
:=
\int_{\Gamma_w}
\bigg[
&p_h(\bv_h\cdot\bn)
-
\bigl(
    \nu(\beps(\bu_h))\beps(\bu_h)\bn\cdot\bn
\bigr)
(\bv_h\cdot\bn)
\\
&+
\alpha_{\mathrm{bc}}h^{1-r}
\left|\bu_h\cdot\bn\right|^{r-2}
(\bu_h\cdot\bn)(\bv_h\cdot\bn)
+
\beta_{\mathrm{bc}}
(\bu_h\cdot\boldsymbol{t})
(\bv_h\cdot\boldsymbol{t})
\bigg]\,\mathrm{d}s.
\end{aligned}
\end{equation}
The first two terms are the normal consistency terms, the third term is the \(r\)-power penalty for \(\bu_h\cdot\bn=0\), and the last term is the tangential Robin contribution. Since \(\beta_{\mathrm{bc}}=0\) for this mesh, the final term is absent in the
actual Arolla bed-slip run.
The discrete problem is: find
\((\bu_h,p_h)\in \boldsymbol{V}_{h}^{\mathrm{A}}\times Q_h^{\mathrm{A}}\) such that
\begin{equation}
\begin{aligned}
&\int_{\Omega_{\mathrm{A}}}
\nu(\beps(\bu_h))\beps(\bu_h):\beps(\bv_h)\,\mathrm{d}x
-
\int_{\Omega_{\mathrm{A}}}
p_h\nabla\!\cdot\bv_h\,\mathrm{d}x
-
\int_{\Omega_{\mathrm{A}}}
q_h\nabla\!\cdot\bu_h\,\mathrm{d}x
\\
&\qquad
+
\mathcal{B}_w(\bu_h,p_h;\bv_h)
=
\int_{\Omega_{\mathrm{A}}}
\rho\boldsymbol{g}_{\mathrm{ice}}\cdot\bv_h\,\mathrm{d}x
\qquad\forall
(\bv_h,q_h)\in \boldsymbol{V}_{h}^{\mathrm{A}}\times Q_h^{\mathrm{A}}.
\end{aligned}
\end{equation}
%The stress-free surface contributes no nonzero boundary integral because its traction is homogeneous.

\subsection{Test results}

\begin{figure}[H]
    \centering

    \begin{subfigure}[t]{0.32\textwidth}
        \centering
        \includegraphics[width=\linewidth]{
            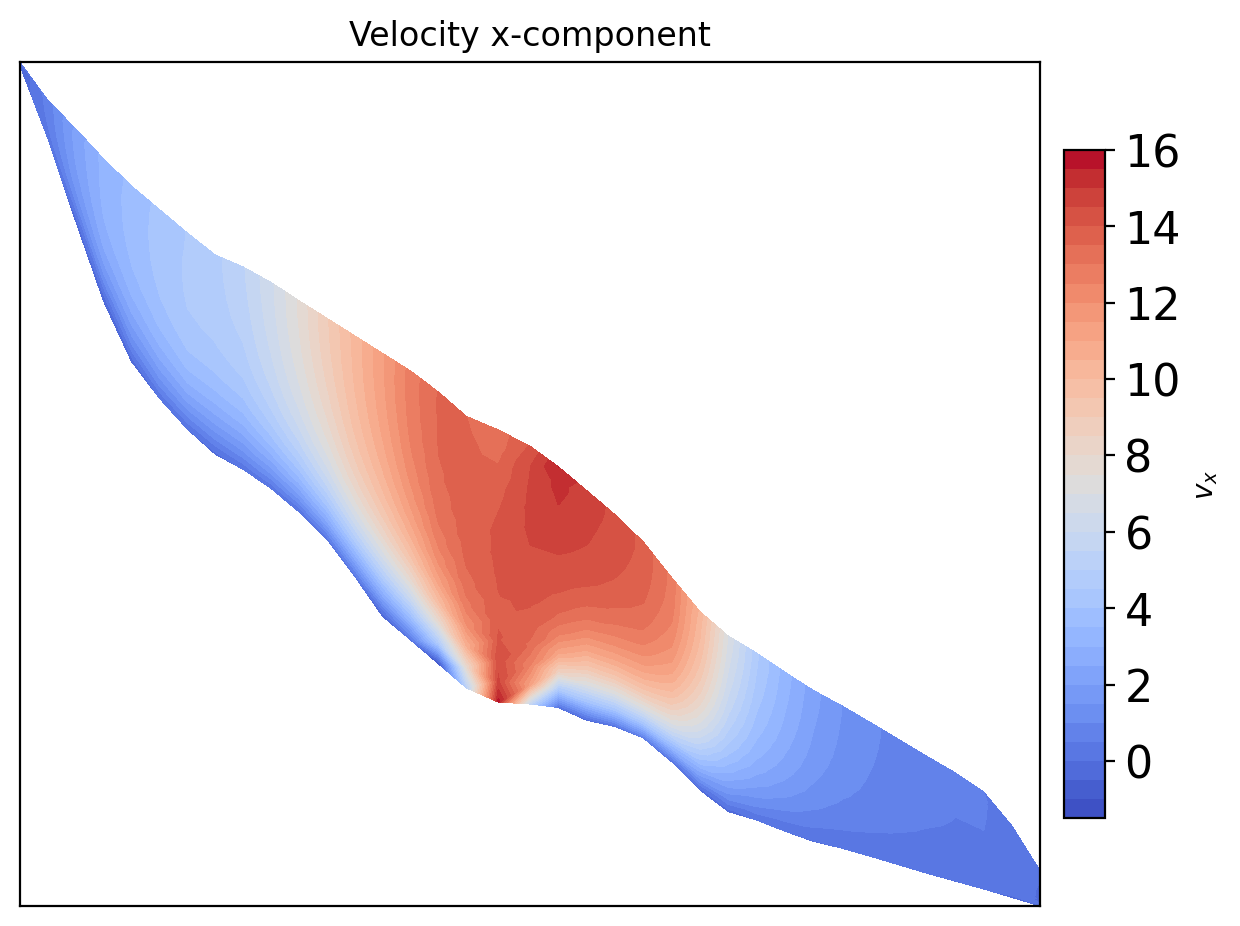
        }
        \caption{Horizontal velocity component \(u_x\).}
    \end{subfigure}
    \hfill
    \begin{subfigure}[t]{0.32\textwidth}
        \centering
        \includegraphics[width=\linewidth]{
            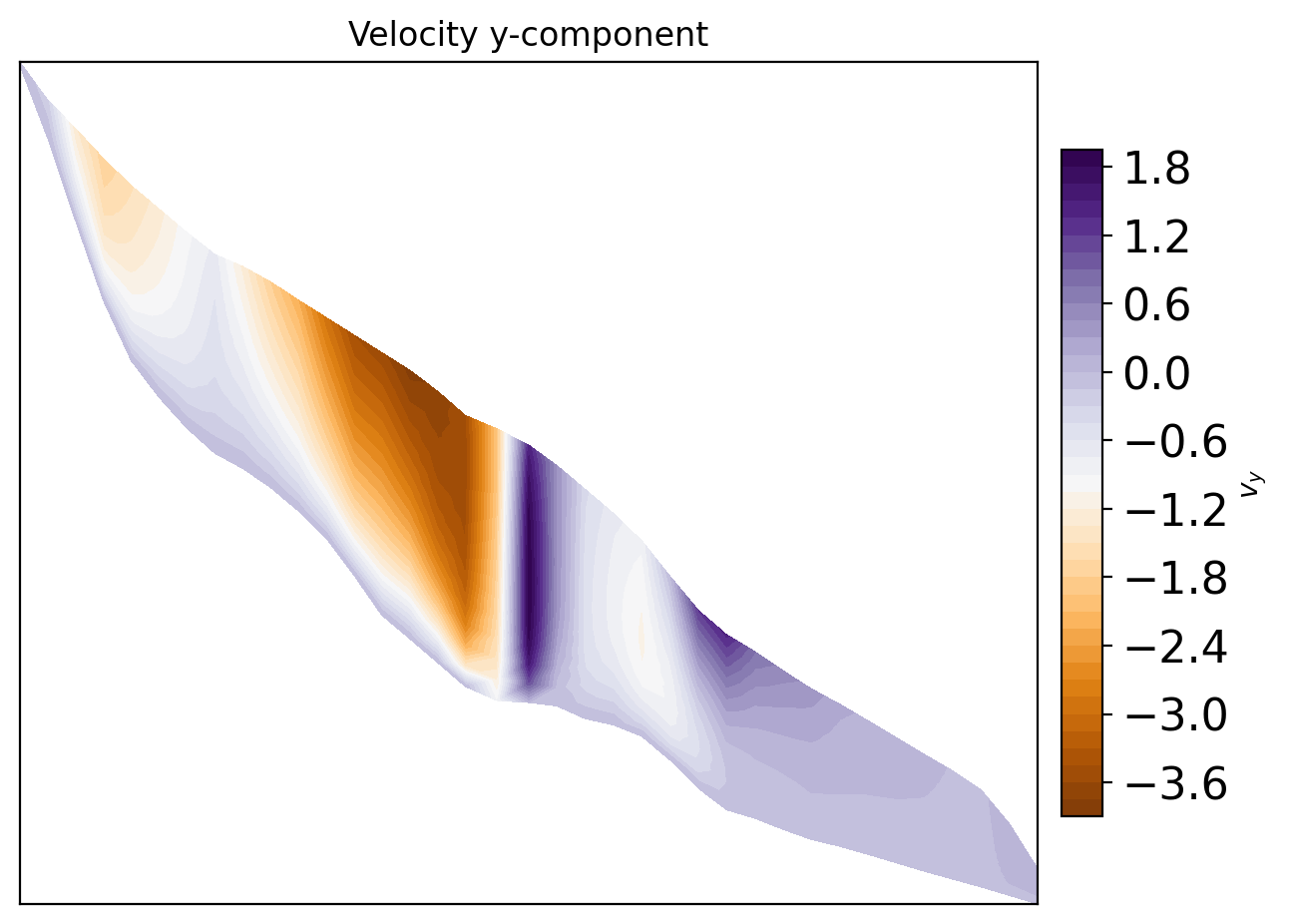
        }
        \caption{Vertical velocity component \(u_y\).}
    \end{subfigure}
    \hfill
    \begin{subfigure}[t]{0.32\textwidth}
        \centering
        \includegraphics[width=\linewidth]{
            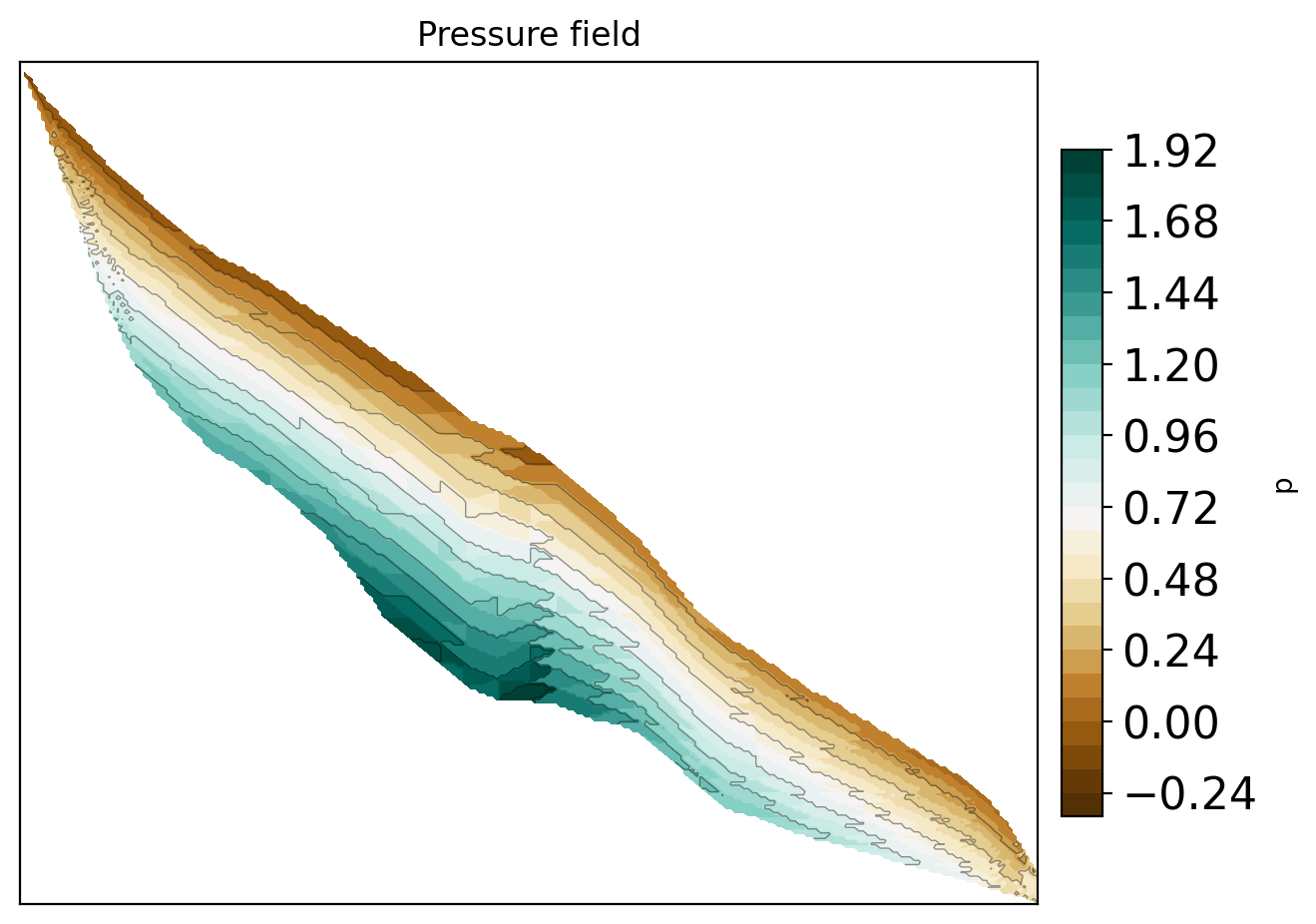
        }
        \caption{Pressure field \(p\).}
    \end{subfigure}

    \caption{
        Computed horizontal velocity \(u_x\), vertical velocity \(u_y\), and pressure \(p\) for the Arolla problem using the primal Kačanov method with \(\delta_h=0\), $n_x=37, n_z=14$, KS, $r=4/3$. The vertical dimension is exaggerated in all panels for ease of visualization.
    }
    \label{fig:arolla-fields}
\end{figure}

In Figure~\ref{fig:arolla-fields} we report the results of the setting mentioned in Section~\ref{subsec:arolla-bed-slip-geometry}. We report the plots of the velocity field, for $x$ and $y$ components, and of the pressure field.

\section{Conclusions}
\label{sec:conclusions}

This work provides a numerical investigation of finite element discretizations and nonlinear solvers for incompressible power-law Stokes flows, covering both shear-thinning and shear-thickening regimes. We assess solver robustness with benchmark configurations, including the power and degenerate tests, designed to isolate two distinct sources of numerical difficulty.

\paragraph{Finite element discretizations.}

On the smooth trigonometric benchmark, the $\mathbb{P}_2$--$\mathbb{P}_1$ element exhibits approximately second-order convergence in velocity and pressure, whereas MINI and Kouhia--Stenberg show a first-order behaviour. The choice of discretization methods covers the full range of exponents and avoids the need to introduce additional stabilization terms. The experiments also show that the finite element discretization can affect the convergence property of the linearization (cf. Figure~\ref{fig:degenerate-summary}). None of the elements considered is pressure robust. In the degenerate benchmark, where the load is dominated by its gradient part, and in the shear-thickening regime of the power benchmark, where the pressure is singular, the resulting pressure-induced error dominates and exceeds the error of the trivial approximation $\bv_h=\mathbf 0$. A first test with a Raviart--Thomas reconstruction of the load for the Kouhia--Stenberg element reduces the velocity error in the degenerate benchmark by up to seven orders of magnitude.

\paragraph{Primal Newton and Ka\v{c}anov iterations.}

We achieve convergence in both strongly nonlinear regimes: down to $r=1.01$ in the shear-thinning case using the primal Ka\v{c}anov iteration, and up to $r=10$ in the shear-thickening case using primal Newton combined with $r$-continuation. Newton's method is more sensitive to the choice of initial guess and finite element space; however, when direct convergence is achieved, it is often more efficient than Ka\v{c}anov. Regularization mitigates the singular behavior of the problem and reduces the number of nonlinear iterations required for convergence, at the price of a regularization error that grows with the regularization parameter, whereas continuation improves robustness without modifying the discrete problem.

\paragraph{Dual strategies.}
We have derived pressure-explicit representations of the primal and dual energy functionals associated with the Fenchel formulation of the problem. Unlike the standard constrained energy formulation, these representations do not require the divergence-free condition to be incorporated directly into the admissible velocity space and naturally accommodate mixed boundary conditions. On the Kouhia--Stenberg discretization, the dual-weighted Ka\v{c}anov iteration and the Ka\v{c}anov iteration applied to the two-field Fenchel-dual formulation generate identical iterates (Proposition~\ref{prop:dual-kacanov-equivalence}). The dual-weighted Ka\v{c}anov iteration converges over the whole sampled shear-thickening range in all three benchmarks, with contraction factors well described by $(r-2)/(r-1)$, mirroring the behaviour $2-r$ of the primal Ka\v{c}anov iteration in the shear-thinning regime. The primal and dual Ka\v{c}anov schemes therefore play complementary roles in the two rheological regimes. For Fenchel-dual Newton, the iterates obtained at large $r$ in the degenerate benchmark did not reach the discrete solution, and its performance in this regime requires further investigation.
Overall, the experiments reveal that no single formulation--solver combination is uniformly superior across all benchmarks.

\paragraph{Future work.}
A natural extension of this study concerns pressure-robust discretizations. Future work will include pressure-robust elements in the systematic comparison of finite element methods and linearization strategies, for instance reconstruction-based variants of the present elements~\cite{Linke2014role,lederer2017divergence,diening2025pressure}, Scott--Vogelius elements on barycentric meshes, or $H(\mathrm{div})$-conforming discontinuous Galerkin methods. The first test reported in Section~\ref{subsec:degenerate-benchmark} suggests that this would separate the effect of the discretization from that of the nonlinear solver more clearly, in particular for large $r$, and it is also relevant for gravity-driven applications such as glacier flow.
A second direction concerns the a priori analysis of higher-order discretizations. The reference
rates in Table~\ref{tab:trig-rates-summary} are those available for low-order elements, and the
$\mathbb{P}_2$--$\mathbb{P}_1$ element clearly exceeds them on smooth solutions, with a reduction of
the observed velocity rate for larger $r$. Deriving error estimates for higher-order finite element discretizations under
suitable regularity assumptions would provide a theoretical reference for these observations.

\section*{Acknowledgments}
MB kindly acknowledges partial financial support by INdAM-GNCS project 2026 CUP E53C25002010001. This research is part of the activities of “Dipartimento di Eccellenza 2023-2027”, funded by MUR, Italy.

\section*{Declaration and code availability}

%While open-source packages exist for related non-Newtonian and glaciological Stokes models, the finite element literature on power-law/p-Stokes discretizations only sporadically provides public, paper-specific implementations of the numerical methods and benchmark studies. Therefore, 
We release our code in a public \textsc{GitHub} repository \href{https://github.com/renonesopocaj/p-stokes}{https://github.com/renonesopocaj/p-stokes}, with the aim of making the experiments reproducible and easy to extend. The repository will be completed and made fully available alongside the final version of the paper.
During the development of the \textsc{Python} repository, the authors used \textsc{OpenAI ChatGPT 5} to assist with generating and debugging portions of code. AI-assisted code was reviewed, tested, and, where necessary, modified by the authors. The authors take responsibility for the correctness and integrity of the final code and reported results.

\bibliographystyle{abbrv}
\bibliography{references}

@article{Botti-Mascotto:2026,
  author  = {Botti, M. and Mascotto, L.},
  title   = {{S}obolev--{P}oincar\'e inequalities for piecewise {$W^{1,p}$} functions over general polytopic meshes},
  journal = {SIAM J. Numer. Anal.},
  volume  = {64},
  number  = {4},
  pages   = {1478--1504},
  year    = {2026},
  doi     = {10.1137/25M1748470}
}

@book{Galdi:2011,
  author    = {Galdi, G. P.},
  title     = {An Introduction to the Mathematical Theory of the {N}avier--{S}tokes Equations: Steady-State Problems},
  series    = {Springer Monographs in Mathematics},
  edition   = {Second},
  publisher = {Springer},
  address   = {New York},
  year      = {2011},
  doi       = {10.1007/978-0-387-09620-9}
}

@article{Barrett.Liu:94,
  author  = {Barrett, J. W. and Liu, W. B.},
  title   = {Quasi-norm error bounds for the finite element approximation of a non-{N}ewtonian flow},
  journal = {Numer. Math.},
  volume  = {68},
  number  = {4},
  pages   = {437--456},
  year    = {1994},
  doi     = {10.1007/s002110050071}
}

@article{BL:1993,
  author  = {Barrett, J. W. and Liu, W. B.},
  title   = {Finite element error analysis of a quasi-{N}ewtonian flow obeying the {C}arreau or power law},
  journal = {Numer. Math.},
  volume  = {64},
  number  = {4},
  pages   = {433--453},
  year    = {1993},
  doi     = {10.1007/BF01388698}
}

@article{Belenki.Berselli.ea:12,
  author  = {Belenki, L. and Berselli, L. C. and Diening, L. and R{\r{u}}{\v{z}}i{\v{c}}ka, M.},
  title   = {On the finite element approximation of {$p$}-{S}tokes systems},
  journal = {SIAM J. Numer. Anal.},
  volume  = {50},
  number  = {2},
  pages   = {373--397},
  year    = {2012},
  doi     = {10.1137/10080436X}
}

@article{Hirn:13,
  author  = {Hirn, A.},
  title   = {Approximation of the {$p$}-{S}tokes equations with equal-order finite elements},
  journal = {J. Math. Fluid Mech.},
  volume  = {15},
  number  = {1},
  pages   = {65--88},
  year    = {2013},
  doi     = {10.1007/s00021-012-0095-0}
}

@article{Kalte.Ruzi:23,
  author  = {Kaltenbach, A. and R{\r{u}}{\v{z}}i{\v{c}}ka, M.},
  title   = {A local discontinuous {G}alerkin approximation for the {$p$}-{N}avier--{S}tokes system, {P}art {I}: {C}onvergence analysis},
  journal = {SIAM J. Numer. Anal.},
  volume  = {61},
  number  = {4},
  pages   = {1613--1640},
  year    = {2023},
  doi     = {10.1137/22M151474X}
}

@article{Kreuzer.Suli:16,
  author  = {Kreuzer, C. and S{\"u}li, E.},
  title   = {Adaptive finite element approximation of steady flows of incompressible fluids with implicit power-law-like rheology},
  journal = {ESAIM Math. Model. Numer. Anal.},
  volume  = {50},
  number  = {5},
  pages   = {1333--1369},
  year    = {2016},
  doi     = {10.1051/m2an/2015085}
}

@article{Sandri:93,
  author  = {Sandri, D.},
  title   = {Sur l'approximation num\'erique des \'ecoulements quasi-newtoniens dont la viscosit\'e suit la loi puissance ou la loi de {C}arreau},
  journal = {RAIRO Mod\'el. Math. Anal. Num\'er.},
  volume  = {27},
  number  = {2},
  pages   = {131--155},
  year    = {1993},
  doi     = {10.1051/m2an/1993270201311}
}

@article{BN:1990,
  author  = {Baranger, J. and Najib, K.},
  title   = {Analyse num\'erique des \'ecoulements quasi-newtoniens dont la viscosit\'e ob\'eit \`a la loi puissance ou la loi de {C}arreau},
  journal = {Numer. Math.},
  volume  = {58},
  number  = {1},
  pages   = {35--49},
  year    = {1990},
  doi     = {10.1007/BF01385609}
}

@article{yasuda1981shear,
  author  = {Yasuda, K. and Armstrong, R. C. and Cohen, R. E.},
  title   = {Shear flow properties of concentrated solutions of linear and star branched polystyrenes},
  journal = {Rheol. Acta},
  volume  = {20},
  number  = {2},
  pages   = {163--178},
  year    = {1981},
  doi     = {10.1007/BF01513059}
}

@article{diening2025guaranteed,
  author  = {Diening, L. and Storn, J.},
  title   = {Guaranteed upper bounds for iteration errors and modified {K}a{\v{c}}anov schemes via discrete duality},
  journal = {Comput. Methods Appl. Math.},
  volume  = {25},
  number  = {3},
  pages   = {587--600},
  year    = {2025},
  doi     = {10.1515/cmam-2025-0017}
}

@article{diening2025pressure,
  author       = {Diening, L. and Hirn, A. and Kreuzer, C. and Zanotti, P.},
  title        = {Pressure robust finite element discretizations of the nonlinear {S}tokes equations},
  journal      = {Math. Models Methods Appl. Sci.},
  volume       = {35},
  number       = {12},
  pages        = {2661--2693},
  doi          = {10.1142/S0218202525500484},
  year         = {2025}
}

@article{MUMPS:1,
  author  = {Amestoy, P. R. and Duff, I. S. and L'Excellent, J.-Y. and Koster, J.},
  title   = {A fully asynchronous multifrontal solver using distributed dynamic scheduling},
  journal = {SIAM J. Matrix Anal. Appl.},
  volume  = {23},
  number  = {1},
  pages   = {15--41},
  year    = {2001},
  doi     = {10.1137/S0895479899358194}
}

@article{MUMPS:2,
  author  = {Amestoy, P. R. and Buttari, A. and L'Excellent, J.-Y. and Mary, T.},
  title   = {Performance and scalability of the block low-rank multifrontal factorization on multicore architectures},
  journal = {ACM Trans. Math. Softw.},
  volume  = {45},
  number  = {1},
  pages   = {2:1--2:26},
  year    = {2019},
  doi     = {10.1145/3242094}
}

@article{Seehanam2023,
  author  = {Seehanam, S. and Chanchareon, W. and Promoppatum, P.},
  title   = {Assessing the effect of manufacturing defects and non-{N}ewtonian blood model on flow behaviors of additively manufactured {G}yroid {TPMS} structures},
  journal = {Heliyon},
  volume  = {9},
  number  = {5},
  pages   = {Paper No. e15711},
  year    = {2023},
  doi     = {10.1016/j.heliyon.2023.e15711}
}

@article{WalburnSchneck1976,
  author  = {Walburn, F. J. and Schneck, D. J.},
  title   = {A constitutive equation for whole human blood},
  journal = {Biorheology},
  volume  = {13},
  number  = {3},
  pages   = {201--210},
  year    = {1976},
  doi     = {10.3233/bir-1976-13307}
}

@article{Gallagher2019,
  author  = {Gallagher, M. T. and Wain, R. A. J. and Dari, S. and Whitty, J. P. and Smith, D. J.},
  title   = {Non-identifiability of parameters for a class of shear-thinning rheological models, with implications for haematological fluid dynamics},
  journal = {J. Biomech.},
  volume  = {85},
  pages   = {230--238},
  year    = {2019},
  doi     = {10.1016/j.jbiomech.2019.01.036}
}

@article{dipietro2020improved,
  author  = {Di Pietro, D. A. and Droniou, J. and Harnist, A.},
  title   = {Improved error estimates for {H}ybrid {H}igh-{O}rder discretizations of {L}eray--{L}ions problems},
  journal = {Calcolo},
  volume  = {58},
  number  = {2},
  pages   = {Paper No. 19},
  year    = {2021},
  doi     = {10.1007/s10092-021-00410-z}
}

@article{Beldzinski2022,
  author  = {Be{\l}dzi{\'n}ski, M. and Galewski, M. and Kossowski, I.},
  title   = {Dependence on parameters for nonlinear equations---abstract principles and applications},
  journal = {Math. Methods Appl. Sci.},
  volume  = {45},
  number  = {3},
  pages   = {1668--1686},
  year    = {2022},
  doi     = {10.1002/mma.7882}
}

@article{Hirn2013,
  author  = {Hirn, A.},
  title   = {Finite element approximation of singular power-law systems},
  journal = {Math. Comp.},
  volume  = {82},
  number  = {283},
  pages   = {1247--1268},
  year    = {2013},
  doi     = {10.1090/S0025-5718-2013-02668-3}
}

@article{Lee2003BallisticSTF,
  author  = {Lee, Y. S. and Wetzel, E. D. and Wagner, N. J.},
  title   = {The ballistic impact characteristics of {K}evlar woven fabrics impregnated with a colloidal shear thickening fluid},
  journal = {J. Mater. Sci.},
  volume  = {38},
  number  = {13},
  pages   = {2825--2833},
  year    = {2003},
  doi     = {10.1023/A:1024424200221}
}

@article{Decker2007StabSTF,
  author  = {Decker, M. J. and Halbach, C. J. and Nam, C. H. and Wagner, N. J. and Wetzel, E. D.},
  title   = {Stab resistance of shear thickening fluid ({STF})-treated fabrics},
  journal = {Compos. Sci. Technol.},
  volume  = {67},
  number  = {3--4},
  pages   = {565--578},
  year    = {2007},
  doi     = {10.1016/j.compscitech.2006.08.007}
}

@article{BoscoCaladoMaia2019,
  author  = {Bosco, A. and Calado, V. and Maia, J.},
  title   = {Rheological parameters of shear-thickening fluids using an experimental design},
  journal = {Mater. Res.},
  volume  = {22},
  number  = {5},
  pages   = {Paper No. e20180631},
  year    = {2019},
  doi     = {10.1590/1980-5373-MR-2018-0631}
}

@article{Manukonda2020,
  author  = {Manukonda, B. H. and Chatterjee, V. A. and Verma, S. K. and Bhattacharjee, D. and Biswas, I. and Neogi, S.},
  title   = {Rheology based design of shear thickening fluid for soft body armor applications},
  journal = {Period. Polytech. Chem. Eng.},
  volume  = {64},
  number  = {1},
  pages   = {75--84},
  year    = {2020},
  doi     = {10.3311/PPch.13626}
}

@article{Glen1955,
  author  = {Glen, J. W.},
  title   = {The creep of polycrystalline ice},
  journal = {Proc. R. Soc. Lond. A},
  volume  = {228},
  number  = {1175},
  pages   = {519--538},
  year    = {1955},
  doi     = {10.1098/rspa.1955.0066}
}

@book{CuffeyPaterson2010,
  author    = {Cuffey, K. M. and Paterson, W. S. B.},
  title     = {The Physics of Glaciers},
  edition   = {Fourth},
  publisher = {Butterworth-Heinemann},
  address   = {Burlington, MA},
  year      = {2010}
}

@article{IsaacStadlerGhattas2015,
  author  = {Isaac, T. and Stadler, G. and Ghattas, O.},
  title   = {Solution of nonlinear {S}tokes equations discretized by high-order finite elements on nonconforming and anisotropic meshes, with application to ice sheet dynamics},
  journal = {SIAM J. Sci. Comput.},
  volume  = {37},
  number  = {6},
  pages   = {B804--B833},
  year    = {2015},
  doi     = {10.1137/140974407}
}

@article{Gagliardini2008,
  author  = {Gagliardini, O. and Zwinger, T.},
  title   = {The {ISMIP-HOM} benchmark experiments performed using the finite-element code {E}lmer},
  journal = {Cryosphere},
  volume  = {2},
  number  = {1},
  pages   = {67--76},
  year    = {2008},
  doi     = {10.5194/tc-2-67-2008}
}

@article{Pattyn2008,
  author  = {Pattyn, F. and Perichon, L. and Aschwanden, A. and Breuer, B. and de Smedt, B. and Gagliardini, O. and Gudmundsson, G. H. and Hindmarsh, R. C. A. and Hubbard, A. and Johnson, J. V. and Kleiner, T. and Konovalov, Y. and Martin, C. and Payne, A. J. and Pollard, D. and Price, S. and R{\"u}ckamp, M. and Saito, F. and Sou{\v{c}}ek, O. and Sugiyama, S. and Zwinger, T.},
  title   = {Benchmark experiments for higher-order and full-{S}tokes ice sheet models ({ISMIP-HOM})},
  journal = {Cryosphere},
  volume  = {2},
  number  = {2},
  pages   = {95--108},
  year    = {2008},
  doi     = {10.5194/tc-2-95-2008}
}

@article{Ahlkrona2021,
  author  = {Ahlkrona, J. and Elfverson, D.},
  title   = {A cut finite element method for non-{N}ewtonian free surface flows in {2D}---application to glacier modelling},
  journal = {J. Comput. Phys. X},
  volume  = {11},
  pages   = {Paper No. 100090},
  year    = {2021},
  doi     = {10.1016/j.jcpx.2021.100090}
}

@article{Millstein2022IceViscosity,
  author  = {Millstein, J. D. and Minchew, B. M. and Pegler, S. S.},
  title   = {Ice viscosity is more sensitive to stress than commonly assumed},
  journal = {Commun. Earth Environ.},
  volume  = {3},
  number  = {1},
  pages   = {Paper No. 57},
  year    = {2022},
  doi     = {10.1038/s43247-022-00385-x}
}

@article{Fan2025FlowLawsIce,
  author  = {Fan, S. and Wang, T. and Prior, D. J. and Breithaupt, T. and Hager, T. F. and Wallis, D.},
  title   = {Flow laws for ice constrained by 70 years of laboratory experiments},
  journal = {Nat. Geosci.},
  volume  = {18},
  number  = {4},
  pages   = {296--304},
  year    = {2025},
  doi     = {10.1038/s41561-025-01661-z}
}

@article{Cuzzone2025PolythermalISSM,
  author  = {Cuzzone, J. K. and Barth, A. and Barker, K. and Morlighem, M.},
  title   = {Ice sheet model simulations reveal that polythermal ice conditions existed across the northeastern {USA} during the {L}ast {G}lacial {M}aximum},
  journal = {Cryosphere},
  volume  = {19},
  number  = {4},
  pages   = {1559--1575},
  year    = {2025},
  doi     = {10.5194/tc-19-1559-2025}
}

@article{Leger2025GreenlandPISM,
  author  = {Leger, T. P. M. and Ely, J. C. and Clark, C. D. and Bradley, S. L. and Archer, R. E. and Zhu, J.},
  title   = {The {Greenland-Ice-Sheet} evolution over the last 24\,000 years: insights from model simulations evaluated against ice-extent markers},
  journal = {Cryosphere},
  volume  = {19},
  number  = {11},
  pages   = {5719--5761},
  year    = {2025},
  doi     = {10.5194/tc-19-5719-2025}
}

@article{alves2013study,
  author  = {Alves, M. M. and Rocha, C. and Gon{\c{c}}alves, M. P.},
  title   = {Study of the rheological behaviour of human blood using a controlled stress rheometer},
  journal = {Clin. Hemorheol. Microcirc.},
  volume  = {53},
  number  = {4},
  pages   = {369--386},
  year    = {2013},
  doi     = {10.3233/CH-121645}
}

@article{SchmidtHumbertSlawig2024,
  author  = {Schmidt, N. and Humbert, A. and Slawig, T.},
  title   = {Assessing the benefits of approximately exact step sizes for {P}icard and {N}ewton solver in simulating ice flow ({FEniCS}-full-{S}tokes v.1.3.2)},
  journal = {Geosci. Model Dev.},
  volume  = {17},
  number  = {12},
  pages   = {4943--4959},
  year    = {2024},
  doi     = {10.5194/gmd-17-4943-2024}
}

@article{HelanowAhlkrona2025,
  author  = {Helanow, C. and Ahlkrona, J.},
  title   = {Theoretical results on a block preconditioner used in ice-sheet modeling: eigenvalue bounds for singular power-law fluids},
  journal = {Comput. Geosci.},
  volume  = {29},
  number  = {1},
  pages   = {Paper No. 9},
  year    = {2025},
  doi     = {10.1007/s10596-025-10339-9}
}

@article{ErvinHowellStanculescu2008,
  author  = {Ervin, V. J. and Howell, J. S. and Stanculescu, I.},
  title   = {A dual-mixed approximation method for a three-field model of a nonlinear generalized {S}tokes problem},
  journal = {Comput. Methods Appl. Mech. Engrg.},
  volume  = {197},
  number  = {33--40},
  pages   = {2886--2900},
  year    = {2008},
  doi     = {10.1016/j.cma.2008.01.022}
}

@article{Howell2009DualMixed,
  author  = {Howell, J. S.},
  title   = {Dual-mixed finite element approximation of {S}tokes and nonlinear {S}tokes problems using trace-free velocity gradients},
  journal = {J. Comput. Appl. Math.},
  volume  = {231},
  number  = {2},
  pages   = {780--792},
  year    = {2009},
  doi     = {10.1016/j.cam.2009.05.002}
}

@article{FarhloulZine2017DualMixed,
  author  = {Farhloul, M. and Zine, A.},
  title   = {A dual-mixed finite element method for quasi-{N}ewtonian flows whose viscosity obeys a power law or the {C}arreau law},
  journal = {Math. Comput. Simulation},
  volume  = {141},
  pages   = {83--95},
  year    = {2017},
  doi     = {10.1016/j.matcom.2016.09.015}
}

@book{EkelandTemam1999,
  author    = {Ekeland, I. and T{\'e}mam, R.},
  title     = {Convex Analysis and Variational Problems},
  series    = {Classics in Applied Mathematics},
  volume    = {28},
  publisher = {SIAM},
  address   = {Philadelphia, PA},
  year      = {1999},
  doi       = {10.1137/1.9781611971088}
}

@article{Rockafellar1968,
  author  = {Rockafellar, R. T.},
  title   = {Integrals which are convex functionals},
  journal = {Pacific J. Math.},
  volume  = {24},
  number  = {3},
  pages   = {525--539},
  year    = {1968},
  doi     = {10.2140/pjm.1968.24.525}
}

@book{BoydVandenberghe2004,
  author    = {Boyd, S. and Vandenberghe, L.},
  title     = {Convex Optimization},
  publisher = {Cambridge University Press},
  address   = {Cambridge},
  year      = {2004},
  doi       = {10.1017/CBO9780511804441}
}

@article{KouhiaStenberg1995,
  author  = {Kouhia, R. and Stenberg, R.},
  title   = {A linear nonconforming finite element method for nearly incompressible elasticity and {S}tokes flow},
  journal = {Comput. Methods Appl. Mech. Engrg.},
  volume  = {124},
  number  = {3},
  pages   = {195--212},
  year    = {1995},
  doi     = {10.1016/0045-7825(95)00829-P}
}

@article{HuSchedensack2019,
  author  = {Hu, J. and Schedensack, M.},
  title   = {Two low-order nonconforming finite element methods for the {S}tokes flow in three dimensions},
  journal = {IMA J. Numer. Anal.},
  volume  = {39},
  number  = {3},
  pages   = {1447--1470},
  year    = {2019},
  doi     = {10.1093/imanum/dry021}
}

@article{John2017divergence,
  author  = {John, V. and Linke, A. and Merdon, C. and Neilan, M. and Rebholz, L. G.},
  title   = {On the divergence constraint in mixed finite element methods for incompressible flows},
  journal = {SIAM Rev.},
  volume  = {59},
  number  = {3},
  pages   = {492--544},
  year    = {2017},
  doi     = {10.1137/15M1047696}
}

@article{Linke2014role,
  author  = {Linke, A.},
  title   = {On the role of the {H}elmholtz decomposition in mixed methods for incompressible flows and a new variational crime},
  journal = {Comput. Methods Appl. Mech. Engrg.},
  volume  = {268},
  pages   = {782--800},
  year    = {2014},
  doi     = {10.1016/j.cma.2013.10.011}
}

@article{lederer2017divergence,
  author  = {Lederer, P. L. and Linke, A. and Merdon, C. and Sch{\"o}berl, J.},
  title   = {Divergence-free reconstruction operators for pressure-robust {S}tokes discretizations with continuous pressure finite elements},
  journal = {SIAM J. Numer. Anal.},
  volume  = {55},
  number  = {3},
  pages   = {1291--1314},
  year    = {2017},
  doi     = {10.1137/16M1089964}
}

@article{ArnoldBrezziFortin1984,
  author  = {Arnold, D. N. and Brezzi, F. and Fortin, M.},
  title   = {A stable finite element for the {Stokes} equations},
  journal = {Calcolo},
  volume  = {21},
  number  = {4},
  pages   = {337--344},
  year    = {1984},
  doi     = {10.1007/BF02576171}
}

@article{CarstensenSchedensack2015,
  author  = {Carstensen, C. and Schedensack, M.},
  title   = {Medius analysis and comparison results for first-order finite element methods in linear elasticity},
  journal = {IMA J. Numer. Anal.},
  volume  = {35},
  number  = {4},
  pages   = {1591--1621},
  year    = {2015},
  doi     = {10.1093/imanum/dru048}
}

@article{HeidSuli2022,
  author  = {Heid, P. and S\"uli, E.},
  title   = {On the convergence rate of the {Ka\v{c}anov} scheme for shear-thinning fluids},
  journal = {Calcolo},
  volume  = {59},
  number  = {1},
  pages   = {Paper No. 4},
  year    = {2022},
  doi     = {10.1007/s10092-021-00444-3}
}

@article{LoggWells2010,
  author  = {Logg, A. and Wells, G. N.},
  title   = {{DOLFIN}: Automated finite element computing},
  journal = {ACM Trans. Math. Softw.},
  volume  = {37},
  number  = {2},
  pages   = {20:1--20:28},
  year    = {2010},
  doi     = {10.1145/1731022.1731030}
}

\appendix

\section{Derivation of the energy functionals and duality gap}\label{app:duality-gap}

Throughout this appendix, we derive the primal and dual energy functionals and establish their duality gap at the continuous level.

\subsection{Derivation of the energy functionals}

We now prove the primal--dual equality stated in Theorem~\ref{thm:primal-dual-equality}. The proof is divided into three parts:
\begin{enumerate}
    \item \textbf{Part 1.} We first rewrite the original Lagrangian $\mathcal{L}$ as $\sup_{\bsigma\in Y^*}\mathcal{M}$ by introducing the dual stress variable $\bsigma$. This representation allows us to eliminate the velocity explicitly at a later stage.

    \item \textbf{Part 2.} We perform the dualization at fixed pressure $p$. The conjugate of the translated bulk energy accounts for the Dirichlet lifting, whereas the conjugate of the linear velocity term enforces the equilibrium condition $(\bsigma,p)\in\mathcal{E}$. Combining these ingredients via Lemma~\ref{lem:conj_translated-functional} yields $\mathcal{J}_D$.

    \item \textbf{Part 3.} Finally, we allow the pressure to vary and show that the infimum over velocity can be exchanged with the joint supremum over pressure and stress, attaining the same saddle point and thus concluding the proof.
\end{enumerate}

\paragraph{Notation.}

We use the notation of Section~\ref{sec:continuouspb}.
Set \(Y=L^r(\Omega;\Sym)\), so that \(Y^*=\Sigma\), and define the symmetric-gradient operator \(K:V_0\to Y\) by \(K\bw=\eps(\bw)\). Define also the bulk functional \(I:Y\to\R\) by
\(
I(\bz):=\int_\Omega \frac{\mu}{r}|\bz|^r\,\dd x.
\)
Define \(B:V_0\to Q^*\) by
\(
\dual{B\bw}{q}_{Q^*,Q}:=\int_\Omega q\,\operatorname{div}\bw\,\dd x.
\)
We write \(K^*:\Sigma\to V_0^*\) and \(B^*:Q\to V_0^*\) for the adjoints. Recall also the pointwise identity \(\mathbf I:\varepsilon(\boldsymbol w)=\operatorname{tr}(\varepsilon(\boldsymbol w))=\operatorname{div}\boldsymbol w\). Using it yields \(p\operatorname{div}\bw=p\mathbf I:\eps(\bw)\), so
\(
\dual{K^*\btau-B^*p}{\bw}_{V_0^*,V_0}=\int_\Omega (\btau-p\mathbf I):\eps(\bw)\,\dd x.
\)
For a proper functional \(F:X\to\R\cup\{+\infty\}\) on a Banach space \(X\), its Fenchel conjugate is
\(
F^*(x^*)=\sup_{x\in X}\left\{\dual{x^*}{x}_{X^*,X}-F(x)\right\}.
\)
See \cite[Chap.~III, Thm.~4.1 and Rem.~4.2]{EkelandTemam1999} for the Fenchel--Rockafellar duality and \cite[Chap.~I, Prop.~4.1]{EkelandTemam1999} for the Fenchel--Moreau representation.

\paragraph{Part 1: from $\mathcal{L}$ to $\mathcal{M}$.}
First, for $1<r<\infty$ and $\mu$ satisfying $0<\mu_-\leq\mu(x)\leq\mu_+<\infty$, we can exchange conjugation and integration \cite[Thm.~2]{Rockafellar1968}; hence, the conjugate of the bulk functional is given, for all $\bsigma\in Y^*$, by
\[
I^*(\bsigma)
=\int_\Omega \varphi^*(x,\bsigma(x))\,\dd x.
\]
By the Fenchel--Moreau theorem, we obtain the dual representation
\[
I(\eps(\bv))
=
\sup_{\bsigma\in Y^*}
\left\{
\int_\Omega
\bsigma:\eps(\bv)\,\dd x
-
I^*(\bsigma)
\right\}.
\]

\begin{proposition}[Primal--dual representation]\label{prop:primal-dual-rep}
For every $(\bv,q)\in V_D\times Q$,
\[
\mathcal{L}(\bv,q)
=
\sup_{\bsigma\in Y^*}
\mathcal{M}(\bv,q,\bsigma).
\]
\end{proposition}

\begin{proof}
Substituting $I$ with its Fenchel--Moreau representation inside $\mathcal{L}$ and using the pointwise identity $\mathbf{I}:\eps(\bv)=\operatorname{div}\bv$ directly yields the claim.
\end{proof}

\noindent
Decomposing $\bu=\bw_D+\bw$ with $\bw\in V_0$, we have $I(\eps(\bu)) = I(K\bw+K\bw_D)$. We therefore define the translated bulk functional $\widetilde{I}$ by $\widetilde{I}(\bz) := I(\bz+K\bw_D)$.

\paragraph{Part 2: from $\mathcal{M}$ to $\mathcal{J}_D$ at fixed pressure.}

For a fixed $p\in Q$, we define
\[
G_p(\bw)
:=
-\ell(\bw)
-
\dual{B\bw}{p}_{Q^*,Q}.
\]

\begin{lemma}[Conjugates]\label{lem:conj_translated-functional}
For every $\bsigma\in Y^*$ and $p\in Q$, 
\[
\widetilde{I}^*(\bsigma)
=
I^*(\bsigma)
-
\int_\Omega \bsigma:\eps(\bw_D)\,\dd x,
\qquad
G_p^*(-K^*\bsigma)<+\infty
\iff
\underbrace{K^*\bsigma-B^*p=\ell}
_{
% \forall\,\bw\in V_0~:~
\int_\Omega
(\bsigma-p\mathbf{I}):\eps(\bw)\,\dd x
=
\ell(\bw)
% \quad \forall\,\bw\in V_0.
}
\]
\end{lemma}

\begin{proof}
First, we prove $\tilde I^*$.
By definition,
\[
\widetilde{I}^*(\bsigma)
=
\sup_{\bz\in Y}
\left\{
\dual{\bsigma}{\bz}
-
I(\bz+K\bw_D)
\right\}.
\]
Setting $\by:=\bz+K\bw_D$, so that $\bz=\by-K\bw_D$, we obtain
\[
\widetilde{I}^*(\bsigma)
=
\sup_{\by\in Y}
\left\{
\dual{\bsigma}{\by}
-
I(\by)
\right\}
-
\dual{\bsigma}{K\bw_D}
=
I^*(\bsigma)
-
\dual{\bsigma}{K\bw_D},
\]
which coincides with the stated expression.

Now we prove $G^*_p$. By the definitions of the Fenchel conjugate and of $B^*$, for any $\xi\in V_0^*$ we have
\[
G_p^*(\xi)
=\sup_{\bw\in V_0}\dual{\xi+\ell+B^*p}{\bw}.
\]
This supremum equals $0$ if $\xi+\ell+B^*p=0$ and $+\infty$ otherwise (cf.~\cite[\S~5.2.1, p.~225]{BoydVandenberghe2004}). Setting $\xi=-K^*\bsigma$ and using the identity $p\,\operatorname{div}\bw = p\mathbf{I}:\eps(\bw)$ yields the stated equilibrium condition.
\end{proof}

\noindent
For $(\bsigma,p)\in\mathcal{E}$, define $\mathcal{J}_D(\bsigma,p):=-\mathcal{J}^*(\bsigma,p).$

\begin{proposition}[Dual problem at fixed pressure]\label{prop:fixed-pressure}
For every fixed $p\in Q$,
\[
  \inf_{\bv\in V_D}\mathcal{L}(\bv,p)
  =
  \sup_{\bsigma\in Y^*}\inf_{\bv\in V_D}
       \mathcal{M}(\bv,p,\bsigma)
  =
  \sup_{\substack{\bsigma\in Y^*\\(\bsigma,p)\in\mathcal{E}}}
       \mathcal{J}_D(\bsigma,p).
\]
\end{proposition}

\begin{proof}
Fix $p\in Q$ and write $\bv=\bw_D+\bw$ with $\bw\in V_0$. Collect the load and pressure contributions depending solely on the lifting in
\[
  C_p:=-\ell(\bw_D)
       -\int_\Omega p\,\operatorname{div}\bw_D\,\dd x.
\]
By the definitions of $\widetilde{I}$, $K$, and $G_p$, we have
\[
  \mathcal{L}(\bw_D+\bw,p)
  =\widetilde{I}(K\bw)+G_p(\bw)+C_p.
\]
Recall from Proposition~\ref{prop:primal-dual-rep} that $\mathcal{L}(\bv,p)=\sup_{\bsigma\in Y^*}\mathcal{M}(\bv,p,\bsigma)$. Since $G_p(0)=0$ and $\widetilde{I}$ is finite and continuous at $K0=0$, Fenchel--Rockafellar duality \cite[Chap.~III, Thm.~4.1 and Rem.~4.2]{EkelandTemam1999} allows us to exchange the velocity infimum and the stress supremum, yielding the first equality.

To prove the second equality, fix $\bsigma\in Y^*$ and group terms of $\mathcal{M}$ according to their dependence on $\bw$:
\[
  \mathcal{M}(\bw_D+\bw,p,\bsigma)
  = C_p
  +\underbrace{\int_\Omega \bsigma:\eps(\bw_D)\,\dd x
                -I^*(\bsigma)}_{\text{translated bulk term}}
  +\underbrace{\dual{\bsigma}{K\bw}+G_p(\bw)}
               _{\text{homogeneous velocity term}}.
\]
By Lemma~\ref{lem:conj_translated-functional}, the first bracket equals $-\widetilde{I}^*(\bsigma)$, while only the second bracket enters the infimum over $\bw$. By the definition of the Fenchel conjugate and Lemma~\ref{lem:conj_translated-functional},
\[
  \inf_{\bw\in V_0}\bigl\{\dual{\bsigma}{K\bw}+G_p(\bw)\bigr\}
  =-G_p^*(-K^*\bsigma)
  % =
  % \begin{cases}
  %   0, &(\bsigma,p)\in\mathcal{E},\\
  %   -\infty, &(\bsigma,p)\notin\mathcal{E}.
  % \end{cases}
\]
Combining the translated bulk term with $C_p$ and using $\mathbf{I}:\eps(\bw_D)=\operatorname{div}\bw_D$ gives, for $(\bsigma,p)\in\mathcal{E}$,
\[
  \inf_{\bv\in V_D}\mathcal{M}(\bv,p,\bsigma)
  =
  \int_\Omega(\bsigma-p\mathbf{I}):\eps(\bw_D)\,\dd x
  -I^*(\bsigma)-\ell(\bw_D)
  =\mathcal{J}_D(\bsigma,p).
\]
Taking the supremum over $\bsigma\in Y^*$ establishes the second equality.
\end{proof}

\paragraph{Part 3: exchanging the optimization order with varying pressure.}

We now pass from fixed pressure to the full coupled problem. Since Proposition~\ref{prop:fixed-pressure} holds for every $p\in Q$, taking the supremum over $p$ yields
\begin{equation}\label{eq:full-dual-value}
  \sup_{(p,\bsigma)\in Q\times Y^*}
  \inf_{\bv\in V_D}\mathcal{M}(\bv,p,\bsigma)
  =
  \sup_{(\bsigma,p)\in\mathcal{E}}\mathcal{J}_D(\bsigma,p).
\end{equation}
Thus, Proposition~\ref{prop:fixed-pressure} already identifies the dual value when the velocity is minimized first.

In the opposite order, the supremum over $\bsigma$ recovers $\mathcal{L}(\bv,p)$ by Proposition~\ref{prop:primal-dual-rep}, and the supremum over $p$ enforces weak incompressibility:
\begin{equation}\label{eq:full-primal-value}
  \inf_{\bv\in V_D}
  \sup_{(p,\bsigma)\in Q\times Y^*}\mathcal{M}(\bv,p,\bsigma)
  =
  \inf_{\bv\in\mathcal{A}}\mathcal{J}(\bv).
\end{equation}
It remains to show that \eqref{eq:full-primal-value} and \eqref{eq:full-dual-value} share the same saddle point, so that
\begin{equation}\label{eq:minimax-exchange}
  \sup_{(p,\bsigma)\in Q\times Y^*}
  \inf_{\bv\in V_D}\mathcal{M}(\bv,p,\bsigma)
  =
  \inf_{\bv\in V_D}
  \sup_{(p,\bsigma)\in Q\times Y^*}\mathcal{M}(\bv,p,\bsigma).
\end{equation}

\begin{proof}[Proof of~\eqref{eq:minimax-exchange} and Theorem~\ref{thm:primal-dual-equality}]
Let $(\bar{\bu},\bar{p})$ be the continuous weak solution, and set $\bar{\bsigma}:=\mu|\eps(\bar{\bu})|^{r-2}\eps(\bar{\bu})$. For every $q\in Q$ and $\btau\in Y^*$, the incompressibility of $\bar{\bu}$ and the Fenchel--Young inequality yield
\[
 \mathcal{M}(\bar{\bu},q,\btau)
 \leq \mathcal{J}(\bar{\bu})
 = \mathcal{M}(\bar{\bu},\bar{p},\bar{\bsigma}).
\]
Moreover, for every $\bv\in V_D$, testing~\eqref{eq:dual-three-field-continuous} with $\bv-\bar{\bu}\in V_0$ gives
\[
 \mathcal{M}(\bv,\bar{p},\bar{\bsigma})
 -\mathcal{M}(\bar{\bu},\bar{p},\bar{\bsigma})
 =
 \int_\Omega(\bar{\bsigma}-\bar{p}\mathbf{I}):
       \eps(\bv-\bar{\bu})\,\dd x-\ell(\bv-\bar{\bu})=0.
\]
Combining these two relations shows that, for every $\bv\in V_D$ and $(q,\btau)\in Q\times Y^*$,
\[
  \mathcal{M}(\bar{\bu},q,\btau)
  \leq \mathcal{M}(\bar{\bu},\bar{p},\bar{\bsigma})
  = \mathcal{M}(\bv,\bar{p},\bar{\bsigma}).
\]
The first inequality implies that $(\bar{p},\bar{\bsigma})$ maximizes $\mathcal{M}$ at fixed velocity $\bar{\bu}$, whereas the equality implies that $\bar{\bu}$ minimizes $\mathcal{M}$ at fixed $(\bar{p},\bar{\bsigma})$. Hence $(\bar{\bu},\bar{p},\bar{\bsigma})$ is a saddle point of $\mathcal{M}$, which establishes~\eqref{eq:minimax-exchange}. Combining~\eqref{eq:full-dual-value}, \eqref{eq:full-primal-value}, and~\eqref{eq:minimax-exchange} with $\mathcal{J}^*=-\mathcal{J}_D$ proves Theorem~\ref{thm:primal-dual-equality}.
\end{proof}

\subsection{Duality gap}
We empirically further confirm the correctness of our Fenchel-dual functionals by showing that the primal-dual gap approaches zero as the dual-weighted Ka\v{c}anov scheme converges.

\begin{definition}[Primal--dual gap]
Let $\mathbf u\in\mathcal A$ be primal-admissible and let
$(\bsigma,p)\in\mathcal E$ be dual-admissible.
We define the primal--dual gap by
\[
\mathcal G_{\mathrm{PD}}(\mathbf u,\bsigma,p)
:=
\mathcal J(\mathbf u)
+
\mathcal J^*(\bsigma,p).
\]
\end{definition}
\noindent
Using the primal and dual admissibility conditions, this can be written as
\[
\mathcal G_{\mathrm{PD}}(\mathbf u,\bsigma,p)
=
\int_\Omega
\left[
\frac{\mu}{r}
\|\boldsymbol{\varepsilon}(\mathbf u)\|^r
+
\frac{1}{r'}
\mu^{-\frac{1}{r-1}}
\|\bsigma\|^{r'}
-
\bsigma:
\boldsymbol{\varepsilon}(\mathbf u)
\right]\,dx.
\]
By the Fenchel--Young inequality,
\[
\mathcal G_{\mathrm{PD}}(\mathbf u,\bsigma,p)\ge0,
\]
with equality at the primal--dual solution.

\begin{figure}[H]
    \centering
    \includegraphics[width=0.5\linewidth]{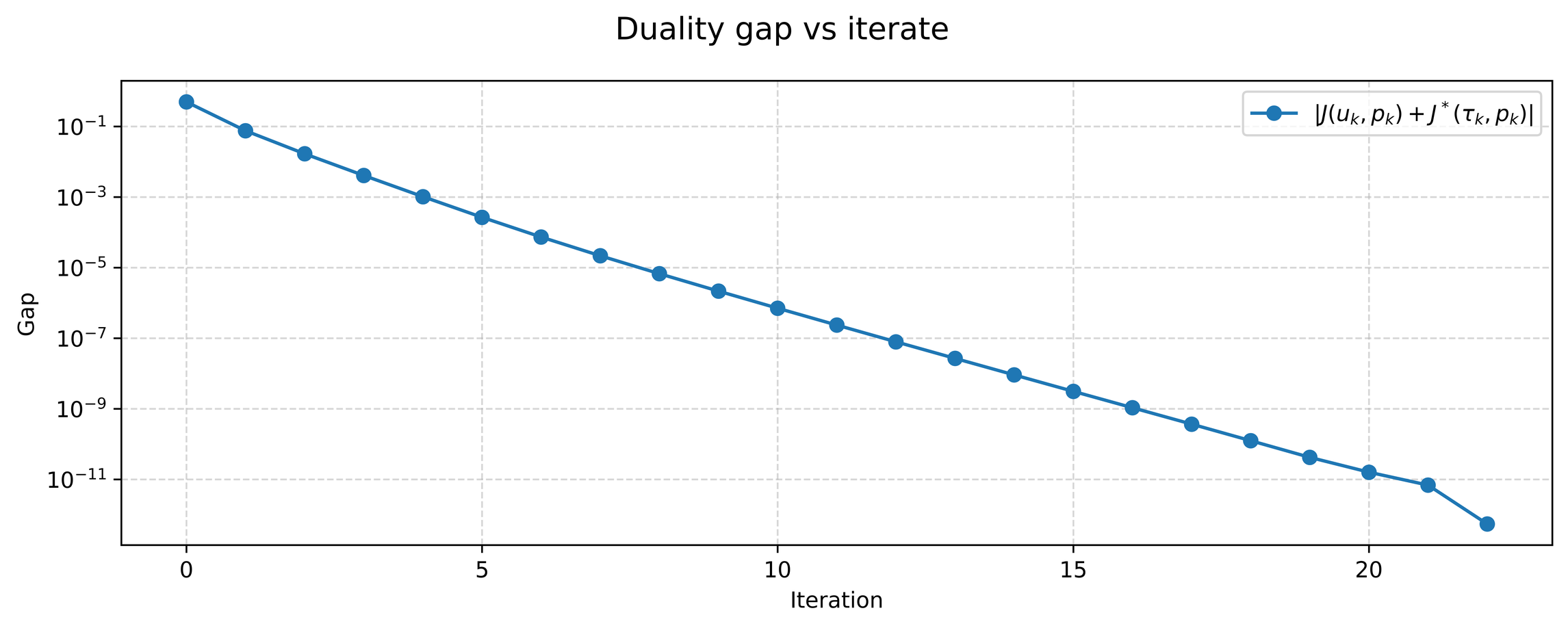}
    \caption{Duality Gap $\mathcal G_{\mathrm{PD}}(\mathbf u_k,\bsigma_k,p_k)$ against iterations for dual-weighted Ka\v{c}anov. KS, $n_h=128$, $r=3.5$, trigonometric benchmark.}
\end{figure}

\FloatBarrier
\section{Equivalence of the two dual Ka\v{c}anov iterations on KS}
\label{app:proof-equivalence-dual-schemes}

We prove Proposition~\ref{prop:dual-kacanov-equivalence}.

\begin{proof}
% As every $\bv_h\in\boldsymbol{V}_h^{\mathrm{KS}}$
% is piecewise affine, on each element, $\beps{\bu_h}$ and $\beps$. 
Since every $\bv_h\in\boldsymbol{V}_h^{\mathrm{KS}}$
is piecewise affine, $\beps(\bv_h)$ and $\nabla\cdot\bv_h=\tr\beps(\bv_h)$ are piecewise constant.
For $\btau\in\mathbb{R}^{2\times2}_{\mathrm{sym}}$ we use $\btau=\dev(\btau)+\tfrac12\tr(\btau)\mathbf I$. We denote with $\chi_K$ the indicator function on the element $K$. Since $\bw_{D,h}$ is chosen discretely divergence-free, i.e.\ $b(\bw_{D,h},q_h)=0$ for all $q_h\in Q_h^{\mathrm{P0}}$, both iterations depend on $\bu_h^{k+1}$ only through the full velocity $\bu_h^{k+1}+\bw_{D,h}$, which is again piecewise affine: in \eqref{eq:dual-kacanov-primal} the incompressibility term can be written as $-\int_\Omega q_h\,\nabla\cdot(\bu_h^{k+1}+\bw_{D,h})\dd x$, and in \eqref{eq:dual-two-field-kacanov} the lifting term on the right-hand side combines with $-\int_\Omega\beps(\bu_h^{k+1}):\mathbf{T}_h\dd x$ into $-\int_\Omega\beps(\bu_h^{k+1}+\bw_{D,h}):\mathbf{T}_h\dd x$. To lighten notation, in the rest of the proof we write $\bu_h^{k+1}$ for $\bu_h^{k+1}+\bw_{D,h}$. We proceed by induction. Remember our chosen finite elements guarantee uniqueness of the solution.

\emph{Base case.} The base case $\dev(\boldsymbol{S}^0_h)=\bsigma_0^h$ is true by hypothesis. 

\emph{Induction.} Assume $\bsigma_h^k=\dev(\mathbf{S}_h^k)$.
By definition of $(\varphi^*)'(s)=\mu^{-\frac{1}{r-1}}s^{r'-1}$, the relationship between the weight in
\eqref{eq:dual-kacanov-primal} and the coefficient $c_k$ in the first line of \eqref{eq:dual-two-field-kacanov} is, since they are both piecewise constant and well defined,
\[
  w_k
  :=
  \frac{\|\bsigma_h^k\|}{(\varphi^*)'(\|\bsigma_h^k\|)}
  =
  \mu^{\frac{1}{r-1}}\|\bsigma_h^k\|^{2-r'}
  =
  \dfrac {1}{
  \mu^{-\frac{1}{r-1}}\|\dev(\mathbf{S}_h^k)\|^{r'-2}
  }
  =
  \dfrac {1}{c_k}
\]
We now show each solution of \eqref{eq:dual-two-field-kacanov} is a solution of \eqref{eq:dual-kacanov-primal} and vice versa.

\emph{\eqref{eq:dual-two-field-kacanov} $\implies$ \eqref{eq:dual-kacanov-primal}.} Let now $(\mathbf{S}_h^{k+1},\bu_h^{k+1})$ solve
\eqref{eq:dual-two-field-kacanov}. Testing its first line with $\mathbf{T}_h=\chi_K\mathbf I$, $K\in\mathcal{T}_h$, using $\dev(\mathbf{S}_h^{k+1}):\mathbf I=0$ and $\boldsymbol\varepsilon(\mathbf u_h^{k+1}):\mathbf I = \nabla\cdot\bu_h^{k+1}$ gives
\[
 \int_\Omega c_k\,\operatorname{dev}(\mathbf S_h^{k+1}):\mathbf T_h\,\dd x 
  = 
  \int_\Omega \boldsymbol\varepsilon(\mathbf u_h^{k+1}):\mathbf T_h\,\dd x 
  \implies
\int_K \nabla\cdot\bu_h^{k+1}\,\dd x=0.
\]
Combining this with the fact that $\nabla\cdot\bu_h^{k+1}$ is constant on $K$, we conclude that
\begin{equation}\label{eq:equiv-divfree}
  \nabla\cdot\bu_h^{k+1}=0\quad\text{on every }K\in\mathcal{T}_h,
\end{equation}
so $\beps(\bu_h^{k+1})$ is trace free. We now let $\bv_h=0$ in \eqref{eq:dual-two-field-kacanov} and test with trace-free $\mathbf{T}_h=\chi_K\mathbf{T}$,
$\mathbf{T}\in\mathbb{R}^{2\times2}_{\mathrm{sym}}$, $\tr\mathbf{T}=0$. Using that both
$c_k\dev(\mathbf{S}_h^{k+1})$ and $\beps(\bu_h^{k+1})$ are trace free on $K$, we may choose in particular $\mathbf{T} =\left. \bigl(c_k\dev(\mathbf{S}_h^{k+1})-\beps(\bu_h^{k+1})\bigr) \right|_K$

% \begin{equation}
\begin{align}
    \int_\Omega 
\mu^{-\frac{1}{r-1}} 
\bigl\|\operatorname{dev}(\mathbf S_h^k)\bigr\|^{r'-2} 
\operatorname{dev}(\mathbf S_h^{k+1}):\mathbf T_h\,
\dd x 
&= 
\int_\Omega 
\boldsymbol\varepsilon(\mathbf u_h^{k+1}):\mathbf T_h\,\dd x 
\\
\implies
  \int_K 
    \bigl\|c_k\dev(\mathbf{S}_h^{k+1})-\beps(\bu_h^{k+1})\bigr\|^2 
    \,\dd x
    &=
    0.
\end{align}
% \end{equation}
Using that both $c_k\dev(\mathbf{S}_h^{k+1})$, $\beps(\bu_h^{k+1})$ are constant on $K$, and multiplying both sides by $w_k$ we can conclude
\begin{equation}\label{eq:equiv-constitutive}
  \dev(\mathbf{S}_h^{k+1})=w_k\,\beps(\bu_h^{k+1})\quad\text{on every }K\in\mathcal{T}_h .
\end{equation}
Fix any $\bv_h\in\boldsymbol{V}_h^{\mathrm{KS}}$ and take 
$\mathbf{T}_h=\mathbf 0$ in \eqref{eq:dual-two-field-kacanov}. 
With this choice, the two integrals involving $\mathbf T_h$ vanish: 
\[ 
\begin{aligned} 
  \underbrace{\int_\Omega 
    c_k\,\dev(\mathbf{S}_h^{k+1}):\mathbf 0\dd x}_{=0} 
  - 
  \underbrace{\int_\Omega 
    \beps(\bu_h^{k+1}):\mathbf 0\dd x}_{=0} 
  +\int_\Omega \mathbf{S}_h^{k+1}:\beps(\bv_h)\dd x 
  =\ell(\bv_h). 
\end{aligned}
  \] 
  % The remaining equation is 
  % \[ 
  %   \int_\Omega \mathbf{S}_h^{k+1}:\beps(\bv_h)\dd x 
  %   =\ell(\bv_h). 
  % \] 
  Substituting $p:=-\tfrac12\tr(\mathbf{S}_h^{k+1})\in Q_h^{\mathrm{P0}}$ and \ref{eq:equiv-constitutive} into the preceding equation, and using   \(
    \mathbf I:\beps(\bv_h) 
    =\tr(\beps(\bv_h)) 
    =\nabla\cdot\bv_h, 
  \) we obtain
  \[ 
    \begin{aligned} 
      \int_\Omega \mathbf{S}_h^{k+1}:\beps(\bv_h)\dd x 
      &= 
      \int_\Omega 
      \bigl(w_k\beps(\bu_h^{k+1})-p\,\mathbf I\bigr) 
      :\beps(\bv_h)\dd x 
      \\ 
      &= 
      \int_\Omega w_k\,\beps(\bu_h^{k+1}):\beps(\bv_h)\dd x 
      - 
      \int_\Omega p \nabla\cdot\bv_h, \dd x. 
    \end{aligned} 
  \] 
% Combining it with \eqref{eq:equiv-constitutive} we have
% $\mathbf{S}_h^{k+1}=w_k\beps(\bu_h^{k+1})-p\,\mathbf I$,
% and the second line of
% \eqref{eq:dual-two-field-kacanov} becomes
% \[
%   \int_\Omega w_k\,\beps(\bu_h^{k+1}):\beps(\bv_h)\,dx
%   -
%   \int_\Omega p\,\nabla\cdot\bv_h\,dx
%   =
%   \ell(\bv_h)
%   \qquad
%   \forall\,\bv_h\in\boldsymbol{V}_h^{\mathrm{KS}}.
% \]
Together with \eqref{eq:equiv-divfree}, this shows that $(\bu_h^{k+1},p)$ solves
\eqref{eq:dual-kacanov-primal}, thus \eqref{eq:dual-two-field-kacanov} $\implies$ \eqref{eq:dual-kacanov-primal}.

\emph{\eqref{eq:dual-kacanov-primal} $\implies$ \eqref{eq:dual-two-field-kacanov}.}
Conversely, if $(\bu,p)$ solves \eqref{eq:dual-kacanov-primal}, then testing with $\bv_h=\mathbf 0$ and $q_h\in Q_h^{\mathrm{P0}}$ gives $\nabla\cdot\bu=0$ elementwise, and $\mathbf{S}:=w_k\beps(\bu)-p\,\mathbf I$ satisfies both lines of \eqref{eq:dual-two-field-kacanov}: the first pointwise, because $c_k\dev(\mathbf{S})=c_k w_k\beps(\bu)=\beps(\bu)$, and the second by \eqref{eq:dual-kacanov-primal} with $q_h=0$. Since $\mathbf{S}$ is determined by $(\dev\mathbf{S},\tr\mathbf{S})$, and the pair $(\bu,p)$ is unique, the solution of \eqref{eq:dual-two-field-kacanov} is unique as well (up to the normalisation of $\tr\mathbf{S}$ when $\Gamma_D=\partial\Omega$).
Hence $\bu_h^{k+1}$ and  $p_h^{k+1}=-\tfrac12\tr(\mathbf{S}_h^{k+1})$ agree in both schemes. Finally, by \eqref{eq:duality-u-sigma} and \eqref{eq:equiv-constitutive},
$\bsigma_h^{k+1}=w_k\beps(\bu_h^{k+1})=\dev(\mathbf{S}_h^{k+1})$.
\end{proof}
\FloatBarrier
% \input{paper_sections/AppendixImplementation}
% \FloatBarrier
\section{Additional Experiments}
\subsection{Obstacle-domain test}\label{app:obstacle-domain-test}

We finally consider a flow problem in a channel containing a rigid circular obstacle. Since no exact solution is available, this test is used to assess the robustness of the discretizations and nonlinear solvers on a more complex geometry.

The configuration also admits an (idealized) hemodynamic interpretation: the channel can be viewed as a vessel-like domain, while the obstacle represents a simplified intraluminal obstruction or device strut. It is not a complete arterial-flow simulation, as it would be beyond the scope of this work.

The body force vanishes, $\ff=(0,0)^\top$, and the boundary conditions are
\begin{equation}
\bu=
\begin{cases}
\bigl(4(1-y^2),\,0\bigr)^\top
    & \text{on } \Gamma_{\mathrm{left}},\\[0.5ex]
(0,0)^\top
    & \text{on } \Gamma_{\mathrm{rigid}},
\end{cases}
\qquad
(\bsigma\bn-p\bn)=(0,0)^\top
\quad \text{on } \Gamma_{\mathrm{right}}.
\end{equation}

The numerical results display stable, consistent convergence on this geometry as well.

\begin{figure}[htbp]
    \centering
    \begin{subfigure}{0.48\linewidth}
        \centering
        \includegraphics[width=\linewidth]{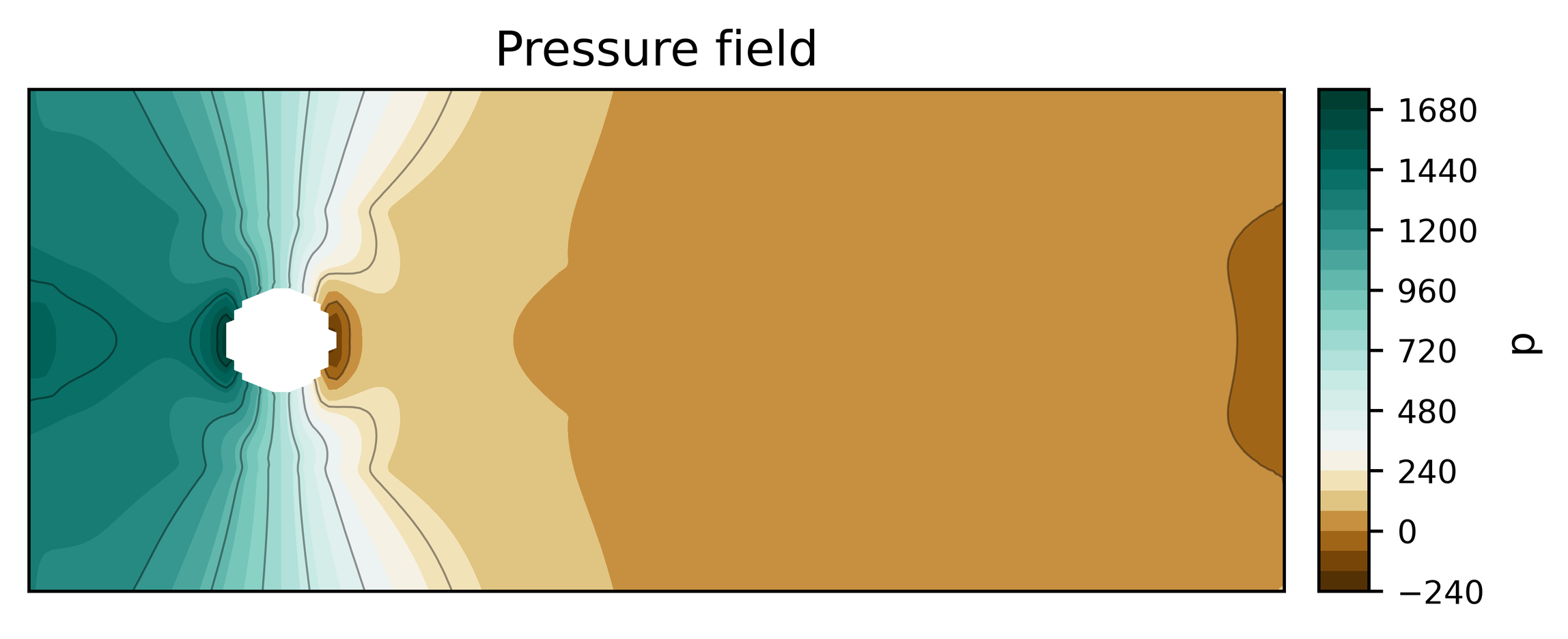}
        \caption{Pressure field $p$.}
    \end{subfigure}
    \hfill
    \begin{subfigure}{0.48\linewidth}
        \centering
        \includegraphics[width=\linewidth]{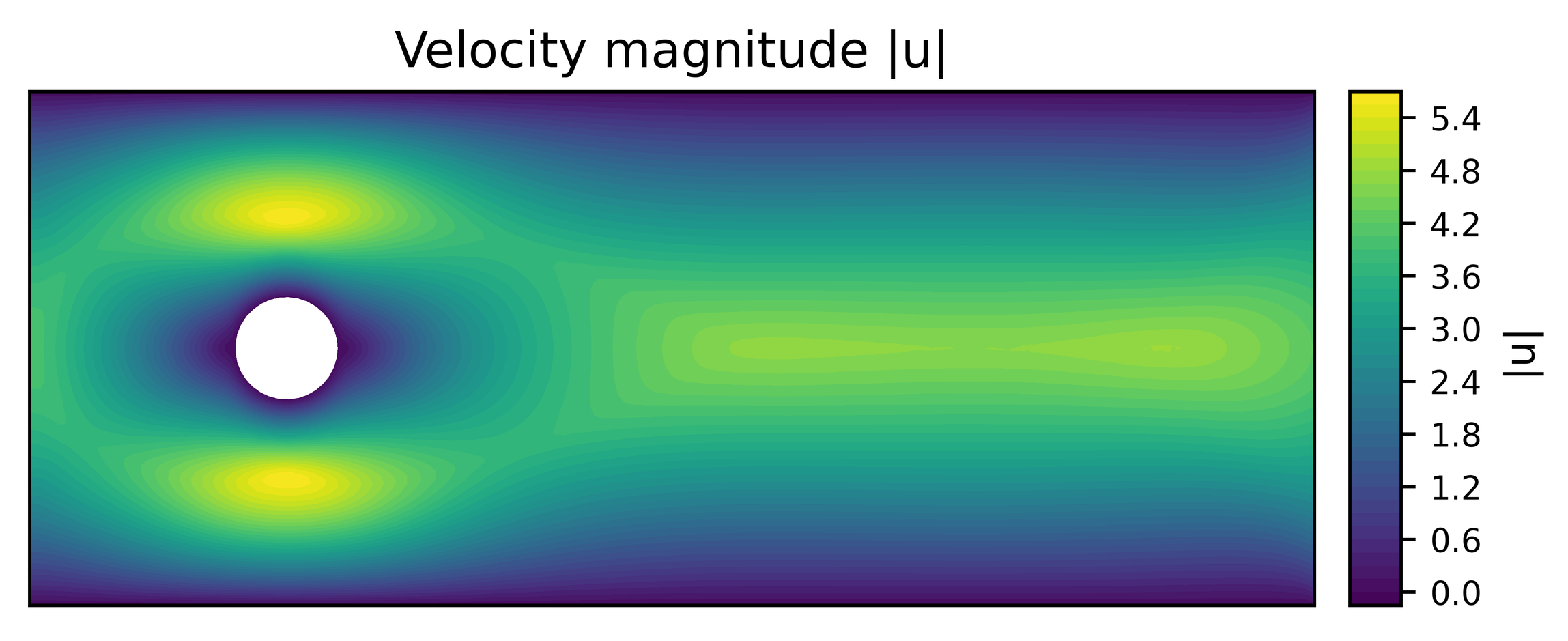}
        \caption{Velocity field magnitude $|\bu|$.}
    \end{subfigure}
    \caption{Pressure field and velocity magnitude for the circular obstacle benchmark problem, $r=4$, $n_h=128$.}
\end{figure}

\begin{figure}[htbp]
    \centering
    \includegraphics[width=0.35\linewidth,trim=0cm 0cm 0cm 0cm, clip]{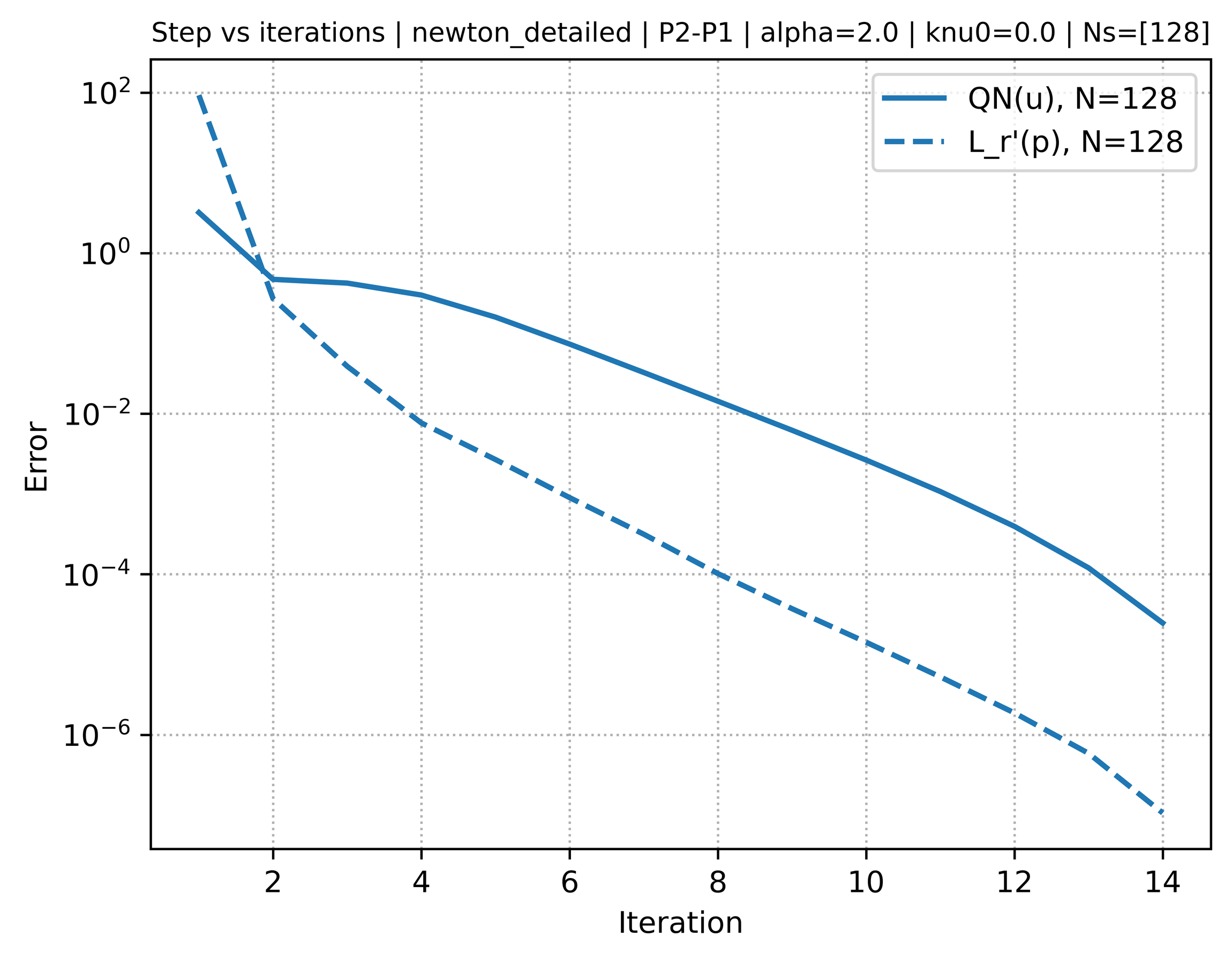}
    \caption{Increment criterion \eqref{eq:increment-criteria} vs iterations for the circular obstacle benchmark problem, $r=4$, $n_h=128$.}
\end{figure}

As $r\to 1$, the velocity profile becomes plug-like: almost constant across the central part of the flow, with significant shear occurring mainly near the walls.

\begin{figure}[htbp]
    \centering
    \begin{subfigure}{0.48\linewidth}
        \centering
        \includegraphics[width=\linewidth]{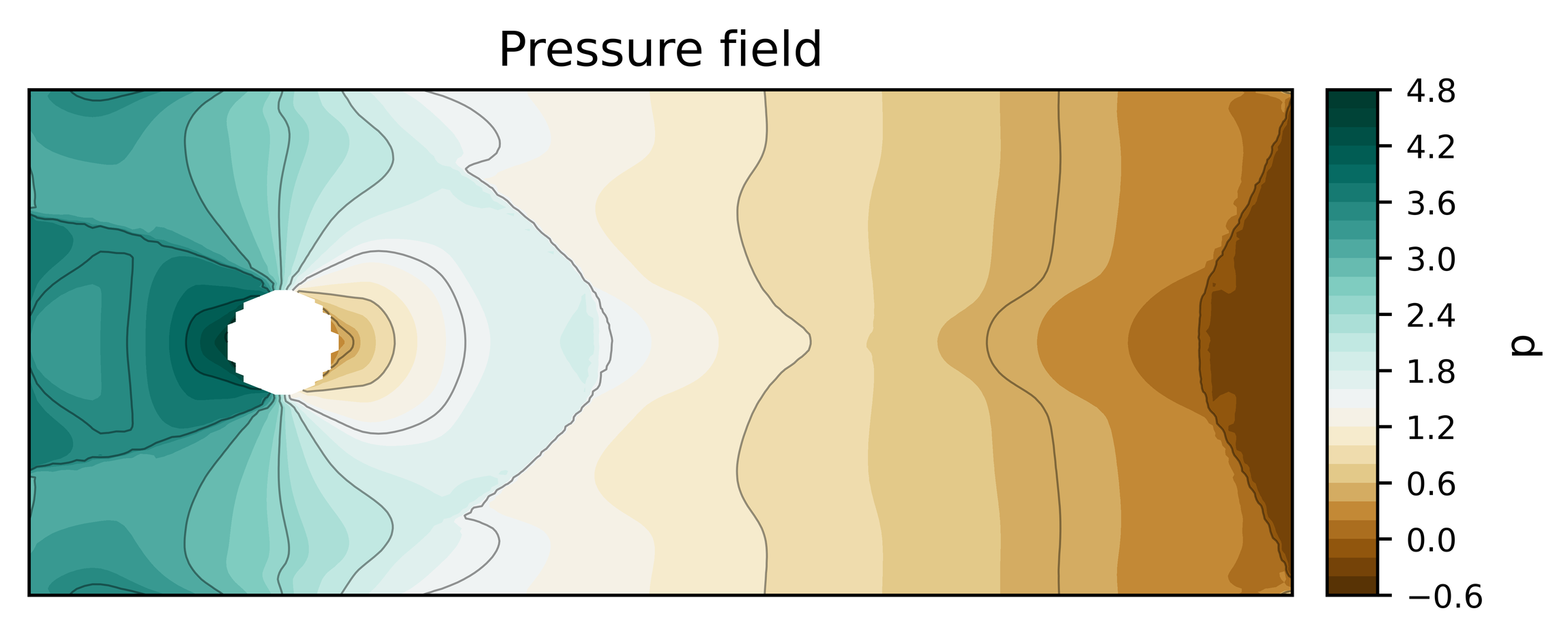}
        \caption{Pressure field $p$.}
    \end{subfigure}
    \hfill
    \begin{subfigure}{0.48\linewidth}
        \centering
        \includegraphics[width=\linewidth]{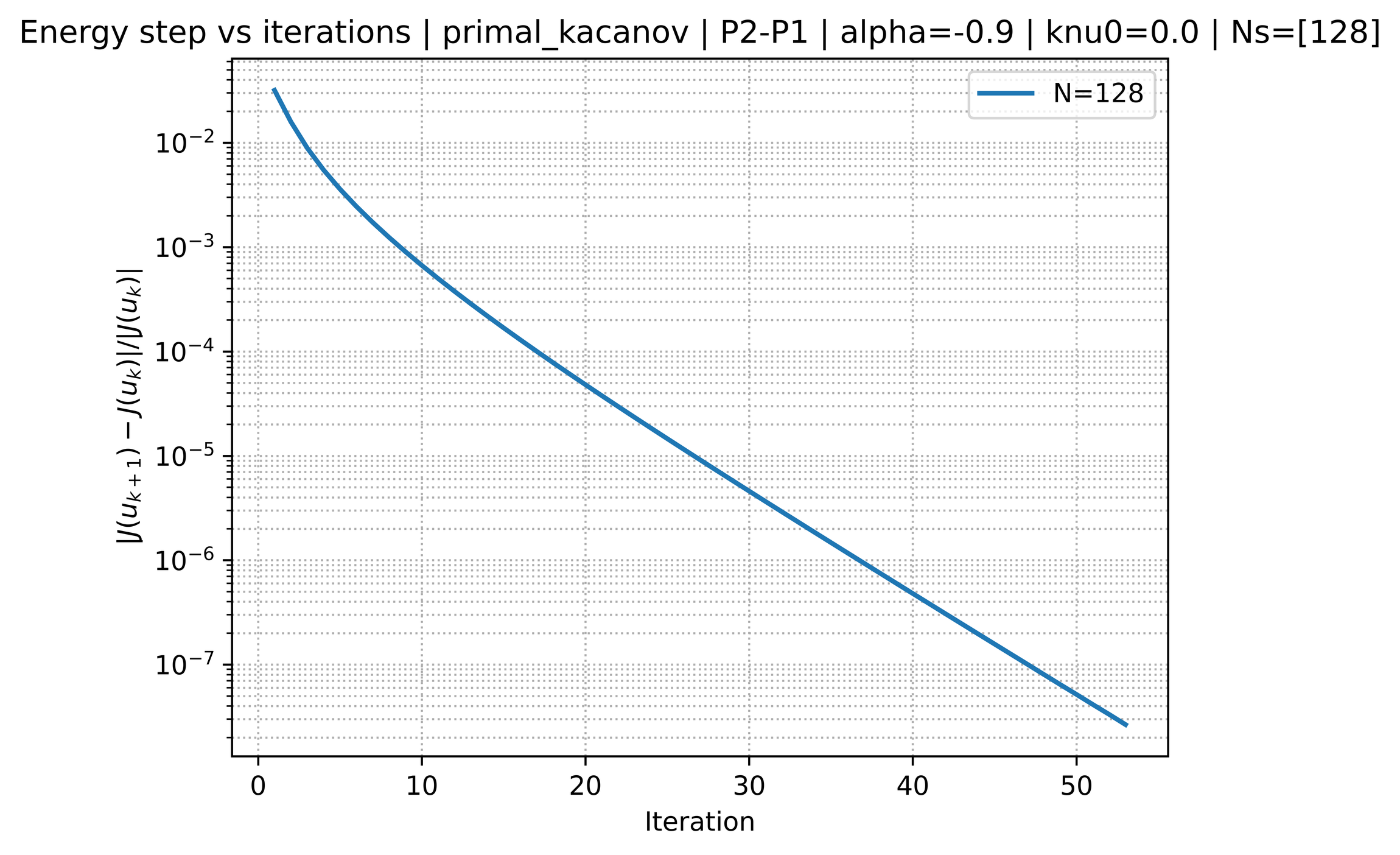}
        \caption{Energy step.}
    \end{subfigure}
    \caption{Pressure field and energy step for the circular obstacle benchmark problem, $r=1.1$, $n_h=128$, $k=0$.}
\end{figure}

\begin{figure}[htbp]
    \centering
    \begin{subfigure}{0.48\linewidth}
        \centering
        \includegraphics[width=\linewidth]{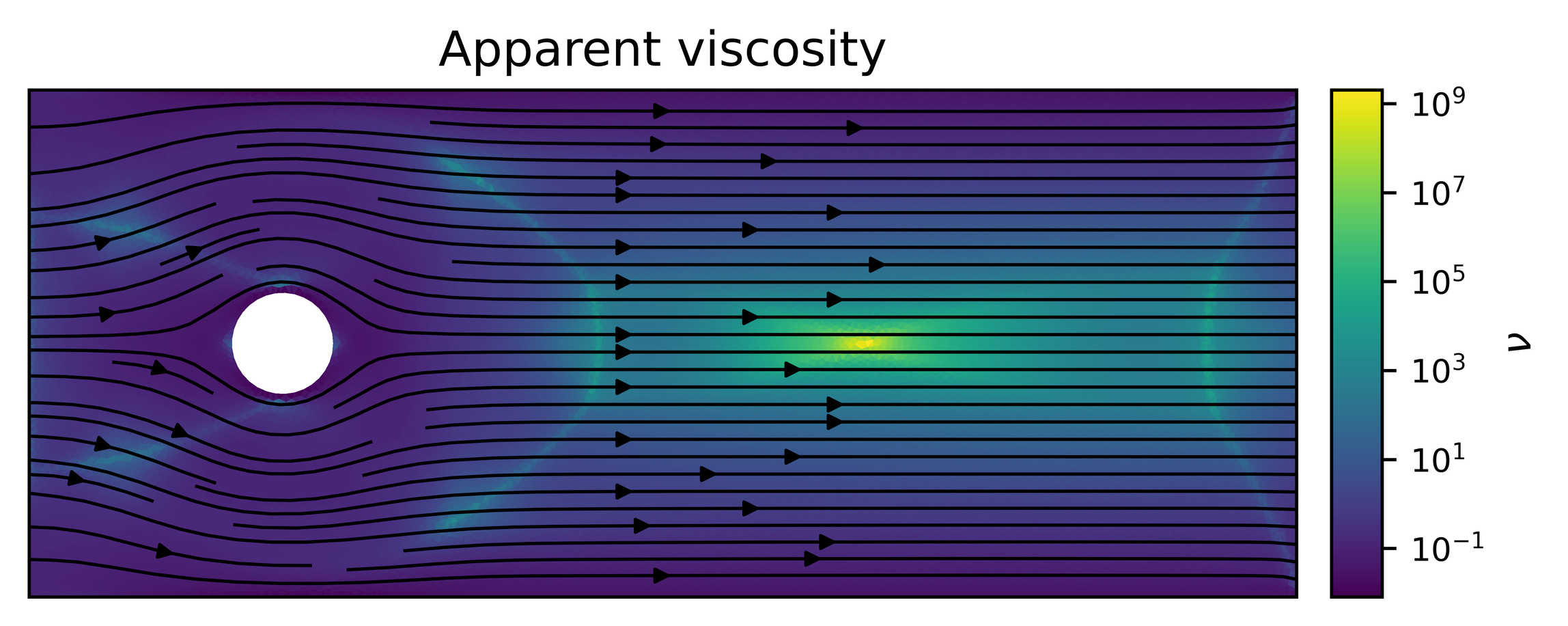}
        \caption{$k=0$.}
    \end{subfigure}
    \hfill
    \begin{subfigure}{0.48\linewidth}
        \centering
        \includegraphics[width=\linewidth]{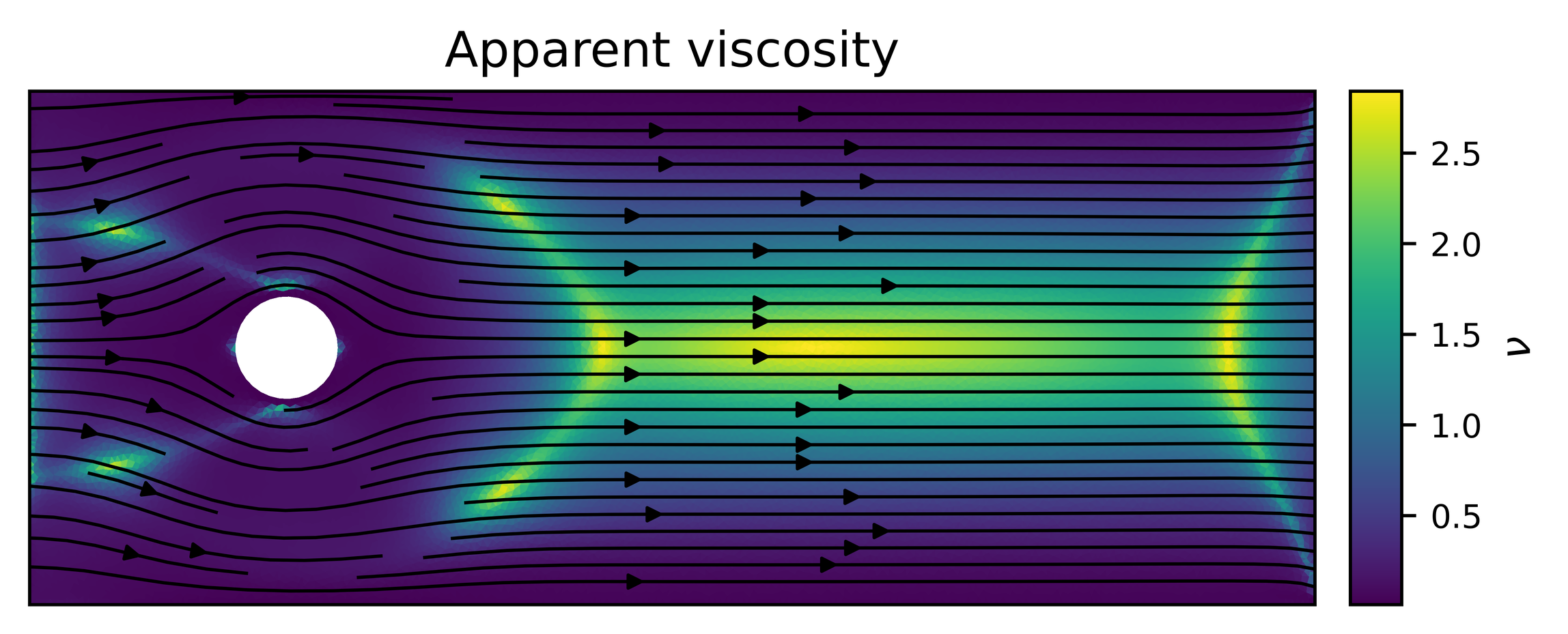}
        \caption{$k=100$.}
    \end{subfigure}
    \caption{Apparent viscosity \(\nu(\beps(\bu))=\mu\left(\delta_h+\|\beps(\bu)\|\right)^{r-2}\), with velocity streamlines, for $k=0,100$, $r=1.1$, $n_h=128$. Note the difference in the scales between the two plots.}
    \label{fig:apparent-viscosity-circular}
\end{figure}

\begin{figure}[htbp]
    \centering
    \begin{subfigure}{0.48\linewidth}
        \centering
        \includegraphics[width=\linewidth]{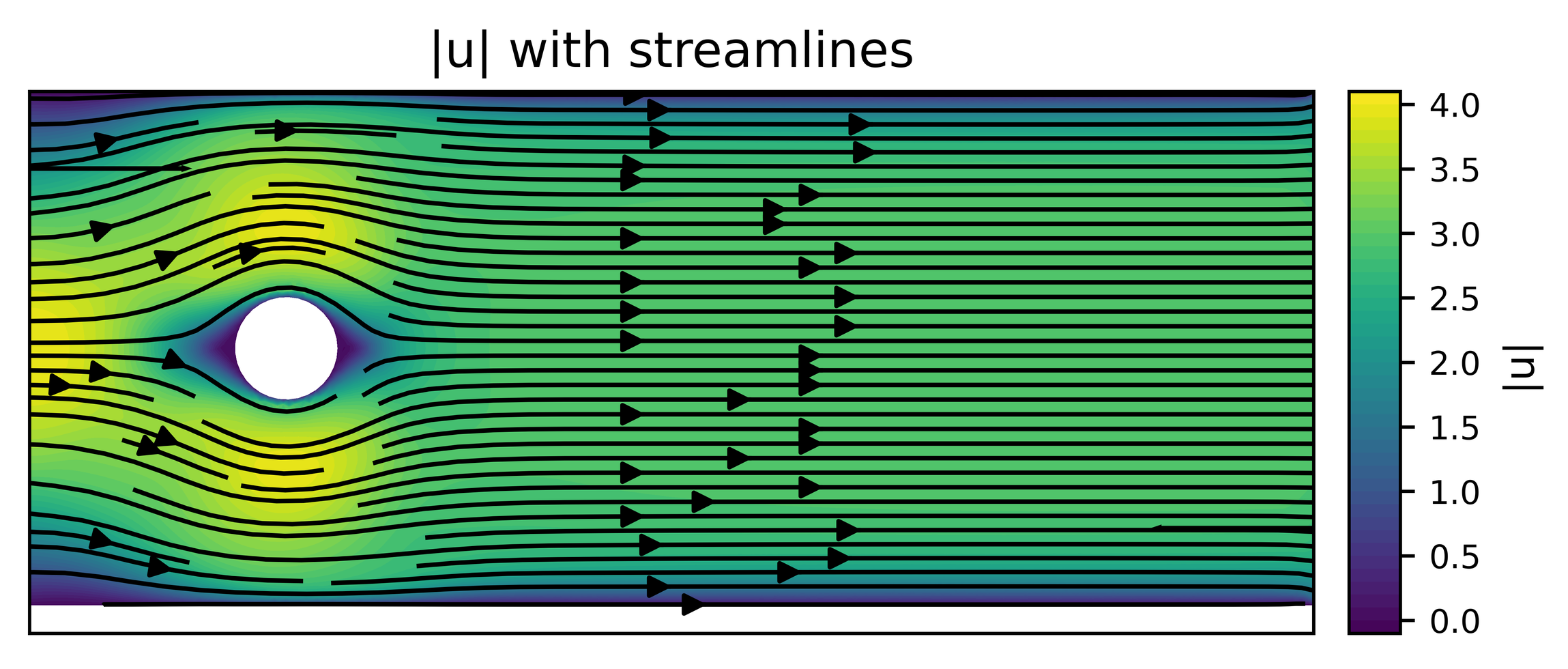}
        \caption{$k=0$.}
    \end{subfigure}
    \hfill
    \begin{subfigure}{0.48\linewidth}
        \centering
        \includegraphics[width=\linewidth]{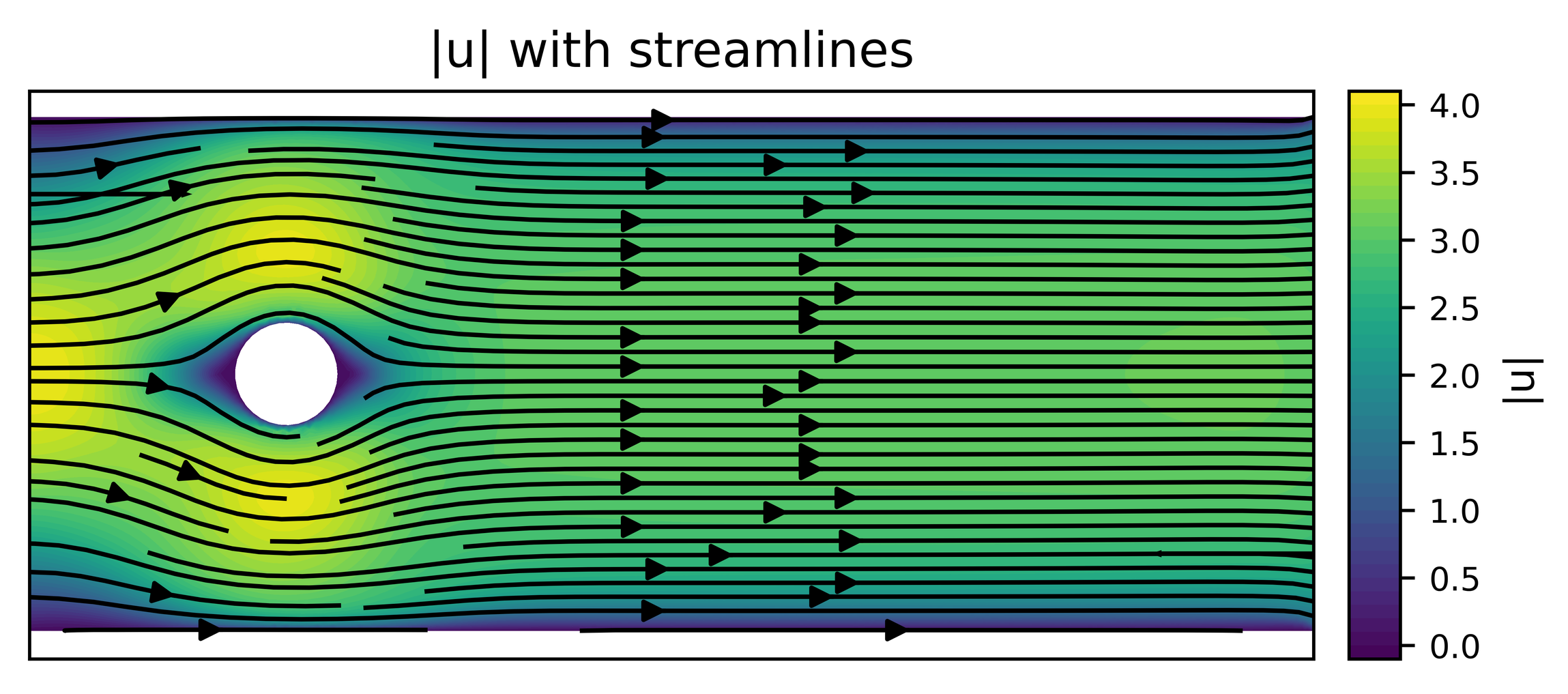}
        \caption{$k=100$.}
    \end{subfigure}
    \caption{Velocity streamlines for the corresponding values of $k$ in Figure~\ref{fig:apparent-viscosity-circular}, $r=1.1$, $n_h=128$.}
\end{figure}

\FloatBarrier
\subsection{Power-law and Carreau--Yasuda comparison}
\label{app:obstacle-blood-comparison}
The obstacle geometry was also used to investigate the effect of the
constitutive model in a hemodynamically motivated setting. The power-law and
Carreau--Yasuda parameters were taken from Alves et
al.~\cite{alves2013study}, using fits obtained from the same experimental blood
flow curve at $37^\circ\mathrm{C}$ with a parallel-plate gap of
$1000\,\mu\mathrm{m}$. This avoids comparing parameter sets obtained from
different experiments.

\begin{figure}[tbp]
    \centering

    \begin{subfigure}{0.48\linewidth}
        \centering
        \includegraphics[width=\linewidth]
        {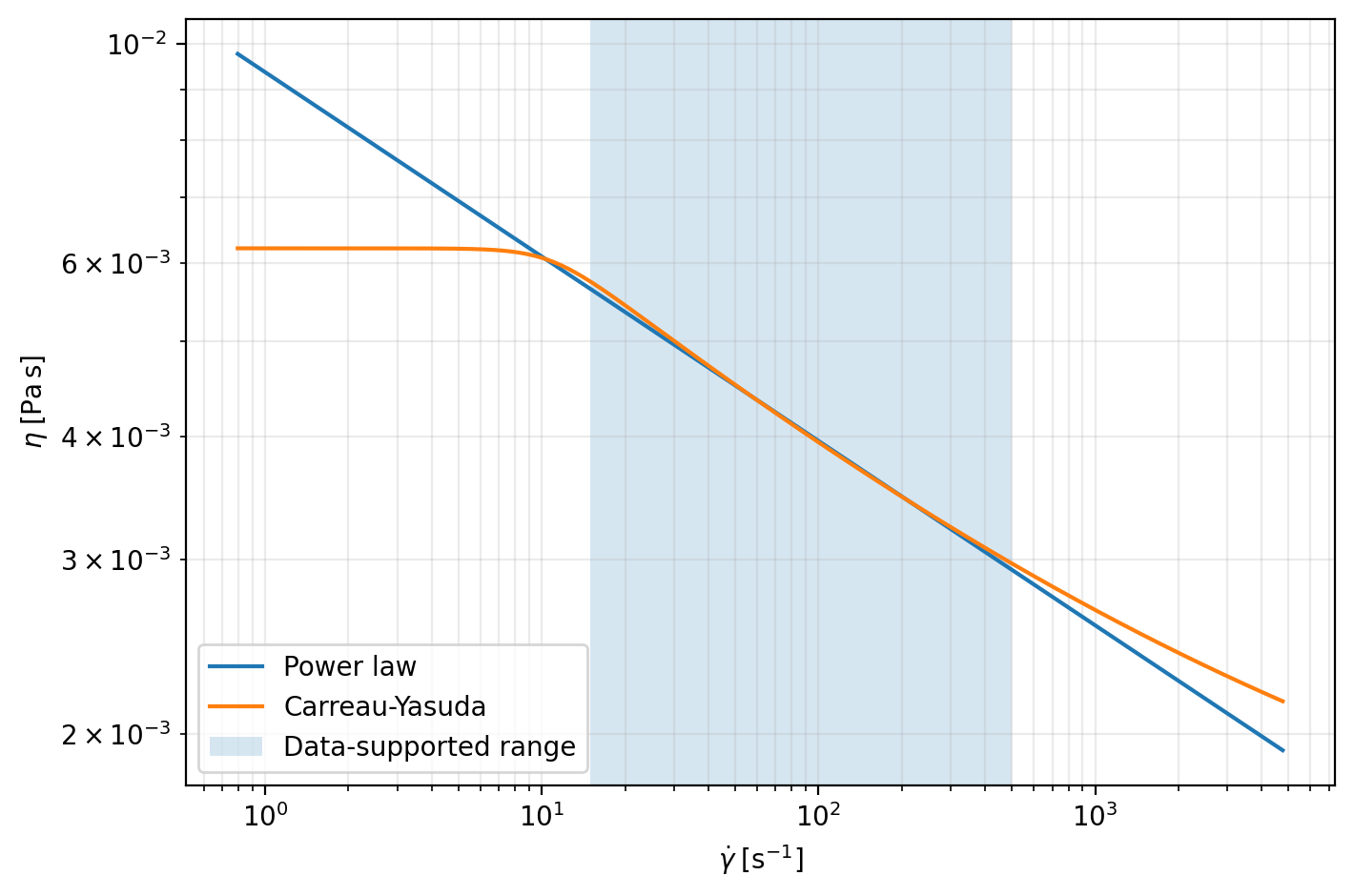}
        \caption{Apparent-viscosity curves.}
    \end{subfigure}
    \hfill
    \begin{subfigure}{0.48\linewidth}
        \centering
        \includegraphics[width=\linewidth]
        {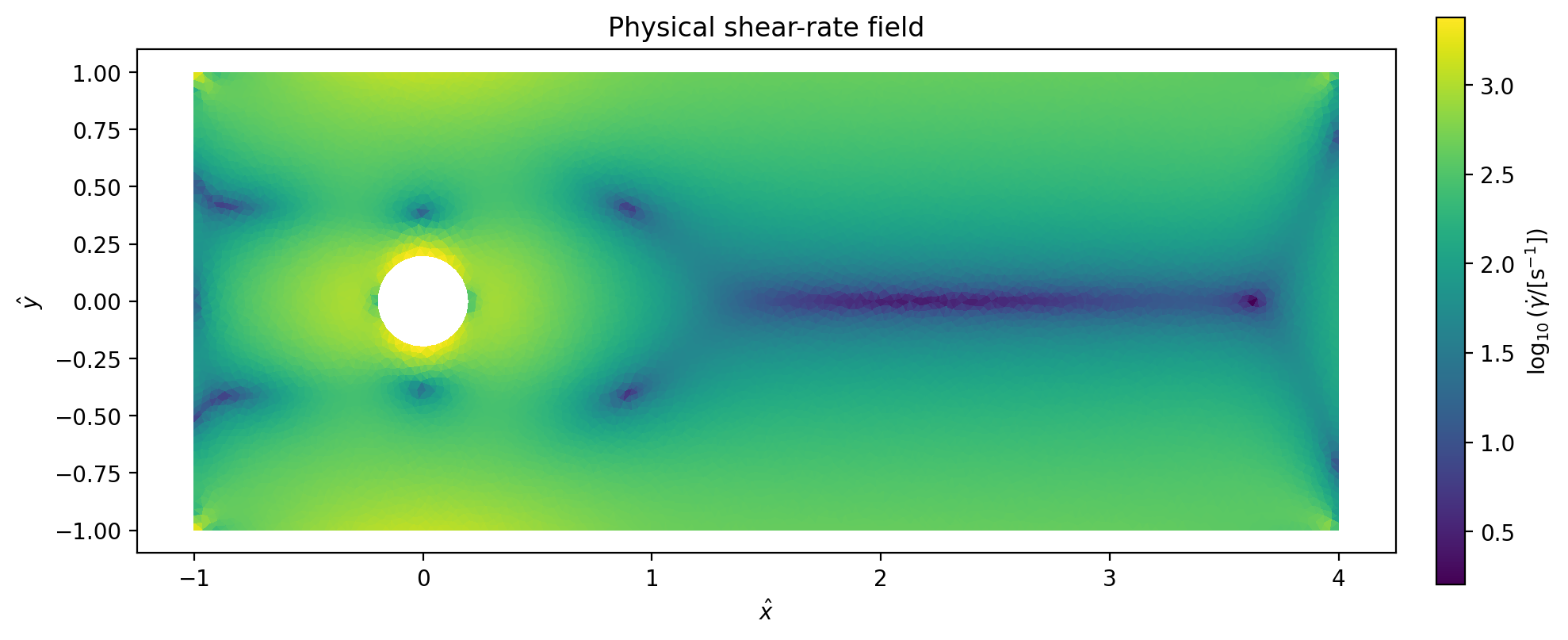}
        \caption{Physical shear-rate field.}
    \end{subfigure}

    \medskip

    \begin{subfigure}{0.48\linewidth}
        \centering
        \includegraphics[width=\linewidth]
        {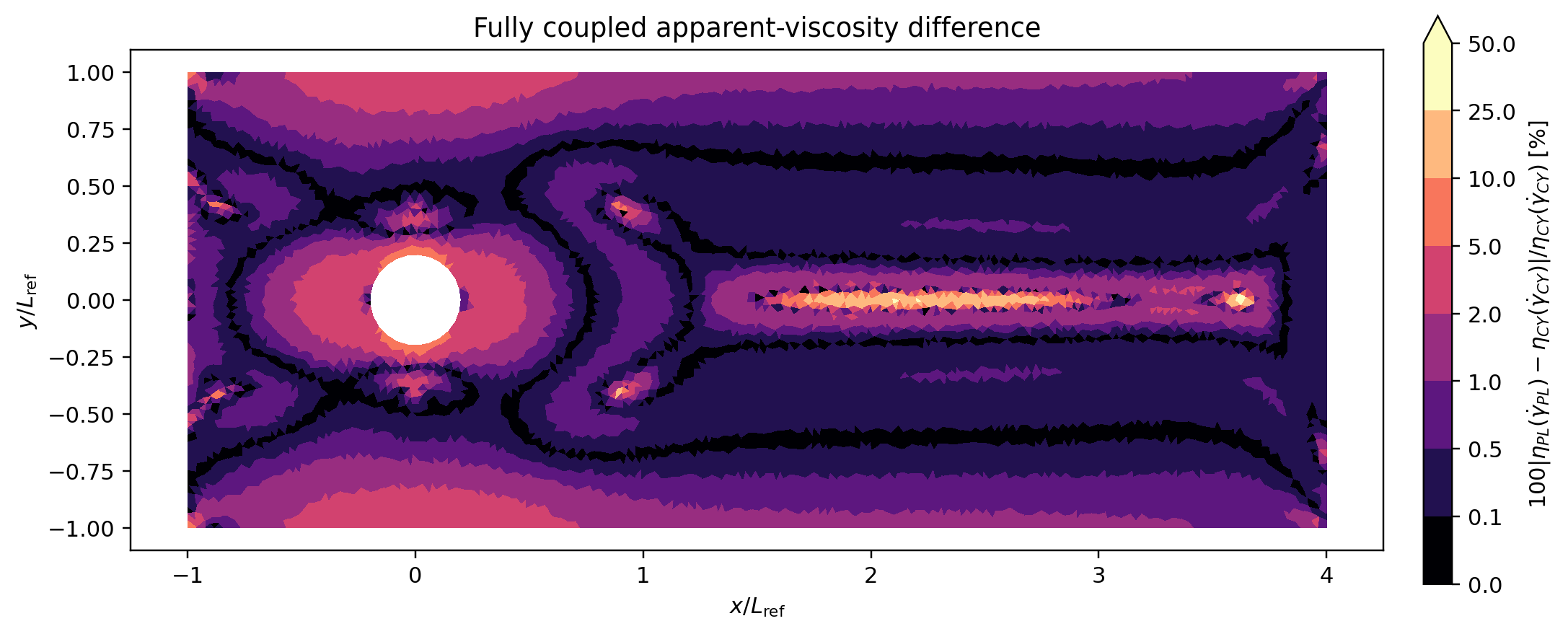}
        \caption{Fully coupled viscosity difference $e_\eta$.}
    \end{subfigure}
    \hfill
    \begin{subfigure}{0.48\linewidth}
        \centering
        \includegraphics[width=\linewidth]
        {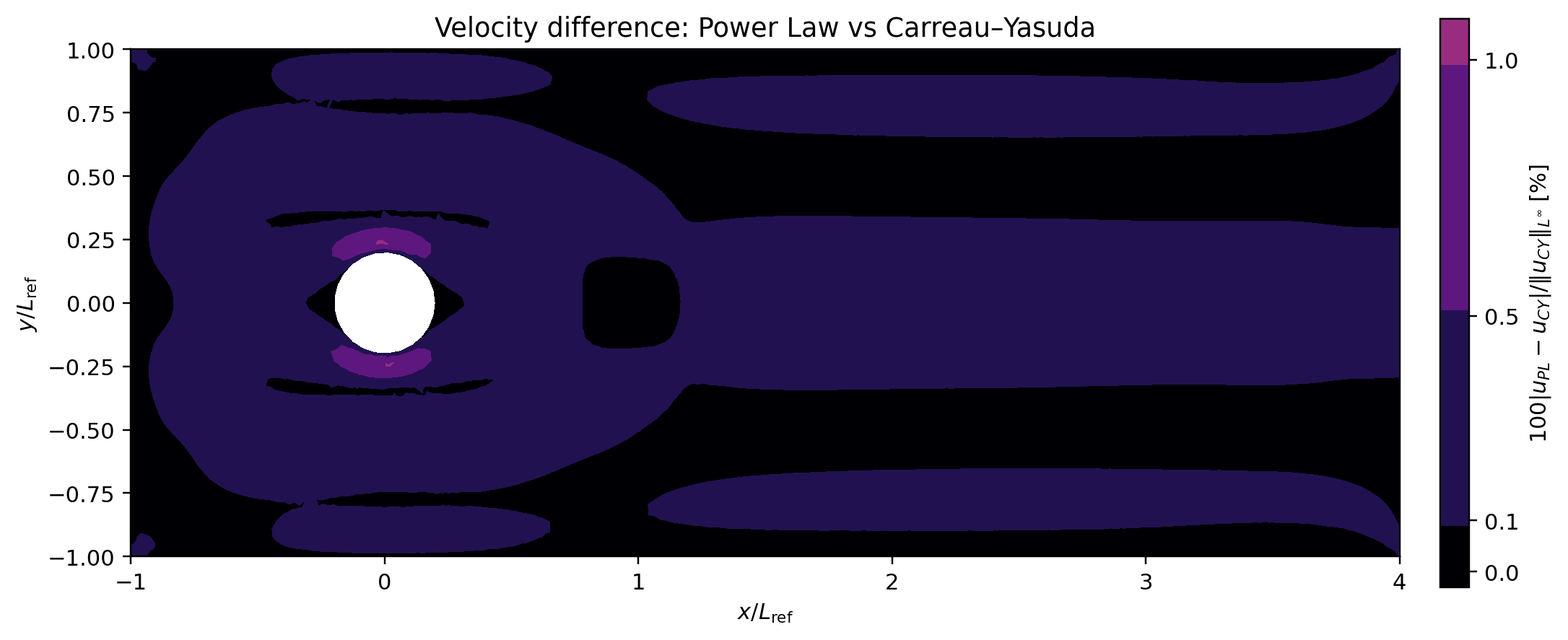}
        \caption{Normalized velocity difference $e_{\bu}$.}
    \end{subfigure}

    \caption{
    Power-law and Carreau--Yasuda comparison on the obstacle domain. The two
    models use parameters fitted to the same experimental flow curve reported
    in~\cite{alves2013study}.
    }
    \label{fig:blood-constitutive-fields}
\end{figure}

The selected power-law parameters are
\[
K=9.37\times10^{-3}\,\mathrm{Pa\,s}^{n_{\mathrm{PL}}},
\qquad
n_{\mathrm{PL}}=0.813,
\]
corresponding to $r=1.813$.
The Carreau--Yasuda parameters in~\eqref{eq:carreau-yasuda} are
\[
\eta_0=6.21\times10^{-3}\,\mathrm{Pa\,s},
\quad
\eta_\infty=1.17\times10^{-3}\,\mathrm{Pa\,s},
\quad
\lambda=9.23\times10^{-2}\,\mathrm{s},
\quad
n_{\mathrm{CY}}=0.732,
\quad
a=5.
\]
The reference scales
\(
L_{\mathrm{ref}}=10^{-3}\,\mathrm{m},
\qquad
U_{\mathrm{ref}}=5\times10^{-2}\,\mathrm{m\,s^{-1}}
\)
give the physical shear rate
\[
\dot\gamma_M
=
\frac{U_{\mathrm{ref}}}{L_{\mathrm{ref}}}
\sqrt{2}\|\beps(\bu_M)\|,
\qquad
M\in\{\mathrm{PL},\mathrm{CY}\},
\]
where \(\bu_M\) is the nondimensional velocity for model \(M\).
Figure~\ref{fig:blood-constitutive-fields} shows the constitutive curves and
the main results of the fully coupled comparison. The apparent-viscosity
difference is evaluated as
\[
e_\eta
=
100
\frac{
\left|
\eta_{\mathrm{PL}}(\dot\gamma_{\mathrm{PL}})
-
\eta_{\mathrm{CY}}(\dot\gamma_{\mathrm{CY}})
\right|
}{
\eta_{\mathrm{CY}}(\dot\gamma_{\mathrm{CY}})
},
\]
whereas the velocity difference is normalized by the global maximum
Carreau--Yasuda velocity,
\[
e_{\bu}
=
100
\frac{
|\bu_{\mathrm{PL}}-\bu_{\mathrm{CY}}|
}{
\|\bu_{\mathrm{CY}}\|_{L^\infty(\Omega)}
}.
\]

Although the two models produce appreciable local differences in the apparent viscosity, their velocity fields remain remarkably close, with discrepancies
of approximately one percent or less over almost the entire domain. Larger
differences occur in quantities depending on the velocity gradient
and constitutive response, particularly around the obstacle and in the
downstream wake.
\begin{figure}[tbp]
    \centering

    \begin{subfigure}{0.48\linewidth}
        \centering
        \includegraphics[width=\linewidth]
        {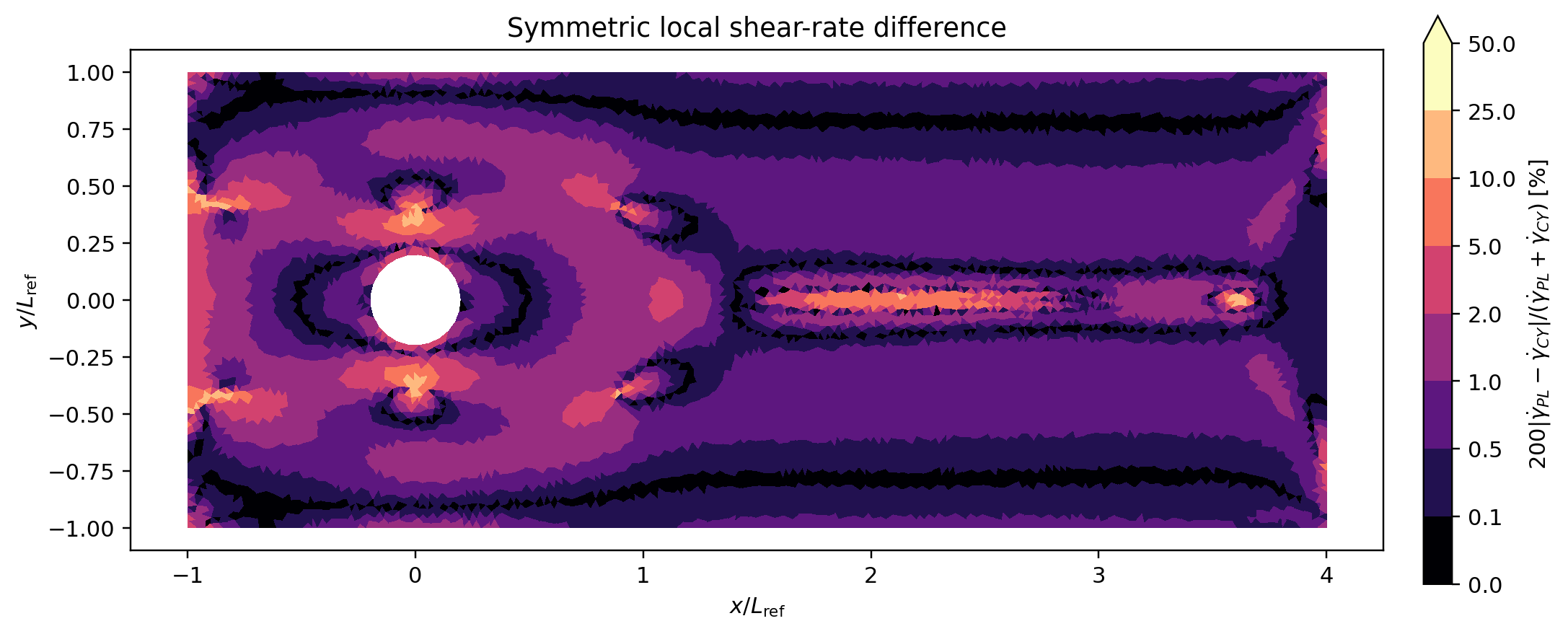}
        \caption{Symmetric shear-rate difference.}
    \end{subfigure}
    \hfill
    \begin{subfigure}{0.48\linewidth}
        \centering
        \includegraphics[width=\linewidth]
        {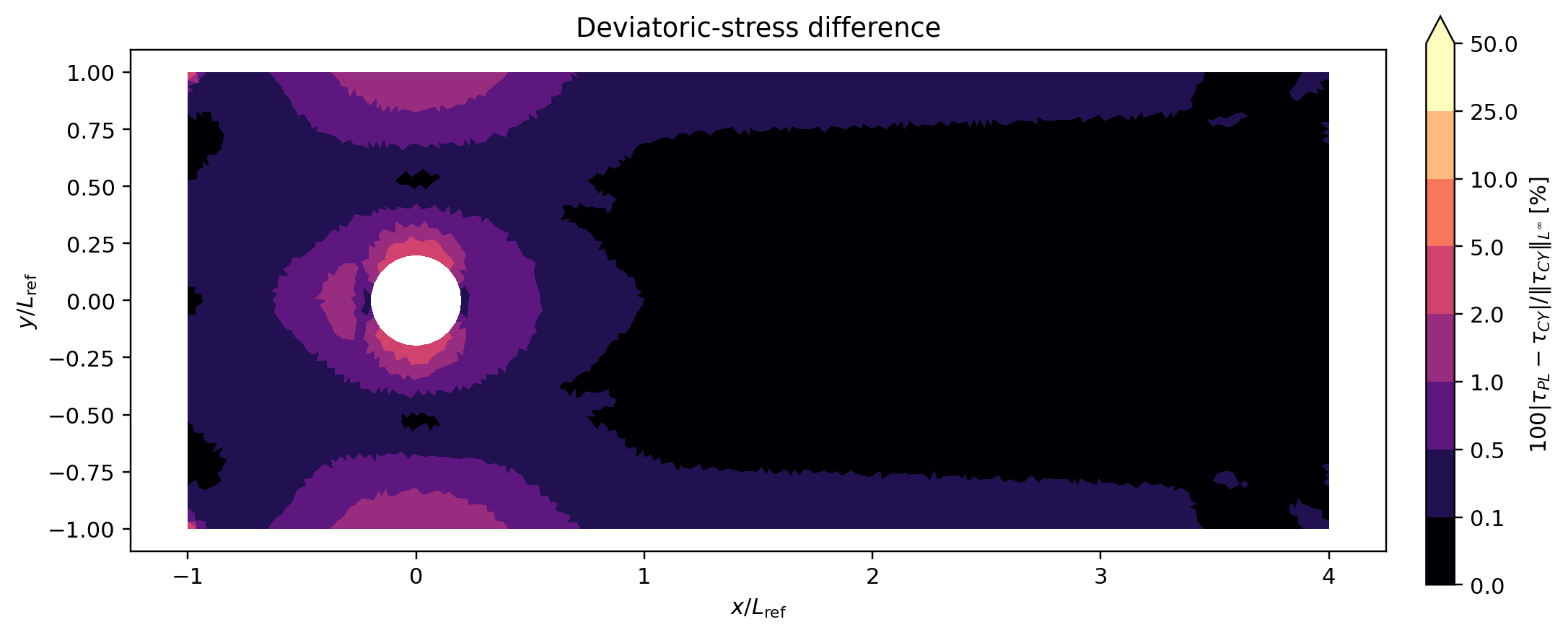}
        \caption{Deviatoric-stress difference.}
    \end{subfigure}

    \medskip

    \begin{subfigure}{0.48\linewidth}
        \centering
        \includegraphics[width=\linewidth]
        {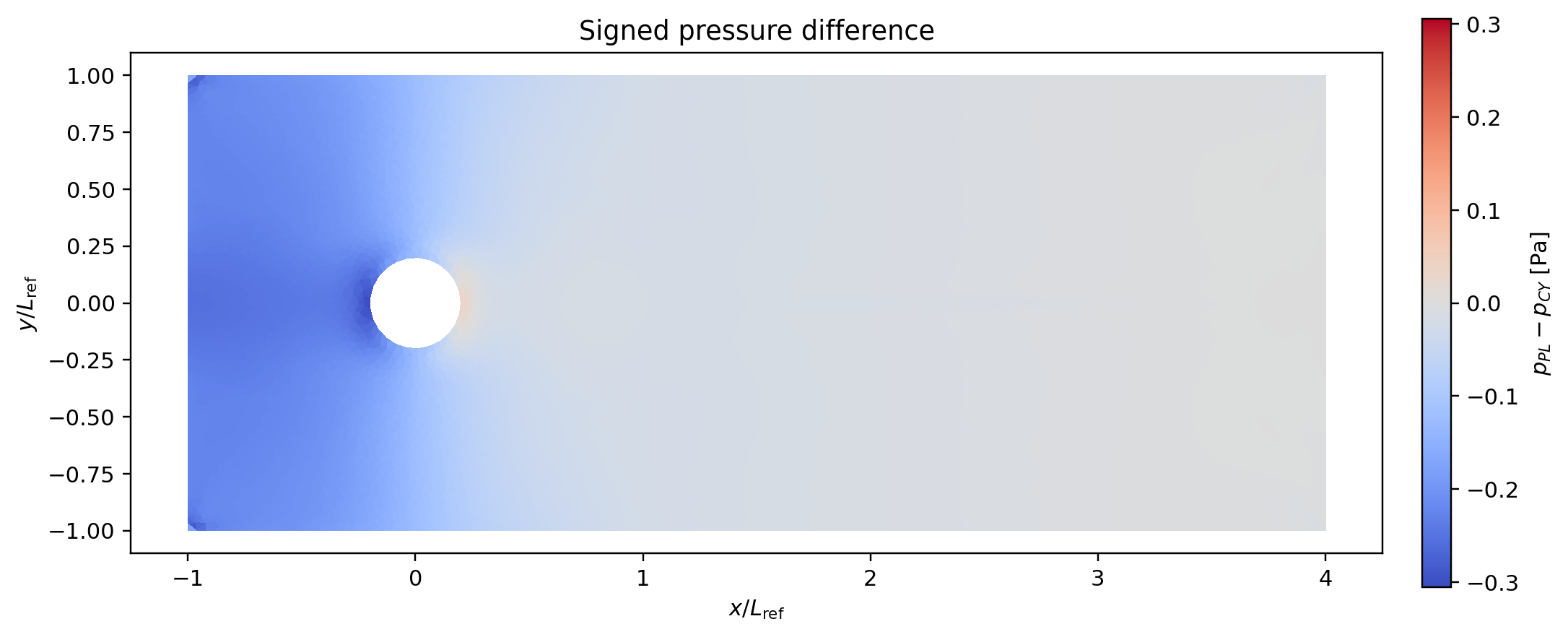}
        \caption{Signed pressure difference
        $p_{\mathrm{PL}}-p_{\mathrm{CY}}$.}
    \end{subfigure}
    \hfill
    \begin{subfigure}{0.48\linewidth}
        \centering
        \includegraphics[width=\linewidth]
        {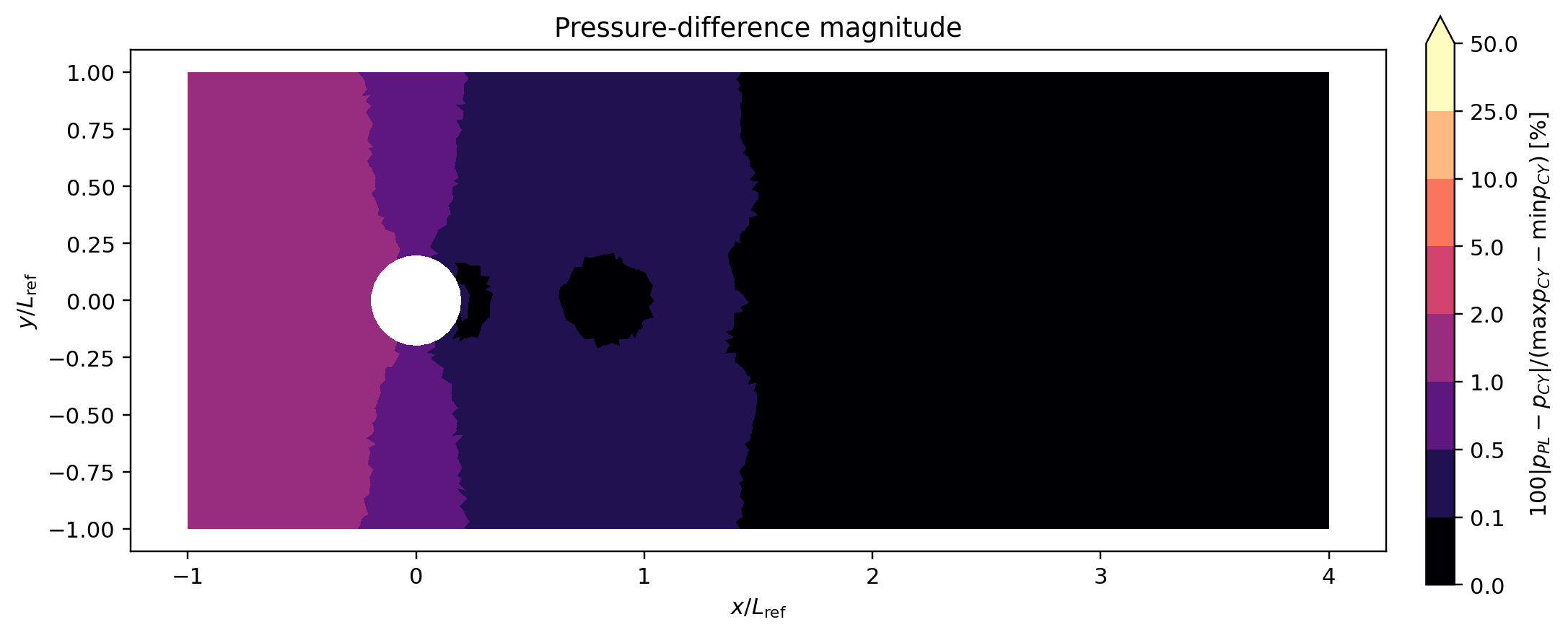}
        \caption{Normalized pressure difference.}
    \end{subfigure}

    \caption{
    Differences in shear rate, deviatoric stress, and pressure between the
    independently computed power-law and Carreau--Yasuda solutions.
    }
\end{figure}
The results indicate that the power law provides a close approximation of the
Carreau--Yasuda velocity and pressure fields for this benchmark. The comparison therefore suggests
the use of the power law for global kinematic predictions.

\end{document}